\RequirePackage{fix-cm}
\documentclass[12pt]{article}

\usepackage[T1]{fontenc}    
\usepackage{hyperref}       
\usepackage{url}            
\usepackage{booktabs}       
\usepackage{amsfonts}       
\usepackage{nicefrac}       
\usepackage[margin = 1in]{geometry}

\usepackage{amsmath, amssymb, amsthm}
\usepackage{cancel, makecell, graphicx, multirow, rotating, enumerate}
\usepackage{diagbox}
\usepackage{enumitem} 
\usepackage{xcolor}
\usepackage[final]{changes}
\usepackage{wrapfig}
\colorlet{Changes@Color}{green!60!black}

\usepackage[linesnumbered, ruled]{algorithm2e}

\SetKwRepeat{Do}{do}{while}
\SetKw{KwGoTo}{go to}
\SetKwProg{Fn}{Function}{}{}
\SetKw{KwBreak}{break}

\newtheorem{theorem}{Theorem}
\newtheorem{defn}{Definition}
\newtheorem{lem}{Lemma}
\newtheorem{rem}{Remark}
\newtheorem{propn}{Proposition}
\newtheorem{cor}{Corollary}

\usepackage{footnote}
\makesavenoteenv{tabular}
\makesavenoteenv{table}

\begin{document}
\title{Optimal Posteriors for Chi-squared Divergence based PAC-Bayesian Bounds and Comparison with KL-divergence based Optimal Posteriors and Cross-Validation Procedure}

\author{Puja Sahu \and Nandyala Hemachandra}
\date{}

\maketitle

\begin{abstract} 
We investigate optimal posteriors for recently introduced \cite{begin2016pac} chi-squared divergence based PAC-Bayesian bounds in terms of nature of their distribution, scalability of computations, and test set performance. For a finite classifier set, we deduce bounds for three distance functions: KL-divergence, linear and squared distances. Optimal posterior weights are proportional to deviations of empirical risks, usually with subset support. For uniform prior, it is sufficient to search among posteriors on classifier subsets ordered by these risks. We show the bound minimization for linear distance as a convex program and obtain a closed-form expression for its optimal posterior. Whereas that for squared distance is a quasi-convex program under a specific condition, and the one for KL-divergence is non-convex optimization (a difference of convex functions). To compute such optimal posteriors, we derive fast converging fixed point (FP) equations. We apply these approaches to a finite set of SVM regularization parameter values to yield stochastic SVMs with tight bounds. We perform a comprehensive performance comparison between our optimal posteriors and known KL-divergence based posteriors on a variety of UCI datasets with varying ranges and variances in risk values, etc. Chi-squared divergence based posteriors have weaker bounds and worse test errors, hinting at an underlying regularization by KL-divergence based posteriors. Our study highlights the impact of divergence function on the performance of PAC-Bayesian classifiers. We compare our stochastic classifiers with cross-validation based deterministic classifier. The latter has better test errors, but ours is more sample robust, has quantifiable generalization guarantees, and is computationally much faster.
\end{abstract}


\section{Introduction}
\label{intro}
In classification algorithms, the choice of the parameter(s) influences the level of accuracy that the generated classifier can achieve. For example, consider the Support Vector Machine (SVM) algorithm for classification with the regularization parameter, $\lambda > 0$. This parameter is a user input which trades off between model complexity and training error. The optimal classifier that we get, depends heavily on the sample $S$ that is used for training and the value of the parameter, $\lambda$. We can control only this parameter value for obtaining a classifier with low (training) error, \textit{but not} the given data. For a given training sample, we can choose the best value of the parameter from a prefixed set of values, which yields a classifier with the lowest error. However, this is a long drawn process. Plus, there is no guarantee that the chosen value will yield a classifier having low(est) error on another sample from the same distribution. This implies that the best parameter value is sample dependent and that there is no unique value which is best for almost all the training samples (Please see \cite{pmlr-v101-sahu19a} and Appendix A, pg. 17 in \cite{sahu2019PACBKLarxiv} for illustration on a UCI dataset). However, if we determine the set of $\lambda$ values with, say, lowest 30\% error rates on each sample, we observe a recurring subset of $\lambda$ values across these samples (Please see \cite{pmlr-v101-sahu19a} and Table 4, in Appendix B in \cite{sahu2019PACBKLarxiv} for illustration on a UCI dataset). Thus, we have an ensemble of $\lambda$ values to pick from. We can combine multiple base classifiers resulting from different parameter values, to build a strong stochastic classifier using PAC-Bayesian framework.

\textbf{PAC-Bayesian Bounds and Optimal Posteriors} PAC-Bayesian approach assumes an arbitrary but fixed prior distribution on the space of classifiers and outputs a posterior distribution on this space, corresponding to a stochastic classifier. This approach provides a probabilistic bound on the difference between the posterior averaged true  and empirical risks of a stochastic classifier as measured by a convex distance function. These bounds on \textit{unknown averaged true risk} offer a trade-off between averaged empirical risk and a term which encompasses model complexity of the stochastic classifier. The bound is computed based on a single sample but with a high probability guarantee over different samples (from the same distribution). For a chosen distance function, we are interested in the `optimal PAC-Bayesian posterior' defined as the posterior distribution which \textit{minimizes the corresponding PAC-Bayesian bound}. By design, these bounds and the resulting optimal posterior are robust to the choice of sample used for training, addressing the sample bias.

\textbf{Relevant Work} A well known form of bounds estimating the unknown true risk of a classifier known as PAC-Bayesian bounds were proposed by \cite{mcallester1998PACBthms,mcallester1999PACBmodelavg,Seeger02theproof} using the idea of Bayesian priors and posteriors on the classifier space, and refined further by \cite{maurer2004note,langford2005tut,mcallester2013pac,McAllesterAkinbiyi2013Festschrift}. Several authors improvised the bounds for the choice of the distance function they considered for evaluating the classifiers. While \cite{maurer2004note} provided a bound for the KL divergence as the distance function, $\phi$, by tightening up the threshold with a factor of $\sqrt{m}$ instead of $m$, \cite{germain2009pac} generalized the framework of PAC-Bayesian bounds for a broader class of convex $\phi$ functions and relaxed the constraints on tail bounds of the empirical risk of the classifiers under consideration. PAC-Bayesian theory has been used to devise margin bounds for linear classifiers such as SVMs \cite{mcallester2003pac,LangfordTaylor2002}. \cite{AmbroladzeNIPS2006} specialized the PAC-Bayesian bounds using spherical Gaussian distributions on the space of linear classifiers and extended the set up for using data-dependent priors \cite{Parrado-HernandezJMLR2012}. More recently, \cite{begin2016pac} introduced PAC-Bayesian bounds based on  R{\'e}nyi divergence between the posterior and the prior distributions. We use a specific case of this R{\'e}nyi divergence which corresponds to $\chi^2$-divergence.

All of the above consider a continuous (SVM) classifier space ($n$-dimensional Euclidean space) and continuous prior as well as posterior distributions on it (spherical Gaussian distributions) whereas we consider a finite set of classifiers such as those generated by a finite set of regularization parameter values for the SVM. Our $\chi^2$-divergence based PAC-Bayesian bounds are derived for this set up with a discrete prior distribution, and three different distance functions between posterior averaged empirical risk and posterior averaged true risk. \added{The motivation for choosing a different divergence function in the PAC-Bayesian framework is to achieve a better test set performance and a tighter risk bound. In order to do so, we \textit{first} need to investigate the nature of these $\chi^2$-squared divergence based PAC-Bayesian bound minimization problems, identify the corresponding optimal PAC-Bayesian posteriors and understand their nature. These posteriors might not be at par with the classical KL-divergence based PAC-Bayesian posteriors or the cross-validation based procedure, but the comparison amongst them brings forth some insightful aspects of the PAC-Bayesian optimal posteriors. We list below the contributions of this paper.}

\textbf{Contributions} We are interested in the \textit{optimal PAC-Bayesian posterior} which minimizes the $\chi^2$-divergence based PAC-Bayesian bound for a given distance function (Section \ref{secn:genericPACBbnd}).  We consider a finite classifier set and three distance functions -- linear distance, squared distance (second degree polynomial) and KL-divergence (infinite degree polynomial). 

\begin{itemize}[label = \textbullet]
    	\item We deduce $\chi^2$-divergence based PAC-Bayesian bounds for the above three distance functions and identify the optimal posteriors for them via respective bound minimization problems. 
	\item The linear distance based bound was considered in \cite{begin2016pac}; we identify the associated bound minimization as a convex program and obtain a closed form expression for the global optimal posterior (Section \ref{secn:linChi2}). 
	\item We also deduce PAC-Bayesian bounds for squared distance and KL-divergence, and show that they are non-convex programs (Sections \ref{secn:sqChi2} and \ref{secn:klChi2}). We further show that the squared distance based bound is quasi-convex under certain conditions. \added{The KL-divergence based bound minimization problem involves a difference of convex (DC) functions and hence is a DC program. Therefore we applied a DC approach known as Convex-Concave Procedure (CCP) \cite{CCP2016LippBoyd} to find its local minimum. In our computations, we observed that the CCP did not work for certain cases, especially when we have almost linearly separable data. (Such cases are illustrated in Table \ref{appdx_tab:Bnd.TestErr.FPnCCP.klChi2} in Appendix D.)}
	\item For deriving optimal posteriors for such non-convex cases, we identify Fixed Point (FP) equations deduced from the partial KKT system with strict positivity constraints. These FP equations converge even when the solver or an alternate approach like CCP fails to identify a solution, and are much faster than the solver. (Some examples of such cases are in Tables \ref{appdx_tab:Bnd.klChi2}, \ref{appdx_tab:TestErr.klChi2} and \ref{appdx_tab:Bnd.TestErr.FPnCCP.klChi2} in the appendix.)
    \item For any of the above 3 distance functions, for the uniform prior distribution, we simplify the search for optimal posteriors on the simplex restricted to subsets of classifiers ordered by empirical risk values (Section \ref{secn:OptQ_unifP}). 
	\item For computational illustration, we consider a comprehensive set of nine UCI datasets \cite{UCI:2017} with small to moderate number of examples and features, balanced and imbalanced classes, and having different ranges and variances in the empirical risk values. Using such datasets helps us compare and understand the performance of optimal posteriors due to different distance functions for the $\chi^2$-squared divergence based optimal posteriors and also compare with the known KL-divergence based optimal posteriors \cite{pmlr-v101-sahu19a}.
    \item We use these approaches on the set of SVMs generated by a finite set of regularization parameter values (Section \ref{secn:SVM_PACB}). This leads us to the notion of a \textit{stochastic SVM} characterized by an optimal posterior on the regularization parameter set. \added{Usually small values of the regularization parameter values are preferred. Keeping this in mind, we used an arithmetic-geometric series of regularization parameter values, $\lambda$ with a logarithmic scale for $\lambda \in (0,  0.1)$ and a linear scale for $\lambda \geq 0.1$. We chose a mixture of logarithmically and linearly spaced values of $\lambda$ so that we cover many different $\lambda$s corresponding to distinct SVM classifiers with low test errors.} 
    \begin{itemize}
        \item Optimal posteriors for KL-divergence give extremely loose bounds and are computationally expensive but have test error rates generally better than linear distance ones. The optimal bound value and test error rate of the squared distance based optimal posterior are remarkably lower than those of linear or KL-distance based posteriors when base classifiers have high variation in empirical risk values. 
        \item This is accompanied by relatively high concentration on low empirical risk values and sparse nature of squared distance based posteriors. For almost separable datasets, posteriors due to these three distance functions  have comparable PAC-Bayesian bound values and test error rates. 
    \end{itemize}
    Table \ref{tab:Chi2.resultoutline} outlines theoretical and computational aspects of optimal posteriors considered in this paper. 
    \item To understand the role of divergence measure on PAC-Bayesian bounds, we conducted a comparative study of these $\chi^2$-divergence based optimal PAC-Bayesian posteriors with the posteriors derived for classical KL-divergence based PAC-Bayesian bounds \cite{pmlr-v101-sahu19a} (Section \ref{secn:compareKLCV}).
\begin{itemize}
	\item We observe that though both the classes of posteriors have weights which are decreasing with the increasing empirical risk values of classifiers, the rate at which they decrease is different in the two classes -- KL-divergence based posteriors decrease exponentially, while $\chi^2$-divergence based posteriors decrease linearly with empirical risk values. 
	\item Another difference is in the size of support set for the two classes of posteriors -- KL-divergence based posteriors take up the full support on the set of base classifiers, whereas those for $\chi^2$-divergence usually depend only on a strict subset as their support set. 
	\item The class of optimal posteriors for $\chi^2$-divergence based PAC-Bayesian bounds is observed to have weaker bounds and higher test set errors than the class of KL-divergence based PAC-Bayesian posteriors on a set of SVM classifiers. Such behaviour can be attributed to $\chi^2$-divergence based posteriors \textit{overfitting} the data by choosing a strict subset support of classifiers with least empirical risk values. 
\end{itemize}
\item We also compared the performance of the stochastic SVM classifier governed by these PAC-Bayesian posteriors with the deterministic SVM classifier obtained via the popular cross-validation procedure for regularization parameter selection (Section \ref{secn:compareKLCV}) as the baseline case. Though the cross-validation procedure gives a classifier with better performance on a test set, the PAC-Bayesian posteriors yield a sample robust classifier with quantifiable guarantees on the unknown true risk. On the computation side, PAC-Bayesian procedure is more than 10 times faster than the cross-validation procedure.  
\end{itemize}


{ 
\begin{table}[htp]
\caption[An outline of theoretical aspects and computational results for the $\chi^2$-divergence PAC-Bayesian bound minimization problem]{\small An outline of theoretical aspects and computational results for optimal posteriors $Q^{\ast}_{\phi, \chi^2} =  \lbrace q^{\ast}_{i, \phi, \chi^2} \rbrace_{i =1}^H$ for minimization of PAC-Bayesian bound $B_{\phi, \chi^2}(Q)$ based on chi-squared divergence, $\chi^2(Q||P) = \sum_{i = 1}^H \frac{q_i^2}{p_i}$ between a posterior $Q$ and a prior $P$ on the classifier space $\mathcal{H}$. We consider three different distance functions, $\phi$:  KL-divergence $kl(\hat{l}, l) = \hat{l}\ln \frac{\hat{l}}{l} + (1 - \hat{l}) \ln \left( \frac{1 - \hat{l}}{1 - l} \right)$, linear $\phi_{\text{lin}} (\hat{l}, l) = l - \hat{l}$ and squared distances $\phi_{\text{sq}}(\hat{l}, l) = (l - \hat{l})^2$ for $l, \hat{l} \in (0,1)$. $H$ denotes the classifiers set size and $H^{\ast}$ denotes the size of the support set of the optimal posterior $Q^{\ast}_{\phi, \chi^2}$. $\hat{l}_i$ denotes empirical risk value of a classifier in $\mathcal{H}$ computed on a sample of size $m$. $\mathbf{\mathcal{I}^{R}_{\phi}} (m)$ is a sample size based constant for a distance function $\phi$. It is a component of the bound function $B_{\phi, \chi^2}(Q)$.   \label{tab:Chi2.resultoutline}}
\begin{center}
\bgroup
\footnotesize
\def\arraystretch{0.99}%
{\setlength{\tabcolsep}{0.2mm}
\begin{tabular}{|c|c|c|c|}
\hline
 \makecell{\textbf{Dist-} \\ \textbf{ance} \\ \textbf{fn $\phi$}} & \multicolumn{3}{c|}{\normalsize\textbf{Theoretical Aspects}}  \\
\cline{2-4} 
  & $\mathbf{\mathcal{I}^{R}_{\phi}} (m)$ & \textbf{Convexity} & \makecell{\textbf{Global min}/\textbf{Fixed Point (FP)} \\ (for uniform $P$)} \\
\hline 
$\phi_{\text{lin}}$ & \makecell{$\frac{1}{4m\delta}$  \\(due to \\ \cite{begin2016pac})} & Convex & \makecell{$q^{\ast}_{i, \text{lin}, \chi^2} (H^{\ast})= \left[ 1 + \frac{\left( \frac{\sum_{i = 1}^{H^{\ast}}\hat{l}_i}{H'} - \hat{l}_{i} \right)}{\sqrt{\frac{H}{H^{\ast}4m\delta} - \hat{var}_{H^{\ast}}(\hat{l})}} \right] \frac{1}{H} $ \\  (Global min)}
\\ \hline
$\phi_{\text{sq}}$ & $\frac{12m - 11}{16m^3 \delta}$ & \makecell{shown non-convex; \\
 Quasi-convex \\ under a condition}& $\begin{array}{lc}
 & q^{FP}_{i, \text{sq}, \chi^2} (H^{\ast}) = \frac{1}{H^{\ast}} \\
  + &\frac{2 \left( \sum_{i = 1}^{H^{\ast}}  (q^{FP}_{i, \text{sq}, \chi^2} (H^{\ast}))^2 \right)^{\frac{3}{4}}}{\sqrt[4]{\frac{(12m -11)H}{16m^3\delta}}} \left(\frac{\sum_{i = 1}^{H^{\ast}} \hat{l}_i}{H^{\ast}} - \hat{l}_i \right)
 \end{array}$
\\[4mm] \hline
$kl$ & \makecell{computed \\ based on form \\given by \cite{begin2016pac}}& \makecell{Non-convex; \\ Difference of \\ convex functions \\ (DC)}& \makecell{\hspace{-1cm}$q^{FP}_{i, \text{kl}, \chi^2} (H^{\ast})$ satisfies: $q_i = p_i \left(\sum_{i=1}^{H^{\ast}} q_i^2 \right) \times $ \\ $\left[ 1 + \frac{\left(\sum_{i=1}^{H^{\ast}}\hat{l}_iq_i - \hat{l}_i \right)}{\sqrt{\frac{H\left(\sum_{i=1}^{H^{\ast}} q_i^2\right) \mathcal{I}^{R}_{\text{kl}}(m, 2)}{\delta}} }  \ln\left( \frac{(1 - r) \sum_{i = 1}^{H^{\ast}} \hat{l}_i q_i}{r(1-\sum_{i = 1}^{H^{\ast}} \hat{l}_i q_i)}\right) \right]$}\\
\hline
\end{tabular}}
\egroup
\end{center}

\begin{center}
\centering
\bgroup
\footnotesize
\def\arraystretch{0.99}
{\setlength{\tabcolsep}{0.2mm}
\begin{tabular}{|c|c|c|c|}
\hline
 \makecell{\textbf{Dist-} \\ \textbf{ance} \\ \textbf{fn $\phi$}} &  \multicolumn{3}{c|}{ \textbf{Computations} (for uniform $P$) }\\
\cline{2-4} 
  &  \textbf{Solver (\texttt{Ipopt}) output}  & \textbf{Global min} & \textbf{Fixed Point (FP)}\\
\hline 
$\phi_{\text{lin}}$ &  Identifies global minima & \makecell{Identified analytically} & Not required\\ \hline
\makecell{$\phi_{\text{sq}}$} & \makecell{Identifies a unique (local) minimum \\
 even with different initializations} & \makecell{closed form \\ may not exist} & Matches solver output \\ \hline
$kl$ & \makecell{Identifies multiple local minima  \\
 with different  initializations; \\
 throws up error for  \\ moderate and large $H$} &  \makecell{closed form \\ may not exist} & \makecell{Identifies a unique \\ stationary point even with \\ different initializations} \\
\hline
\end{tabular}}
\egroup
\end{center}
\end{table}
}

\section{PAC-Bayesian Bound Minimization, Optimal Posteriors and the Fixed Point Approach} \label{secn:genericPACBbnd}

Classical version of PAC-Bayesian theorem is derived using Donsker-Varadhan inequality
 \footnote{The Donsker-Varadhan inequality can be stated as below:
 \begin{lem}[KL divergence change of measure \cite{begin2016pac}] 
 For any set $\mathcal{H}$, for any two distributions $P$ and $Q$ on $\mathcal{H}$, and for any measurable function $\phi: \mathcal{H} \rightarrow \mathbb{R}$, we have: $\mathbb{E}_{h \sim Q} \phi(h) \leq KL[Q||P] + \ln \left( \mathbb{E}_{h \sim P} e^{\phi(h)} \right)$.
 \end{lem}
 }
\cite{begin2016pac,seldin2012PACBMartingales}
for change of measure which is based on KL-divergence between the two distributions. 
A new version of PAC-Bayesian results has been discovered by \cite{begin2016pac} which involves change of measure guided by R{\'e}nyi divergence between the two distributions:

\begin{theorem}\cite{begin2016pac} For any data distribution $\mathcal{D}$ over input space $\mathcal{X} \times \mathcal{Y}$, the following bound holds for any prior $P$ over the set of classifiers $\mathcal{H}$, for any $\alpha > 1$ and any $\delta \in (0, 1)$, where the probability is over random i.i.d. samples $S_m = \{(x_i, y_i) | i = 1, \ldots , m\}$ of size $m$ drawn from $\mathcal{D}$, for any convex function $\phi : [0,1] \times [0, 1] \rightarrow \mathbb{R}$:
\begin{equation}
\mathbb{P}_{S_m} \left \lbrace \phi\left(\mathbb{E}_Q [\hat{l}], \mathbb{E}_Q[l] \right) \leq \left[ \mathbb{E}_{h \sim P} \left( \frac{Q(h)}{P(h)} \right)^{\alpha} \right]^{\frac{1}{\alpha}} \left[ \frac{\mathcal{I}^{R}_{\phi}(m, \alpha')}{\delta} \right]^{\frac{1}{\alpha'}} \right \rbrace \geq 1 - \delta. \label{eqn:expoPACB_Renyi}
\end{equation}
where $\alpha' = \frac{\alpha}{\alpha - 1}$ and $
\mathcal{I}^{R}_{\phi}(m, \alpha') := \sup_{l \in [0, 1]} \left[ \sum_{k = 0}^m \binom{m}{k} l^k (1- l)^{m - k}  \phi \left(\frac{k}{m}, l \right)^{\alpha'}\right].$
Here, $Q$ is an arbitrary posterior distribution on $\mathcal{H}$, which may depend on the sample $S$ and on the prior $P$. $\mathbb{E}_Q [\hat{l}] := \mathbb{E}_{h \sim Q} \sum_{i = 1}^m \frac{1}{m}[l(h, \mathbf{x}_i, y_i)]$ denotes the averaged empirical risk  and $\mathbb{E}_Q [l] := \mathbb{E}_{h \sim Q} \mathbb{E}_{(\mathbf{x}, y) \sim \mathcal{D}} [l]$ denotes averaged true risk of classifiers in $\mathcal{H}$ computed using a loss function, $l(h, \mathbf{x}, y): \mathcal{H} \times \mathcal{X} \times \mathcal{Y} \rightarrow [a, b)$ (here, $0 \leq a < b$).
\end{theorem} 

We are interested in identifying the optimal posteriors for different choices of distance functions for the case of $\alpha = 2$, which can be related to the chi-squared divergence measure, $\chi^2(Q||P) := \mathbb{E}_{h \sim P} \left[  \left( \frac{Q(h)}{P(h)}\right)^2 - 1 \right]$ between distributions $Q$ and $P$ \cite{begin2016pac}.

\subsection{Optimal posteriors via PAC-Bayesian bound minimization}
PAC-Bayesian theorem \eqref{eqn:expoPACB_Renyi} gives a high probability upper bound on averaged true risk, $\mathbb{E}_Q[l]$ assuming distance function $\phi(\mathbb{E}_Q[\hat{l}], \cdot)$ is invertible for given $\mathbb{E}_Q[\hat{l}]$:
\begin{equation}
\begin{split}
B_{\phi, \chi^2} (Q) &\equiv B_{\phi, \chi^2} (\mathbb{E}_Q[\hat{l}], S_m, \delta, P) \\ 
&= f_{\phi}\left(\mathbb{E}_Q[\hat{l}], \phi^{-1}_{\mathbb{E}_Q[\hat{l}]} \left( \sqrt{ (\chi^2(Q||P) + 1) \frac{\mathcal{I}^{R}_{\phi}(m, 2)}{\delta} } \right) \right), \label{eqn:BphiChi2}
\end{split}
\end{equation}
where $\phi^{-1}_{\mathbb{E}_Q[\hat{l}]}(K) = b$ implies $\phi(\mathbb{E}_Q[\hat{l}], b) = K$ for some $b \in (0, 1)$ and a given $K > 0$. Generally $f_{\phi}(\cdot, \cdot)$ is the sum of its arguments except when $\phi$ is KL-distance function. That is, bound function $B_{\phi, \chi^2} (Q)$ is the sum of averaged empirical risk, $\mathbb{E}_Q[\hat{l}]$, and a model complexity term which depends on system parameters, $S_m, \delta$ and $P$.  We are interested in the determining an optimal posterior distribution $Q^{\ast}_{\phi, \chi^2}$ which minimizes bound $B_{\phi, \chi^2}(Q)$ for a given distance function $\phi$. 

\subsection{The fixed point approach to characterize PAC-Bayesian optimal posterior}
To characterize the minimum of $B_{\phi, \chi^2}(Q)$, we make use of the first order KKT conditions which are necessary for a stationary point of a non-convex problem. These KKT conditions require the objective function and the active constraints to be differentiable at the local minimum. We derive fixed point (FP) equations for the optimal posterior using the partial KKT system. These FP equations use KKT system with strict positivity constraints due to which complementary slackness conditions are automatically satisfied; hence called `\textit{partial}' KKT system. The computations illustrate that these FP equations always converge to a stationary point at a very fast rate, even for a large classifier set when a non-convex solver fails to identify a solution. (Please see Table \ref{appdx_tab:Bnd.klChi2} and Table \ref{appdx_tab:TestErr.klChi2} in Appendix D for an illustration of such cases.)

\paragraph{Framework \label{secn:framework}}
We work with a finite set  of classifiers: $\mathcal{H} = \lbrace h_i \rbrace_{i = 1}^H$ of size $H$. The prior, $P = \lbrace p_i \rbrace_{i =1}^H$ and posterior, $Q =\lbrace q_i \rbrace_{i =1}^H$ are discrete distributions on $\mathcal{H}$, where $p_i, q_i \geq 0 \; \forall i = 1, \ldots, H$ with $\sum_{i = 1}^H p_i = 1$ and $\sum_{i = 1}^H q_i = 1$. For differentiability required by KKT conditions, our objective function should have open domain, that is, the interior of the $H$-dimensional probability simplex: $int(\Delta^H) = \lbrace (q_1, \ldots, q_H) \vert q_i > 0 \; \forall i = 1, \ldots, H; \sum_{i = 1}^H q_i = 1 \rbrace$. In computations, we consider $q_i \geq \epsilon \; \forall i = 1, \ldots, H$ for $\epsilon > 0$ to ensure existence of a minimizer in $int(\Delta^H)$. Our FP equations are derived using partial KKT system on $int(\Delta^H)$.

\section{Optimal posterior, $Q^{\ast}_{\phi, \chi^2}$,  for uniform prior \label{secn:OptQ_unifP}}
We consider the special case of uniform prior on entire $\mathcal{H}$. We want to identify the optimal posterior $Q^{\ast}_{\phi, \chi^2}$ with the $H$-dimensional probability simplex, $\Delta^H$, as the feasible region. We show below that it is enough to restrict the search space to certain subsets of $\Delta^H$. This reduces the computational complexity of the search \textit{from exponential scale to linear scale}.

\begin{theorem} \label{thm:increasing.subsets}
Consider a uniform prior distribution on the set $\mathcal{H}$ of classifiers, and a given set of posterior weights $Q = \lbrace q_j \rbrace_{j = 1}^{H'}$. We have three choices of distance function $\phi \in \lbrace \phi_{\text{lin}}, \phi_{\text{sq}}, kl \rbrace$. Then among all subsets $\mathcal{H}' \subset \mathcal{H}$ of size $H'$, the smallest bound value $B_{\phi, \chi^2}(Q, \mathcal{H}')$ corresponding to the given posterior weights $Q$ is achieved when $\mathcal{H}'$ is the subset formed by the first $H'$ elements of the ordered set of classifiers ranked by non-decreasing empirical risk values, $\hat{l}_1 \leq \hat{l}_2 \leq \ldots \leq \hat{l}_{H}$.
\end{theorem}
\begin{proof}
We first consider the case of linear and squared distance based bounds. Under the given set up, these bound functions are defined as follows:
\begin{align}
B_{\text{lin}, \chi^2}(Q, \mathcal{H}') &:= \sum_{i \in \mathcal{H}'} \hat{l}_iq_i + \sqrt{\frac{ \left(\sum_{i \in \mathcal{H}'} q_i^2 \right) H}{4m\delta}}. \label{eqn:BlinChi2_unifP.genQ} \\
B_{\text{sq}, \chi^2}(Q, \mathcal{H}') &:= \sum_{i\in \mathcal{H'}} \hat{l}_i q_i + \sqrt[4]{H\left( \sum_{i \in \mathcal{H'}} q_i^2 \right) \left(\frac{12m -11}{16m^3\delta} \right)}. \label{eqn:BsqChi2_unifP.genQ}
\end{align} 

For a given set of posterior weights $\lbrace q_j \rbrace_{j = 1}^{H'}$, the second terms of $B_{\text{lin}, \chi^2}(Q, \mathcal{H}')$ and $B_{\text{sq}, \chi^2}(Q, \mathcal{H}')$ are invariant of the support set $\mathcal{H}'$ as long as its cardinality is $H'$. Thus the value of the bound depends on the common first term which is a sum of positive quantities. For given weights $\lbrace q_j \rbrace_{j = 1}^{H'}$, the bounds \eqref{eqn:BlinChi2_unifP.genQ} and \eqref{eqn:BsqChi2_unifP.genQ} are the smallest when the sum $\sum_{i \in \mathcal{H}'} \hat{l}_i q_i$ is minimized. This will happen when $\mathcal{H}'$ consists of classifiers with smallest $H'$ values in the set $\lbrace \hat{l}_i \rbrace_{i = 1}^{H}$.
Furthermore, if the elements of $\mathcal{H}'$ are ordered by non-decreasing empirical risk values, $\hat{l}_1 \leq \hat{l}_2 \leq \ldots \leq \hat{l}_{H'}$, the posterior weights should be ordered non-increasingly. Hence, the claim of the theorem holds true.

Now, for the KL-divergence as a distance function, the bound value, $r$, is the solution to following two equations:
\begin{align}
&kl\left( \sum_{i \in \mathcal{H}'} \hat{l}_i q_i , r \right) = \sqrt{\frac{H \left(\sum_{i \in \mathcal{H}'} q_i^2\right) \mathcal{I}^{R}_{\text{kl}}(m, 2)}{\delta}} \label{eqn:klDCcons.unifP}\\
& r \geq \sum_{i \in \mathcal{H}'} \hat{l}_iq_i \label{eqn:r>EQl.unifP}
\end{align}
The right hand side term of \eqref{eqn:klDCcons.unifP} is invariant of support $\mathcal{H}'$ as long as it is of size $H'$. Let $\hat{L} := \sum_{i \in \mathcal{H}'} \hat{l}_iq_i$, then \eqref{eqn:klDCcons.unifP} is an implicit function of variables $\hat{L}$  and $r$. Using implicit function theorem, we have
\begin{equation}
\frac{dr}{d\hat{L}} = \frac{- \partial kl/\partial \hat{L}}{ \partial kl/\partial r} = \frac{\ln \frac{\hat{L}}{r} - \ln \frac{1 - \hat{L}}{1 - r}}{\frac{\hat{L}}{r} - \frac{1 - \hat{L}}{1 - r}}
\end{equation}
Using \eqref{eqn:r>EQl.unifP} and strict monotonicity of natural logarithm function, we can claim that $\frac{dr}{d\hat{L}} > 0$. That is, the bound $r$ is a strictly increasing function of $\hat{L} := \sum_{i \in \mathcal{H}'} \hat{l}_iq_i$ under the given set up. To find the least $r$ for a given $Q(\mathcal{H}') = \lbrace q_j \rbrace_{j = 1}^{H'}$, we need to find the least $\sum_{i \in \mathcal{H}'} \hat{l}_iq_i$ on all possible subsets $\mathcal{H}'$. This happens when $\mathcal{H}'$ is the subset formed by the first ordered $H'$ elements $\hat{l}_1 \leq \hat{l}_2 \leq \ldots \leq \hat{l}_{H'}$. Hence proved.
\end{proof}

\begin{cor}
As a consequence of the above Theorem \ref{thm:increasing.subsets}, for determining the (global) optimal posterior $Q^{\ast}_{\phi, \chi^2}$, it is sufficient to compare the bound values corresponding to the best posteriors on ordered subsets of $\mathcal{H}$, ranked by non-decreasing $\hat{l}_i$ values. These ordered subsets can be uniquely identified by their size. An ordered subset of size 1 is $\lbrace \hat{l}_1 \rbrace$, of size 2 is $\lbrace \hat{l}_1, \hat{l}_2 \rbrace$ and so on. Thus there exists an isomorphism between the set $\lbrace 1, \ldots, H \rbrace$ (which denote the subset size) and the family of ordered increasing subsets of $\mathcal{H}$.
\end{cor}

\begin{minipage}{0.48\textwidth}
\begin{algorithm}[H]
\DontPrintSemicolon
\KwIn{$ m, \delta, H, \lbrace \hat{l}_i \rbrace_{i =1}^H$}
\KwOut{$Q^{\ast}_{\phi, \chi^2}$}
Define an array $B^{\ast}_{\phi, \chi^2}[\ldots]$ of size $H$\;
$B^{\ast}_{\phi, \chi^2}[1] \gets f_{\phi}\left(\hat{l}_1, \phi^{-1}_{\hat{l}_1} \left( \sqrt{ \frac{\mathcal{I}^{R}_{\phi}(m, 2)}{p_1\delta} } \right) \right)$ \;
\For{$H' = 2, \ldots, H$}{
flag $\gets 0$ \;
Identify $Q^{\ast}_{\phi, \chi^2}$ for $H'$ via \eqref{eqn:OptQlinChi2.subset.unifP} or \eqref{eqn:OptQsqChi2.subset.unifP} or \eqref{eqn:qFP_klChi2.unifP} \;
\For{$i=1, \ldots, H'$}{
\If{$q^{\ast}_{i, \phi, \chi^2} < 0$ \label{step:infeasibleQphiChi2}}{
flag $\gets 1$\;
\KwBreak
}
}
\If{flag $= 1$}{\KwBreak}
Compute $B^{\ast}_{\phi, \chi^2}[H']$ using $Q^{\ast}_{\phi, \chi^2}$ in \eqref{eqn:BphiChi2} \;
}
$H^{\ast} \gets \arg \min_{H'}B^{\ast}_{\phi, \chi^2}[H']$\;
Identify $Q^{\ast}_{\phi, \chi^2}$ for $H^{\ast}$ via \eqref{eqn:OptQlinChi2.subset.unifP} or \eqref{eqn:OptQsqChi2.subset.unifP} or \eqref{eqn:qFP_klChi2.unifP} \;
\Return $Q^{\ast}_{\phi, \chi^2}$ 
\caption{\textsc{OptQ $\phi \mathrm{-}\chi^2$ For Uniform Prior}: Algorithm for finding optimal posterior for the PAC-Bayesian bound based on $\chi^2$-divergence when prior is uniform distribution \label{algo:optQsubset.unifP}}
\end{algorithm}
\end{minipage}
~
\begin{minipage}{0.48\textwidth}
\paragraph{Correctness of Algorithm \textsc{OptQ $\phi-\chi^2$ For Uniform Prior}}
We want to determine the globally optimal posterior $Q^{\ast}_{\phi, \chi^2}$ that has the minimum bound value $B_{\phi, \chi^2}(Q)$ over the $H$-dimensional probability simplex, $\Delta^H$. Using the result of Theorem \ref{thm:increasing.subsets}, we can confine the search to a much smaller space of posteriors with support  on a family of increasing ordered subsets of $\mathcal{H}$. These ordered subsets are defined by their size. For example, an ordered subset of size $H' \in [H]$ comprises of the \textit{lowest} $H'$ values in the set $\lbrace\hat{l}_i \rbrace_{i =1}^{H}$. This restricted space of posteriors, say $\Delta^{\text{ord}} \subset \Delta^H$, is a union of convex sets of posteriors with supports on ordered subsets defined above. Due to increasing subset relation between consecutive supports, this union is also a convex set. The search space $\Delta^{\text{ord}}$ is a restriction of $\Delta^H$, yet consists of uncountably many posteriors. We refine the search further by localizing to optimal posteriors $Q^{\ast}_{\phi, \chi^2}(H')$ on each increasing ordered subset and comparing their bound values, $B^{\ast}_{\phi, \chi^2} (H')$ to find the minimum. Thus, an exponential search on restricted posterior space is simplified to a finite \textit{linear} search on the support size. When using FP scheme, we also need to verify that $Q^{\ast}_{\phi, \chi^2}(H')$ satisfies positivity constraints. We denote the support size of  the optimal posterior $Q^{\ast}_{\phi, \chi^2}$ by $H^{\ast} \in [H]$. Therefore, for determining $Q^{\ast}_{\phi, \chi^2}$ in $\Delta^{\text{ord}}$, it is sufficient to search for $H^{\ast}$ in the set $\lbrace 1, \ldots, H\rbrace$.
\end{minipage}

\section{Optimal PAC-Bayesian Posterior using Linear Distance Function \label{secn:linChi2}} 
As a basic case, we can consider linear distance function, $\phi_{\text{lin}}(\hat{l}, l) = l - \hat{l}$ for $\hat{l}, l \in [0,1]$. The optimal posterior $Q^{\ast}_{\text{lin}, \chi^2}$ that bounds the \textit{unknown} averaged true risk of a stochastic classifier, is obtained via the following minimization problem for the bound $B_{\text{lin},\chi^2}(Q)$ identified by \cite{begin2016pac}. 
\begin{equation} \label{eqn:Blinchi2OP}
\begin{split}
\min_{Q = (q_1, \ldots, q_H ) \in \Delta^H} &B_{\text{lin},\chi^2}(Q):= \sum_{i = 1}^H \hat{l}_i q_i + \sqrt{\frac{\sum_{i =1}^H \frac{q_i^2}{p_i}}{4m\delta}} 
\end{split}
\end{equation}
\begin{theorem} \label{thm:convexity.BlinChi2}
The bound function $B_{\text{lin},\chi^2}(Q)$ (identified by \cite{begin2016pac}) is a strictly convex function and hence, \eqref{eqn:Blinchi2OP} is a convex program with a unique global minimum.
\end{theorem}
Proof (in Section \ref{appdx_secn:BlinChi2OP} in the appendix) uses first order convexity property. 
\subsection{Optimal posterior, $Q^{\ast}_{\text{lin}, \chi^2}$,  for uniform prior \label{secn:QlinChi2_unifP}}
We identify the optimal posterior and the optimal bound value for linear distance function by exploiting convexity of $B_{\text{lin}, \chi^2}(Q)$. Proofs for Theorem \ref{thm:OptQlinChi2.subset} and Theorem \ref{thm:decreasing.optBlinChi2.Hdash} stated below are in Sections \ref{appdx_secn:optQlinChi2} and \ref{appdx_secn:QlinChi2_unifP} in the appendix.
\begin{theorem}[Optimal posterior on an ordered subset support] \label{thm:OptQlinChi2.subset}
When prior is uniform distribution on $\mathcal{H}$, among all the posteriors with support as subset of $\mathcal{H}$ of size exactly $H'$, the best posterior denoted by $Q^{\ast}_{\text{lin}, \chi^2}(H')$ has the support on the ordered subset $\mathcal{H'}_{\text{ord}} = \lbrace \lbrace \hat{l}_i \rbrace_{i=1}^{H'} \vert \hat{l}_1 \leq \hat{l}_2 \leq \ldots \leq \hat{l}_{H'}\rbrace$ consisting of smallest $H'$ values in $\mathcal{H}$. The optimal posterior weights are determined as follows:
\begin{equation}
    q^{\ast}_{i, \text{lin}, \chi^2} (H')=
    \begin{cases} \left[ 1 + \frac{\left( \frac{\sum_{i = 1}^{H'}\hat{l}_i}{H'} - \hat{l}_{i} \right)}{\sqrt{\frac{H}{H'4m\delta} - \hat{var}_{H'}(\hat{l})}} \right] \frac{1}{H} \quad &i=1, \ldots, H' \\
0 \quad &i= H'+1, \ldots, H,
    \end{cases}
\label{eqn:OptQlinChi2.subset.unifP}
\end{equation}
where $\hat{var}_{H'}(\hat{l})=  \frac{1}{H'}\sum\limits_{i =1}^{H'} \left( \frac{\sum_{i =1}^{H'} \hat{l}_i}{H'}- \hat{l}_i\right)^2 = \frac{\sum_{i =1}^{H'} \hat{l}^2_i}{H'} - \left(\frac{\sum_{i =1}^{H'} \hat{l}_i}{H'}\right)^2$ is the variance of the values in $\mathcal{H'}_{\text{ord}}$. We require that $H'$ is such that $\frac{H}{H'4m\delta} - \hat{var}_{H'}(\hat{l}) > 0$ so that $Q^{\ast}_{\text{lin}, \chi^2}(H')$ is defined and for feasibility, $q^{\ast}_{i, \text{lin}, \chi^2}(H') > 0$ for $i=1, \ldots, H'$.
\end{theorem}
Using the closed form expression \eqref{eqn:OptQlinChi2.subset.unifP}, we can identify the optimal posterior $Q^{\ast}_{\text{lin}, \chi^2}$ via Algorithm \ref{algo:optQsubset.unifP}.
\begin{rem}
For given values of $H, H', m$ and $\hat{var}_{H'}(\hat{l})$, the upper bound on $\delta$ is related to sparseness of the optimal posterior $Q^{\ast}_{\text{lin}, \chi^2}$. A higher $\delta$ diminishes the effect of divergence term $\sum_{i =1}^H \frac{q_i^2}{p_i}$ and allows sparse solutions. 
\end{rem}
\begin{theorem} \label{thm:decreasing.optBlinChi2.Hdash}
The bound value of the best posterior $Q^{\ast}_{\text{lin}, \chi^2}(H')$ on an ordered subset of size $H'$,
\begin{equation}
B^{\ast}_{\text{lin}, \chi^2}(H') := B_{\text{lin}, \chi^2} \left( Q^{\ast}_{\text{lin}, \chi^2}(H') \right) = \frac{\sum_{i =1}^{H'} \hat{l}_i}{H'} + \sqrt{\frac{H}{H'4m\delta} - \hat{var}_{H'}(\hat{l})}, \label{eqn:optBlinChi2.Hdash}
\end{equation} 
 is decreasing function of $H' \leq H^{\ast}$, the support size of globally optimal posterior $Q^{\ast}_{\text{lin}, \chi^2}$.
\end{theorem}

\section{Optimal PAC-Bayesian Posterior using Squared Distance Function \label{secn:sqChi2}}
PAC-Bayesian bound for squared distance, $\phi_{\text{sq}} \left(\hat{l}, l \right) =  \left(\hat{l} - l \right)^2$ for $\hat{l}, l \in [0,1]$ is identified below. 
We first to need to identify $\mathcal{I}^{R}_{\text{sq}}(m, 2)$ for a given sample size $m$. Details of the derivation are in Appendix \ref{appdx_secn:sqChi2}.

\begin{lem} \label{lem:I_R_sq}
For a given sample size, $m$, $\mathcal{I}^{R}_{\text{sq}}(m, 2) := \sum_{k=0}^{m}{m \choose k} 0.5^m e^{2m\left(\frac{k}{m}-0.5\right)^4} = \frac{12m -11}{16m^3}.$ 
\end{lem}

\begin{theorem}
For a finite set of classifiers, $\mathcal{H}$, PAC-Bayesian bound, $B_{sq, \chi^2}(Q)$, based on squared distance function with $\chi^2$-divergence measure is given by:
\begin{equation}
  B_{sq, \chi^2}(Q) := \sum_{i = 1}^H \hat{l}_i q_i + \sqrt[4]{\left( \sum_{i = 1}^H \frac{q_i^2}{p_i} \right) \left(\frac{12m -11}{16m^3\delta} \right)} \label{eqn:BsqChi2}
\end{equation}
\end{theorem}
 \begin{proof}
 Using the PAC-Bayesian statement in \eqref{eqn:expoPACB_Renyi} with $\phi_{\text{sq}} \left(\hat{l}, l \right) =  \left(\hat{l} - l \right)^2$ for finite $\mathcal{H}$, we can obtain $B_{sq, \chi^2}(Q)$ in using \eqref{eqn:BphiChi2} and $\mathcal{I}^{R}_{\text{sq}}(m, 2)$ identified via Lemma \ref{lem:I_R_sq}.
 \end{proof}
We want to determine the optimal posterior $Q^{\ast}_{sq, \chi^2}$ which minimizes $B_{sq, \chi^2}(Q)$ over $\Delta^H$. This bound function turns out to be non-convex in $Q$.  

\begin{theorem} \label{thm:nonconvexBsqChi2}
The bound function, $B_{sq, \chi^2}(Q) = \sum_{i = 1}^H \hat{l}_i q_i + \sqrt[4]{\left( \sum_{i = 1}^H \frac{q_i^2}{p_i} \right) \left(\frac{12m -11}{16m^3\delta} \right)} $ is non-convex.
\end{theorem}
We show this non-convexity even when $P \sim Unif(\mathcal{H})$ via counter examples violating first order convexity property in Section \ref{appdx_secn:nonconvexity.BsqChi2} in the appendix.

\begin{rem} Computationally this bound minimization problem for \eqref{eqn:BsqChi2} is observed to have a single solution. We used bordered Hessian test to verify that the solution obtained is a local minimum. This motivates the following. The bound $B_{\text{sq}, \chi^2}(Q)$ is shown to be quasi-convex under a condition on system parameters. 
\end{rem}
\begin{propn} \label{propn:quasiconvex.condn.BsqChi2}
The bound function $B_{\text{sq},\chi^2}(Q)$ is strictly quasi-convex if the following condition holds for any $Q, Q'$ for each $\alpha \in (0,1)$:
\begin{multline*}
 \bigg(\sqrt[4]{\frac{12m-11}{16m^3 \delta}}\bigg)\left[ \sqrt[4]{\sum_{i=1}^{H} \frac{(\alpha q_i + (1-\alpha)q'_i)^2}{p_i}} - \sqrt[4]{\sum_{i=1}^{H} \frac{q_i^2}{p_i}} \right] < (1-\alpha)(E_{Q} [\hat{l}]-E_{Q'} [\hat{l}])
\end{multline*}
and hence a local minimum to the bound minimization problem for the bound \eqref{eqn:BsqChi2} is also a global minimum.
\end{propn}
The proof of this proposition is in Appendix \ref{appdx_secn:nonconvexity.BsqChi2}.

\subsection{The posterior based on fixed point scheme, $Q^{FP}_{\text{sq}, \chi^2}$, for uniform prior \label{secn:QsqChi2_unifP}}
For uniform prior set up, we derive FPE for minimizing \eqref{eqn:BsqChi2} on an ordered subset support of size $H'$.
\begin{theorem}[Optimal posterior on an ordered subset support] \label{thm:OptQsqChi2.subset}
When prior is uniform distribution on $\mathcal{H}$, among all the posteriors with support as subset of size exactly $H'$, the best posterior denoted by $Q^{\ast}_{\text{sq}, \chi^2}(H')$ has the support on the ordered subset $\mathcal{H'}_{\text{ord}} = \lbrace \lbrace \hat{l}_i \rbrace_{i=1}^{H'} \vert \hat{l}_1 \leq \hat{l}_2 \leq \ldots \leq \hat{l}_{H'}\rbrace$ consisting of smallest $H'$ values in $\mathcal{H}$. The optimal posterior weights $\lbrace q^{\ast}_{i, \text{sq}, \chi^2} (H') \rbrace$ are determined as the solution to the following fixed point equation in $\lbrace q_{i}(H') \rbrace_{i = 1}^{H'}$:
\begin{equation}
    q_{i} (H')=
    \begin{cases} \frac{1}{H'} + \frac{2 \left( \sum_{i = 1}^{H'}  (q_{i} (H'))^2 \right)^{\frac{3}{4}}}{\sqrt[4]{\frac{(12m -11)H}{16m^3\delta}}} \left(\frac{\sum_{i = 1}^{H'} \hat{l}_i}{H'} - \hat{l}_i \right)  \quad &i=1, \ldots, H' \\
0 \quad &i= H'+1, \ldots, H.
    \end{cases}
\label{eqn:OptQsqChi2.subset.unifP}
\end{equation}
under the feasibility condition that $q_{i}(H') > 0$ for $i=1, \ldots, H'$.
\end{theorem}
Proof of this theorem is in Section \ref{appdx_secn:optQsqChi2} in the appendix. In our computations, FP iterates in \eqref{eqn:OptQsqChi2.subset.unifP} do converge to a solution. Using them, we can identify the optimal posterior $Q^{\ast}_{\text{sq}, \chi^2}$ via Algorithm \ref{algo:optQsubset.unifP}.

\section{Optimal PAC-Bayesian Posterior using KL-distance Function, $kl(\cdot, \cdot)$} \label{secn:klChi2}
Chi-squared divergence based PAC-Bayesian bound using the distance function $kl(\hat{l} , l) = \hat{l} \ln \left( \frac{\hat{l}}{l}\right) + (1 - \hat{l}) \ln \left( \frac{1 - \hat{l}}{1 - l}\right)$ (for $\hat{l}, l \in [0, 1]$) is:
\begin{equation}
B_{\text{kl}, \chi^2}(Q) = \sup_{r \in (0, 1)} \left\lbrace r : kl\left(\mathbb{E}_Q [\hat{l}], r \right) \leq \sqrt{\frac{(\chi^2(Q||P) + 1) \mathcal{I}^{R}_{\text{kl}}(m, 2)}{\delta}} \right\rbrace
\end{equation}
 \added{where $\mathcal{I}^{R}_{\text{kl}}(m, 2) := \sum\limits_{k=0}^{m}{m \choose k} l^k (1-l)^{m-k} \left(kl\left(\frac{k}{m},l\right)\right)^2$ should be computed first. For $m > 1028$, computation is difficult due to storage limitations in the range of floating point numbers. We notice that $\mathcal{I}^R_{\text{kl}} (m, 2)$ decreases with $m$ and hence, we can use $I^R_\text{kl}(1028, 2)$ as an upper approximation for $I^R_\text{kl}(m, 2)$ for $m > 1028$. Please refer to Table \ref{appdx_tab:I_R_kl} and Figure \ref{appdx_fig:I_R_kl} in Appendix \ref{appdx_secn:klChi2} for details.}

$kl( \cdot,  \cdot)$ is not a monotone function and so its inverse does not exist. Thus, $B_{\text{kl}, \chi^2}(Q) $ does not have an explicit form. However, we can employ a numerical root finding algorithm such as that described in \cite{PACBIntervals} (Algo. (\textsc{KLroots})) to obtain $B_{\text{kl}, \chi^2}(Q) $ for given system parameter values.

For a finite classifier space $\mathcal{H} = \lbrace h_i\rbrace_{i =1}^H$ with empirical risk values $\lbrace \hat{l}_i \rbrace_{i =1}^H$, the KL-distance bound minimization problem is:
\begin{subequations}  \label{eqn:BklChi2OP}
\begin{align}
&\min_{\substack{(q_1, \ldots, q_H) \in \Delta^H \\ r \in (0,1)}} r\\
\hspace{-4cm}
\text{s.t.} \quad & \left(\sum_{i = 1}^H \hat{l}_iq_i \right) \ln\left( \frac{\sum\limits_{i = 1}^H \hat{l}_iq_i}{r}\right) + \left(1 - \sum_{i = 1}^H \hat{l}_iq_i \right) \ln\left( \frac{ 1  - \sum\limits_{i = 1}^H \hat{l}_iq_i}{1 - r}\right) =\sqrt{\frac{\left(\sum\limits_{i=1}^H \frac{q_i^2}{p_i}\right) \mathcal{I}^{R}_{\text{kl}}(m, 2)}{\delta}} \label{eqn:klDCcons}  \\
& r \geq \sum_{i = 1}^H \hat{l}_iq_i  \label{eqn:r>EQl} 
\end{align} 
\end{subequations}
Here, $r$ is the right root of \eqref{eqn:klDCcons} for a given $\mathbb{E}_Q [\hat{l}]$. The above is known to be a non-convex problem with a difference of convex (DC) equality constraint \eqref{eqn:klDCcons}; and has multiple stationary points. This fact is illustrated in our computations, where
different initializations led to different stationary points. Using bordered Hessian test, we verified that these stationary points computed on our datasets by the solver or the FP equation \eqref{eqn:qFP_klChi2.unifP} given below are either local minima or saddle points.The constraint \eqref{eqn:r>EQl} is a strict inequality which is relaxed to ensure a solution in the closed domain. The iterative root finding algorithm adds to the computational complexity of the bound minimization algorithm.

\added{The objective function and constraints of the above bound minimization, \eqref{eqn:BklChi2OP} are either linear or difference of convex (DC) functions, hence it falls into the category of a DC program. We can make use of the convex-concave procedure (CCP), which is a powerful heuristic method used to find local solutions to DC programming problems \cite{CCP2016LippBoyd}. This procedure makes a linear approximation via supporting hyperplane to the second convex function in the DC function of the optimization program. This helps us convert the DC constraint into a convex constraint and hence the original DC program is reduced to a convex program which can be easily solved. The details of the CCP for solving \eqref{eqn:BklChi2OP} are in Appendix \ref{appdx_secn:klChi2CCP} and the related computations are in Table \ref{appdx_tab:Bnd.TestErr.FPnCCP.klChi2} in the appendix.  Our general observation is that the fixed point scheme that we derive outperforms CCP.}

\subsection{The posterior based on fixed point scheme, $Q^{FP}_{\text{kl}, \chi^2}$, for uniform prior}
We derive FPE for \eqref{eqn:BklChi2OP} for uniform prior when the support of $Q$ is an ordered subset of size $H'$. Proof is in Section \ref{appdx_secn:optQklChi2} of the appendix. This FPE is used in Algorithm \ref{algo:optQsubset.unifP} for determining $Q^{FP}_{\text{kl}, \chi^2}$.
\begin{theorem}[Optimal posterior on an ordered subset support] \label{thm:qFP_klChi2}
When prior is uniform distribution on $\mathcal{H}$, among all the posteriors with support as subset of size exactly $H'$, the best posterior denoted by $Q^{\ast}_{\text{sq}, \chi^2}(H')$ has the support on the ordered subset $\mathcal{H'}_{\text{ord}} = \lbrace \lbrace \hat{l}_i \rbrace_{i=1}^{H'} \vert \hat{l}_1 \leq \hat{l}_2 \leq \ldots \leq \hat{l}_{H'}\rbrace$ consisting of smallest $H'$ values in $\mathcal{H}$. A stationary point $Q^{\text{FP}}_{\text{kl}, \chi^2}(H')$ for \eqref{eqn:BklChi2OP} can be obtained as the solution to the following fixed point equation in $\lbrace q_{i} \rbrace_{i = 1}^{H'}$:
\begin{equation}
q_i = \frac{1}{Z_{\text{kl}, \chi^2}} \left(\sum_{i=1}^{H'} q_i^2 \right) \left \lbrace 1 + \frac{\left(\hat{l}_i  -\sum_{i=1}^{H'}\hat{l}_iq_i \right)}{\sqrt{\frac{H\left(\sum_{i=1}^{H'} q_i^2\right) \mathcal{I}^{R}_{\text{kl}}(m, 2)}{\delta}} } \left[ \ln\left( \frac{(1 - r) \sum_{i = 1}^{H'} \hat{l}_i q_i}{r(1-\sum_{i = 1}^{H'} \hat{l}_i q_i)}\right) \right] \right \rbrace 
\label{eqn:qFP_klChi2.unifP}
\end{equation}
for $i = 1, \ldots, H'$, where $Z_{\text{kl}, \chi^2}$ is a suitable normalization constant and $r$ is the solution to \eqref{eqn:klDCcons} and \eqref{eqn:r>EQl} for a given $Q = (q_1, \ldots, q_{H'}) \in \text{interior}( \Delta^{H'} )$.
\end{theorem}


\section{Choice of Regularization Parameter for SVMs} \label{secn:SVM_PACB}
For computations, we included nine datasets from UCI repository \cite{UCI:2017} with small to moderate number of examples (306 examples to 5463 examples) and small to moderate number of 
features (3 features to 57 features). The details about the number of features, numberof examples and class distribution of these datasets are listed in Table \ref{tab:dataset.details} \cite{sahu2019PACBKLarxiv}. These datasets span a variety ranging from almost linearly separable (Banknote, Mushroom and Wave datasets) to moderately inseparable (Wdbc, Mammographic and Ionosphere datasets) to inseparable data (Spambase, Bupa and Haberman datasets). SVMs on these datasets have varying ranges and degrees of variation in their empirical risk values. 

We consider a finite set of SVM regularization parameter values $\Lambda = \lbrace \lambda_i \rbrace_{i = 1}^{H}$, say, between $0$ and an upper bound $\lambda_0 > 0$, since small values of $\lambda_i$'s are preferable. The set $\Lambda$ is an arithmetic-geometric progression (AGP) with a logarithmic scale for $\lambda \in (0, 0.1)$ and a linear scale for $\lambda \geq 0.1$. \added{The logarithmic subset of $\Lambda$ is a union of 3 geometric series with ratios $\frac{1}{2}, \frac{1}{3}$ and $\frac{1}{5}$ each with elements truncated between $1e-10$ and 0.1. We use a lower bound $1e-10$ which is slightly away from 0 since the SVM classifiers become very close (almost indistinguishable) and give similar training and test errors for very small values of $\lambda$ which are almost zero. For values beyond 0.1, we use a arithmetic progression with spacing of 0.05. We use an upper bound $\lambda = 5$ on this arithmetic series, since on most datasets, SVMs generated by $\lambda \geq 5$ do not have good training and test error rates. (Graphical illustration of ranges and variation of training errors and test errors of the nine UCI datasets we considered on the chosen $\Lambda$ range is depicted in Section G in \cite{sahu2019PACBKLarxiv}.) }

SVM QP (with RBF kernels) was implemented using \texttt{ksvm} function in \texttt{kernlab} package \cite{kernlab} in { \em R (version 3.1.3 (2015-03-09))}. The Gaussian width parameter was estimated by \texttt{kernlab} using \texttt{sigest} function which estimates 0.1 and 0.9 quantiles of squared distance between the data points.

\begin{table}[]
\small
\setlength{\tabcolsep}{0.1em}
\begin{tabular}{|c|c|c|c|c|c|c|}
\hline
\textbf{Dataset} & \makecell{ \textbf{Number of} \\ \textbf{features}, $n$} & \makecell{\textbf{Number of } \\ \textbf{examples}}& \textbf{Pos/Neg} & \makecell{ \textbf{Training} \\ \textbf{set size}, $m$} & \makecell{ \textbf{Validation} \\ \textbf{set size}, $v$} & \makecell{ \textbf{Test} \\ \textbf{set size}, $t$} \\ \hline
\textbf{Spambase} & 57 & 4601 & 2788/1813 & 1840 & 1840 & 921 \\ \hline
\textbf{Bupa} & 6 & 345 & 176/169 & 138 & 138 & 69 \\ \hline
\textbf{Mammographic} & 5 & 830 & 427/403 & 332 & 332 & 166 \\ \hline
\textbf{\makecell{Wdbc}} & 30 & 569 & 357/212 & 227 & 227 & 115 \\ \hline
\textbf{Banknote} & 4 & 1372 & 610/762 & 548 & 549 & 275 \\ \hline
\textbf{Mushroom} & 22 (\textit{116 \footnote{after one-hot encoding for categorical features}})  & 5643 \footnote{after removing the rows with missing values from the data}& 3488/2155 & 2257 & 2257 & 1129 \\ \hline
\textbf{Ionosphere} & 34 & 351 & 225/126 & 140 & 140 & 71 \\ \hline
\textbf{Waveform} & 40 & 3308 \footnote{number of examples when class `0' is removed} & 1653/1655 & 1323 & 1323 & 662 \\ \hline
\textbf{Haberman} & 3 & 306 & 225/81 & 122 & 122 & 62\\ \hline
\end{tabular}
\caption[Details of various UCI datasets used for computational experiments]{\textbf{Details of various UCI datasets used for computational experiments} (Table 5 in Section G of \cite{sahu2019PACBKLarxiv}). We list the number of features $n$, total number of examples with distribution into positive and negative classes for each dataset. We also give the number of examples in training, validation and test sets, according to the random partition created by 0.4:0.4:0.2 ratio of the total dataset size.}
\label{tab:dataset.details}
\end{table}

Each of these datasets was partitioned such that 80\% of the examples formed a composition of training set and validation set (in equal proportion) used for constructing the set $\mathcal{H} = \lbrace h(\lambda_i) | \lambda_i \in \Lambda \rbrace_{i = 1}^H$ of SVM classifiers and remaining 20\% used for computing their test error rates. The training set size ($m$), validation set size ($v$) and test set size ($t$) are in the ratio $m:v:t = 0.4:0.4:0.2$. The role of the validation set is to compute the empirical risk $\hat{l}_i$ of the SVM $h(\lambda_i) \in \mathcal{H}$ which will be used for deriving the PAC-Bayesian bound.  We follow the scheme provided in \cite{begin2016pac,thiemann2016quasiconvexPACB} to generate the set $\mathcal{H}$. Each classifier $h(\lambda_i) \in \mathcal{H}$ is trained on $m$ training examples subsampled from this composite set and validated on the remaining $v$ examples. Overlaps between training sets of different classifiers are allowed. Same is true for their validation sets. (For further details about the dataset categorization, please refer to Section G.1 in \cite{sahu2019PACBKLarxiv}.)

The PAC-Bayesian bound minimization for finding the optimal posterior was implemented in AMPL Interface and solved using \texttt{Ipopt} software package \linebreak {\em (version 3.12 (2016-05-01))} \cite{ipopt}, a library for large-scale nonlinear optimization (\texttt{http://projects.coin-or.org/Ipopt}). All the computations were done on a machine equipped with 4 Intel Xeon 2.13 GHz cores and 64 GB RAM. %

\begin{table}[ht]
\caption{\textbf{PAC-Bayesian bounds and averaged test error rates for $Q^{\ast}_{\phi,\chi^2}$.} We compare bound values $B^{\ast}_{\phi,\chi^2}$ and average test error rates $T_{\phi,\chi^2}$ of optimal posteriors $Q^{\ast}_{\phi,\chi^2}$ for 3 distance functions:  KL-divergence $kl$, linear $\phi_{\text{lin}}$ and squared distances $\phi_{\text{sq}}$. For large sample size ($m \geq 1028$), the constant $\mathcal{I}^R_\text{kl}(m, 2)$ cannot be computed due to floating point storage limitations on the machine. So, we use an upper approximation: for $m > 1028$, $\mathcal{I}^R_\text{kl}(m, 2) \leq \mathcal{I}^R_\text{kl}(1028, 2) = 0.00037$ since $\mathcal{I}^R_\text{kl}(m, 2)$ is decreasing with $m$. Please see Table \ref{appdx_tab:TestErr.klChi2} in Appendix \ref{appdx_secn:SVM_PACB} for details. \textcolor{magenta}{$\star$} refers to values obtained using fixed point equation because the solver \texttt{Ipopt} does not converge to a solution for reasons like local infeasibility, Restoration Phase Failed, etc. Please see Appendix \ref{appdx_secn:SVM_PACB} for more such examples.
Lowest 10\% bound values and test error rates for each dataset are denoted in bold face. For a given posterior, we measure sparsity by the number of classifiers needed by its cumulative distribution function (CDF) to achieve a certain significance level. Concentration of a posterior is quantified in terms of its $\ell_2$ norm, which is equivalent to HHI score used for measuring market share of a firm in the industry \cite{Hirschman1945HHI,HHIwiki}. For almost separable datasets, the CDFs of the three posteriors are close to each other. $Q^{\ast}_{\text{kl}, \chi^2}$ have almost full support and low concentration. These posteriors give extremely loose bounds and are computationally expensive but have test error rates better than linear distance ones, $Q^{\ast}_{\text{lin}, \chi^2}$, on most datasets. Squared distance based posteriors, $Q^{\ast}_{\text{sq}, \chi^2}$, are sparse with relatively high concentration (relatively high $\ell_2$ norm) and have the tightest bounds and lowest test error rates. The contrast in test error rates is striking when the dataset yields classifiers with high variation in empirical risk values. See Section \ref{appdx_secn:sparsityconc} in the appendix for details. \label{tab:phiChi2Compare}}
\centering 
{ \footnotesize
\begin{tabular}{ |c | c c c ||c c c|}
\hline
\textbf{Dataset} & \multicolumn{3}{c||}{\textbf{PAC-Bayesian Bound}, $B^{\ast}_{\phi, \chi^2}$} & \multicolumn{3}{c|}{\textbf{Average Test Error}, $T_{\phi, \chi^2}$} \\
\cline{2-7}
  & $B^{\ast}_{\text{lin}, \chi^2}$ & $B^{\ast}_{\text{sq}, \chi^2}$ &$B^{\ast}_{\text{kl}, \chi^2}$ & $T_{\text{lin}, \chi^2}$ & $T_{\text{sq}, \chi^2}$ &$T_{\text{kl}, \chi^2}$
\\ 
\hline
Spambase  & 0.38054 & \textbf{0.30581} & 0.57082{\textcolor{magenta}{$\star$}}& 0.27190 & \textbf{0.19437} & 0.27345{\textcolor{magenta}{$\star$}}\\
Bupa   & 0.82183 & \textbf{0.60348} &  0.83864{\textcolor{magenta}{$\star$}}& 0.44564& \textbf{0.29715} & 0.42343{\textcolor{magenta}{$\star$}} \\
\makecell{Mammographic}  & 0.50276 & \textbf{0.38440} & 0.62780{\textcolor{magenta}{$\star$}}& \textbf{0.23308} & \textbf{0.23036} & \textbf{0.23266}{\textcolor{magenta}{$\star$}} \\
Wdbc  & 0.41631 & \textbf{0.27465} & 0.49508{\textcolor{magenta}{$\star$}} & \textbf{0.08152} & \textbf{0.07664} & \textbf{0.07808}{\textcolor{magenta}{$\star$}} \\
Banknote  & 0.22283 & \textbf{0.13487} & 0.26038{\textcolor{magenta}{$\star$}} & \textbf{0.00668} & \textbf{0.00661} & \textbf{0.00662}{\textcolor{magenta}{$\star$}} \\
Mushroom  & 0.10785 & \textbf{0.06449} & 0.18660{\textcolor{magenta}{$\star$}} & \textbf{0.00095} & \textbf{0.00092} & \textbf{0.00096}{\textcolor{magenta}{$\star$}}\\
Ionosphere  & 0.64273 & \textbf{0.38167} & 0.67611{\textcolor{magenta}{$\star$}} & 0.20960 & \textbf{0.06763} & 0.07974{\textcolor{magenta}{$\star$}}\\
Waveform  & 0.18565 & \textbf{0.12900} & 0.28978 &  \textbf{0.05120} & \textbf{0.05122} & \textbf{0.05127}\\
Haberman  & 0.70477 & \textbf{0.51731} & 0.73791{\textcolor{magenta}{$\star$}} & \textbf{0.28985} & \textbf{0.28979} & \textbf{0.28983}{\textcolor{magenta}{$\star$}}\\ 
\hline
\end{tabular}
}
\end{table}

\section{Comparison with KL-divergence based PAC-Bayesian posteriors and cross-validation method} \label{secn:compareKLCV}
We present a comparative study of $\chi^2$-divergence based optimal posteriors considered in this paper with optimal posteriors for KL-divergence based PAC-Bayesian bounds considered by \cite{pmlr-v101-sahu19a}. We also analyze the performance of these stochastic SVM classifiers governed by PAC-Bayesian posteriors with respect to the deterministic SVM classifier identified via cross-validation procedure.

\subsection{Comparison with KL-divergence based optimal PAC-Bayesian posteriors}
Given a set of base classifiers and having computed their empirical risk values, we observe the following differences and similarities between the optimal PAC-Bayesian posteriors due to the two divergence functions on this classifier set.
\begin{enumerate}[label = \roman*]
\item \textbf{Nature of optimal posteriors}: Both KL-divergence and $\chi^2$-divergence based optimal posteriors that exhibit \textit{decreasing} trend with respect to the empirical risk values. That is to say that higher the empirical risk of a classifier, lower its optimal posterior weight. However, the rate at which these posterior weights decrease is influenced by the choice of divergence function in the PAC-Bayesian bound. In the case of KL-divergence based PAC-Bayesian bounds, optimal posterior weights decrease \textit{exponentially} with the empirical risk values. Whereas, the optimal posteriors that minimize $\chi^2$-divergence based PAC-Bayesian bounds have \textit{linearly} decreasing weights with respect to the empirical risk values. 

When the prior is uniform distribution, the optimal posterior weights in both the cases are directly proportional to the empirical risk values (no role for prior weight). Consider the case of linear distance function with $\chi^2$-divergence whose optimal posterior weights are determined in Theorem \ref{thm:OptQlinChi2.subset} given by Equation \eqref{eqn:OptQlinChi2.subset.unifP}. Similarly, when we have linear distance function with KL-divergence, the optimal posterior weights are given as (Equation (21) in \cite{pmlr-v101-sahu19a}):
\begin{equation} \label{eqn:optQ_linKL}
q^{\ast}_{i, \text{lin, KL}} = \frac{p_ie^{-m\hat{l}_i}}{\sum_{i=1}^{H}p_i e^{-m\hat{l}_i}} ~ \forall i = 1, \ldots, H 
\end{equation}
Using a uniform prior, $p_i = \dfrac{1}{H}$, in above, we get the following expression for the optimal posterior weights, which are directly proportional to the empirical risk values, $\hat{l}_i$s on exponential scale:
\begin{equation} \label{eqn:optQ_linKL}
q^{\ast}_{i, \text{lin, KL}} = \frac{e^{-m\hat{l}_i}}{\sum_{i=1}^{H}e^{-m\hat{l}_i}} ~ \forall i = 1, \ldots, H 
\end{equation}

\item \textbf{Size of support set}: KL-divergence based optimal posteriors take into account all the classifiers in the base set, though the posterior weight associated with a high risk classifier is infinitesimally small. On the other hand, $\chi^2$-divergence based optimal posteriors select only a \textit{strict} subset of base classifiers comprising of the ones with low empirical risk values.
\item \textbf{Test set performance}: Optimal posteriors for the KL-divergence based bounds have relatively lower test error rates than their $\chi^2$-divergence based counterparts. This hints at an underlying regularization phenomenon involving support set and tightness of the bound. The $\chi^2$-divergence based posteriors might be overfitting because they concentrate on a strict subset support of classifiers. In contrast, the KL-divergence based posteriors have the whole classifier set as their support and a better test set performance. This phenomenon is supported by the fact that KL-divergence based optimal posteriors yield tighter bounds than the ones derived for the case for $\chi^2$-divergence.
\end{enumerate}

\subsection{Comparison with cross-validation, (CV), method}
We performed 5-fold cross-validation (CV) on the datasets by setting aside 20\% of the data as a test set for each dataset. The set of $\lambda$ values is an arithmetic-geometric progression (AGP) with a logarithmic scale for $\lambda \in (0, 0.1)$ and a linear scale for $\lambda \geq 0.1$. This is the same set of $\lambda$ values which has been used for the proposed PAC-Bayesian technique. We report the test error of the ``best" $\lambda$ identified by the CV method and compare it with the PAC-Bayesian method using the (sq, $\chi^2$) pair (since this pair gives the lowest test error obtained by using different distance functions). These values are reported in Table \ref{tab:PACBvsCV} below. In terms of relative test error, the CV method is significantly better than the proposed method on Spambase, Bupa, Ionosphere and Haberman datasets, while the proposed PAC-Bayesian method is significantly better than CV method on Wdbc dataset. The difference in the two test errors is small (less than 20\%) on other datasets. Thus, CV method has better test error performance than $\chi^2$-divergence based PAC-Bayesian posterior. But the CV method takes remarkably longer time (between 2 to 10 hours) to identify the ``best" $\lambda$ than the time taken by the PAC-Bayesian method to identify the optimal posterior the classifier space (which generally takes about 20 to 300 seconds). Thus, computational complexity of the CV method is \textit{much higher} than that of PAC-Bayesian method.

\begin{table}[htb]
\centering
\setlength{\tabcolsep}{0.5mm}
\begin{tabular}{|c|c|c|c|c|c|c|}
\hline
\textbf{Dataset} & $\boldsymbol{\lambda^*}$ &  \textbf{\makecell{CV \\ Test Error}} & \textbf{\makecell{sq-$\chi^2$ PAC-B \\ Test Error}} & $\Delta$ Test Error & \makecell{Relative \\Test Error} & \textbf{\makecell{sq-$\chi^2$ PAC-B \\ Bound}}\\ \hline
\textbf{Spambase} & 0.15 &  0.13464 & 0.19437 & -0.05973 & -44.36\% \textcolor{blue}{$\star$} & 0.30581 \\ \hline
\textbf{Bupa} & 0.35 & 0.11594 & 0.29715 & -0.18121 & -156.29\% \textcolor{blue}{$\star$} & 0.60348 \\ \hline
\textbf{Mammographic} & 9.02 & 0.25904 & 0.23036 & 0.02868 & 11.07\%  & 0.38440 \\ \hline
\textbf{Wdbc} & 0.14 & 0.10434 & 0.07664 & 0.02770 & 26.55\% \textcolor{magenta}{$\bullet$} & 0.27465 \\ \hline
\textbf{Banknote} & 0.1 & 0.0 & 0.00661  & -0.00661  & NA & 0.13487 \\ \hline
\textbf{Mushroom} & 0.1 & 0.0 & 0.00092 & -0.00092 & NA & 0.06449\\ \hline
\textbf{Ionosphere} & 0.1 & 0.04225 & 0.06763 & -0.02538 & -60.07\% \textcolor{blue}{$\star$} & 0.38167 \\ \hline
\textbf{Waveform} & 2.31 & 0.06193 & 0.05122 & 0.01071 & 17.29\% & 0.12900 \\ \hline
\textbf{Haberman} & 3.16 & 0.22581 & 0.28979 & -0.06398 & -28.33\% \textcolor{blue}{$\star$} &0.51731 \\ \hline
\end{tabular}
\caption{Comparing the generalization performance of the simple cross-validation (CV) technique with the proposed PAC-Bayesian (PAC-B) technique on a set of SVMs obtained by using a logarithmic scale for $\lambda \in (0, 0.1)$ and a linear scale for $\lambda \geq 0.1$. $\Delta$ test error is the amount by which test error of CV method is smaller or larger than that of the PAC-Bayesian method for sq-$\chi^2$ pair. Relative test error is the ratio of $\Delta$ test error to the CV test error, signifying the relative difference between the test errors of the two methods. \textcolor{magenta}{$\bullet$} denotes the instances where the PAC-Bayesian method has a significantly lower test error rate than the CV method, while \textcolor{blue}{$\star$} denotes the instances where the CV method has a significantly lower test error rate. 
\label{tab:PACBvsCV}}
\end{table}

Apart from the \textit{computational benefits}, the proposed PAC-Bayesian method also has \textit{statistical advantages} over the CV method, as noted below:

	\begin{enumerate}[label = \roman*]
	\item \textbf{Sample robustness}: CV method is not as sample robust as the PAC-Bayesian method even though it trains on multiple sub-samples (partitions) and reports the averaged training error as CV error for choosing the best $\lambda$. For example, a 5-fold CV that we performed uses multiple training samples with a diminished sample size of 20\% of the dataset size. Whereas the PAC-Bayesian method uses a single and much larger training sample (60\% of dataset size in our scheme) to give an upper bound on the true risk.
	\item \textbf{Point estimate versus interval estimate}: CV method gives a point estimate of the true risk by averaging CV error over multiple folds, but there are no guarantees associated with it. On the other hand, the PAC-Bayesian method gives an interval estimate of the form $[0, Bnd]$, where $Bnd$ denotes the upper bound given by the PAC-Bayesian theorem and intrinsically has a high-probability guarantee associated.
	\item \textbf{Deterministic versus stochastic classifier}: CV method outputs a deterministic classifier in terms of the ``best" $\lambda$ value which has good test performance. The PAC-Bayesian technique is a committee method that outputs an optimal distribution on the set of classifiers, which yields a stochastic classifier. The classifier determined by the CV method may have better performance on a single test set, but the stochastic classifier obtained via PAC-Bayesian technique will have comparable performance when used on multiple test set instances.
	\item \textbf{An upper bound on true risk versus its point estimate}: The PAC-Bayesian method gives a tight upper bound on the true risk of the stochastic classifier for all datasets. This upper bounds holds even for the ``best" deterministic classifier obtained by CV method as can be seen in Table \ref{tab:PACBvsCV} above. Thus the high-probability PAC-Bayesian upper bound is a more useful quantity than the estimate given by CV test error for the true risk which can be under-biased or over-biased depending on the training sample and the folds created from this sample for cross-validation.
	\end{enumerate}
These advantages strengthen the usefulness of the PAC-Bayesian method for constituting a stochastic classifier and perhaps in tuning hyperparameters of other classification algorithms.

\section{Discussion \label{secn:conclusion}}
We determine optimal posteriors for PAC-Bayesian bound minimization problem with bounds derived using $\chi^2$-divergence function. The distance functions that we considered are: linear distance, squared distance (second degree polynomial) and KL-divergence (infinite degree polynomial). We first show that, in the uniform prior set up, minimizers of these PAC-Bayesian bounds can be obtained by a restricted search on subsets of the classifier set ordered by empirical risks. The bound minimization problem for linear distance case is shown to be a convex program and we also derive a closed form expression for its optimal posterior, while the other two distance functions result in non-convex programs. \added{We further show that the squared distance results in a quasi-convex bound under certain conditions, and it is computationally observed to have single local minimum.} We propose a convergent and computationally cheap fixed point based approach to identify the optimal posteriors for these bound minimization problems.

\added{Our computational exercise is comprehensive. The nine UCI datasets we have considered take into account small to moderate number of examples and features, balanced and imbalanced classes, and having different ranges and variance in the empirical risk values. Using this set of datasets helps us compare and understand the performance of optimal PAC-Bayesian posteriors due to different distance functions for a given divergence function, and also across KL- and $\chi^2$-divergence functions.}

Based on the computations on SVM classifiers, we observe that the squared distance based posteriors perform the best among the three  distance functions  in terms of bound values as well as average test error rates. The optimal posteriors for linear and squared distances have subset support, especially as the size of the classifier set increases. On the other hand KL-distance based posteriors usually have full support but do not perform well on the test set. This could be because they overfit the data while training. These chi-squared divergence based optimal posteriors do not have a high measure of  concentration, implying less bias towards classifiers with low empirical risks. 

Comparing with KL-divergence based optimal PAC-Bayesian posteriors, we observe that both groups of PAC-Bayesian posteriors have weights decreasing with respect to the empirical risk values of the classifiers. The difference lies in the rate of decrement -- KL-divergence based posterior weights decrease exponentially whereas $\chi^2$-divergence bases ones show a linear decrease.  Also, the former have a full support of the classifier set, while the later usually pick a strict subset of the base classifiers as their support. On the test set performance, KL-divergence based posteriors are better than those based on $\chi^2$-divergence. This phenomenon hints at an underlying regularization by KL-divergence based posteriors; perhaps, an implicit regularization.

We also provide a comparison of these PAC-Bayesian posteriors with the widely used cross-validation procedure as the baseline case. While the CV has lower test error sets, the PAC-Bayesian method has the advantages of sample robustness and a much lower computational cost over the cross-validation method. Also, it provides a reliable high probability upper bound on the true risk rather than a point estimate given by cross-validation method. 

\added{The significance of this work is in understanding the importance of choosing a divergence function for the PAC-Bayesian bound and its influence on the resulting optimal posteriors which are used to design stochastic classifiers. The challenges associated with this study are -- deducing the form of the bound for a given distance function, identifying the nature of the corresponding bound minimization problem, and obtaining a closed-form or a fixed point equation of the optimal posterior which minimizes this bound. Another challenge is identifying the support set of the optimal posterior for an arbitrary prior distribution on the classifier set. Therefore we considered the uniform prior which provided us with a structure to address the problem of identifying the support set in a linear fashion. On the computational side, we had to carefully choose the set of regularization parameter values to be used for generating the base SVMs, since this would greatly influence the performance of the stochastic SVMs built on them. To achieve a right mix of base classifiers with good risks, we considered an arithmetic-geometric progression of the regularization parameter values in the interval $[0, 5]$.} 

As a part of the future work, we can extend our results to a non-uniform prior on the classifier space, where we do not have such a structure to the feasible region and may need to do a full simplex search. We would also like to understand the nature of our results on high dimensional datasets such as image segmentation \cite{UCI:2017}, where each image is represented in matrix form rather than a vector or a record.


\bibliographystyle{plain}
\bibliography{NeurIPS2019bib}   

\appendix

\section{Optimal PAC-Bayesian Posterior using Linear Distance Function \label{appdx_secn:linChi2}} 
As a basic case, we can consider linear distance function, $\phi_{\text{lin}}(\hat{l}, l) = l - \hat{l}$ for $\hat{l}, l \in [0,1]$. The PAC-Bayesian bound in this case takes the following simplified form \cite{begin2016pac}:
\begin{equation} 
\mathbb{P}_S \left \lbrace \mathbb{E}_Q[l] \leq \mathbb{E}_Q [\hat{l}] + \sqrt{\frac{\chi^2(Q||P) + 1 }{4m\delta}} \right \rbrace \geq 1 - \delta.
\end{equation}
Thus, the upper bound on the true risk of a stochastic classifier governed by a distribution $Q$, when using the linear distance function with $\chi^2$-divergence between prior and posterior, is:
\begin{equation} \label{appdx_eqn:BlinChi2}
B_{\text{lin},\chi^2}(Q) = \sum_{i =1}^H \hat{l}_iq_i + \sqrt{\frac{\sum_{i =1}^H \frac{q_i^2}{p_i}}{4m\delta}}
\end{equation}

\subsection{The bound minimization problem \label{appdx_secn:BlinChi2OP}}
The corresponding bound optimization problem is:
\begin{equation} \label{appdx_eqn:Blinchi2OP}
\begin{split}
\min_{Q = \lbrace q_1, \ldots, q_H \rbrace} &\sum_{i = 1}^H \hat{l}_i q_i + \sqrt{\frac{\sum_{i =1}^H \frac{q_i^2}{p_i}}{4m\delta}} \\
\text{s. t.} \; & \sum_{i = 1}^H q_i = 1 \\
& q_i \geq 0 \quad \forall i = 1, \ldots, H.
\end{split}
\end{equation}
We are interested in the distribution $Q^{\ast}_{\text{lin}, \chi^2}$ which is optimal for the above bound minimization problem since that corresponds to the tightest PAC-Bayesian upper bound on the true risk of an stochastic classifier. 

\begin{theorem} \label{appdx_thm:convexity.BlinChi2}
The bound function $B_{\text{lin},\chi^2}(Q)= \sum_{i =1}^H \hat{l}_iq_i + \sqrt{\frac{\sum_{i =1}^H \frac{q_i^2}{p_i}}{4m\delta}}$ is a strictly convex function and hence the optimization problem \eqref{appdx_eqn:Blinchi2OP} is a convex program with a unique global minimum.
\end{theorem}
\begin{proof}$B_{\text{lin},\chi^2}(Q)$ is a differentiable function of $Q = \lbrace q_i \rbrace_{i =1}^H$. Hence we can prove its convexity if we can show that following first order condition holds for any $Q, Q'$:
\begin{align}
&B_{\text{lin},\chi^2}(Q') \geq B_{\text{lin},\chi^2}(Q) + \langle\nabla B_{\text{lin},\chi^2}(Q), Q'-Q\rangle \nonumber\\
\Rightarrow \; &\sum_{i = 1}^H \hat{l}_i q'_i + \sqrt{\frac{\sum_{i =1}^H \frac{q{'^2}_i}{p_i}}{4m\delta}}
\geq 
\sum_{i=1}^{H} \hat{l}_iq^{'}_i - \sum_{i=1}^{H} \hat{l}_iq_i + \frac{1}{\sqrt{4m\delta}} \frac{\sum_{i=1}^{H}\frac{q_iq'_i}{p_i} - \sum_{i=1}^{H}\frac{q^2_i}{p_i}}{\sqrt{\sum_{i =1}^H \frac{q_i^2}{p_i}}}  \nonumber\\
&\hspace{8cm} + \sum_{i = 1}^H \hat{l}_i q_i + \sqrt{\frac{\sum_{i =1}^H \frac{q_i^2}{p_i}}{4m\delta}} \nonumber\\
\Rightarrow \; &\sqrt{\frac{\sum_{i =1}^H \frac{q{'^2}_i}{p_i}}{4m\delta}} \geq \frac{1}{\sqrt{4m\delta}} \cdot \frac{\sum_{i=1}^{H}\frac{q_iq'_i}{p_i} - \sum_{i=1}^{H}\frac{q^2_i}{p_i} + \sum_{i=1}^{H}\frac{q^2_i}{p_i}}{\sqrt{\sum_{i =1}^H \frac{q_i^2}{p_i}}} \nonumber\\
\Rightarrow \; &\sqrt{\sum_{i =1}^H \frac{q{'^2}_i}{p_i}} \geq  \frac{\sum_{i=1}^{H}\frac{q_iq'_i}{p_i} }{\sqrt{\sum_{i =1}^H \frac{q_i^2}{p_i}}} \nonumber\\
\Rightarrow \; &\left(\sqrt{\sum_{i =1}^H \frac{q{'^2}_i}{p_i}} \right) \left( \sqrt{\sum_{i =1}^H \frac{q_i^2}{p_i}} \right) \geq  \sum_{i=1}^{H}\frac{q_iq'_i}{p_i} 
\end{align}
Using Cauchy-Schwarz inequality, the above holds for any pair of distributions $Q, Q'$ for a given prior distribution $P$ with equality if and only if $Q \equiv Q'$. This implies that the bound function $ B_{\text{lin}, \chi^2}(Q)$ is strictly convex. Therefore, the optimization problem \eqref{appdx_eqn:BlinChi2} has a unique global minimum.
\end{proof}

\subsection{The optimal posterior, $Q^{\ast}_{\text{lin}, \chi^2}$ via partial KKT system \label{appdx_secn:optQlinChi2}}
\begin{theorem} \label{appdx_thm:optQlinChi2.subset}
The global minimum of the bound minimization problem \eqref{appdx_eqn:Blinchi2OP} can be obtained via:
\begin{equation}
    q^{\ast}_{i, \text{lin}, \chi^2} = \frac{\left(\sum_{i = 1}^H \hat{l}_ip_i - \hat{l}_{i} \right)}{\sqrt{\frac{1}{4m\delta} - \hat{var}_P(\hat{l})}} \quad \forall i=1, \ldots, H 
\end{equation}
if the following two conditions are satisfied:
\begin{align}
q^{\ast}_{i, \text{lin}, \chi^2} &\geq 0 \quad \forall i=1, \ldots, H \\
\frac{1}{4m\delta} &> \hat{var}_P(\hat{l})
\end{align}
 where $\hat{var}_P(\hat{l}) := \sum_{i=1}^H p_i\left(\sum_{i = 1}^H \hat{l}_ip_i - \hat{l}_i \right)^2 $ is the variance of empirical risk values $\hat{l}_i$s under the prior distribution $P$.
\end{theorem}
\begin{proof}
The Lagrangian function corresponding to the optimization problem \eqref{appdx_eqn:Blinchi2OP} is:
\begin{equation}
\mathcal{L}_{\text{lin}, \chi^2}(Q, \mu) := \sum_{i = 1}^H\hat{l}_i q_i + \sqrt{\frac{\sum_{i =1}^H \frac{q_i^2}{p_i}}{4m\delta}} - \mu \left( \sum_{i = 1}^H q_i - 1 \right)
\end{equation}
At optimality, posterior $Q$ should set the derivatives of this Lagrangian function $\mathcal{L}_{\text{lin}, \chi^2}(Q, \mu)$ to zero. Setting the derivative of $\mathcal{L}_{\text{lin}, \chi^2}$ with respect to $q_i$'s as zero, we get:
\begin{align}
&\frac{\partial \mathcal{L}_{\text{lin}, \chi^2}(Q, \mu)}{\partial q_i} = 0 \nonumber\\
\Rightarrow \, 
&\hat{l}_i +  \frac{1}{\sqrt{4m\delta}} \frac{1}{\sqrt{\sum_{i =1}^H \frac{q_i^2}{p_i}}} \frac{q_i}{p_i} - \mu = 0 \nonumber\\
\Rightarrow \, &\frac{q_i}{\sqrt{\sum_{i =1}^H \frac{q_i^2}{p_i}}} = (\mu - \hat{l}_i)p_i\sqrt{4m\delta} \label{appdx_eqn:qLagrange_linChi2}
\end{align}
Based on the primal feasibility condition $\sum_{i = 1}^H q_i = 1$, we should have:
\begin{align}
&\frac{\partial \mathcal{L}_{\text{lin}, \chi^2}(Q, \mu)}{\partial \mu} = 0 \nonumber\\
& \sum_{i = 1}^H q_i = 1 \nonumber\\
\Rightarrow \, 
& \sqrt{4m\delta}\sqrt{\sum_{i =1}^H  \frac{q_i^2}{p_i}} \left[ \sum_{i = 1}^H (\mu - \hat{l}_i)p_i \right] = 1 \nonumber\\
\Rightarrow \, & \mu = \sum_{i = 1}^H \hat{l}_ip_i + \frac{1}{\sqrt{4m\delta}} \frac{1}{\sqrt{\sum_{i =1}^H \frac{q_i^2}{p_i}}} \label{appdx_eqn:Lagrangemultiplier_linChi2}
\end{align}

Thus, combining the results in \eqref{appdx_eqn:qLagrange_linChi2} and \eqref{appdx_eqn:Lagrangemultiplier_linChi2}, we get the following relation between $q_i$s and $p_i$s:
\begin{align*}
q_i = p_i \left[ 1 + \sqrt{4m\delta} \sqrt{\sum_{i =1}^H \frac{q_i^2}{p_i}} \left(\sum_{i = 1}^H \hat{l}_ip_i - \hat{l}_i \right) \right]\; \forall i = 1 \text{ to } H
\end{align*}
Using the transformation $z_i = \left(\frac{q_i}{p_i} - 1\right)^2$, the above can be reduced to a linear system of equations in $z_i$s:
\begin{align*}
&z_i = 4m\delta \left(\sum_{i = 1}^H \hat{l}_ip_i - \hat{l}_i \right)^2 \left( \sum_{i =1}^H z_ip_i + 1\right)  \; \forall i = 1 \text{ to } H \\
\Leftrightarrow \,& \sum_{i =1}^H z_ip_i - \frac{z_i}{ 4m\delta \left(\sum_{i = 1}^H \hat{l}_ip_i - \hat{l}_i \right)^2} + 1 = 0 \; \forall i = 1 \text{ to } H.
\end{align*}
The solution to this linear system is:
\begin{align}
& z_i^{\ast} = \frac{\left(\sum_{i = 1}^H \hat{l}_ip_i - \hat{l}_{i} \right)^2}{\frac{1}{4m\delta} - \hat{var}_P(\hat{l})} \quad i=1, \ldots, H \nonumber\\
\Rightarrow \; &\tilde{q}_i^{\ast} = \frac{\left(\sum_{i = 1}^H \hat{l}_ip_i - \hat{l}_{i} \right)}{\sqrt{\frac{1}{4m\delta} - \hat{var}_P(\hat{l})}} \quad i=1, \ldots, H \nonumber \\
\text{which implies } &q^{\ast}_{i, \text{lin}, \chi^2} = p_i \left[ 1  + \frac{\left(\sum_{i = 1}^H \hat{l}_ip_i - \hat{l}_{i} \right)}{\sqrt{\frac{1}{4m\delta} - \hat{var}_P(\hat{l})}} \right] \quad i=1, \ldots, H 
\end{align}
For $q^{\ast}_{i, \text{lin}, \chi^2}$ to be a real number, we need: $\frac{1}{4m\delta} - \hat{var}_P(\hat{l}) > 0$.
\end{proof}

\begin{rem}
We suspect that this upper bound on $\delta$ is related to the sparseness of the optimal posterior that minimizes the bound $B_{\text{lin}, \chi^2} (Q) = \sum_{i = 1}^H \hat{l}_i q_i + \sqrt{\frac{\sum_{i =1}^H \frac{q_i^2}{p_i}}{4m\delta}}$. A higher $\delta$ diminishes the effect of divergence $\chi^2[Q||P] = \sum_{i =1}^H \frac{q_i^2}{p_i}$. Hence it allows sparse solutions (where some components of posterior $Q$ take value zero) which have higher divergence from the prior compared to a non-sparse solution. 
\end{rem}

\subsection{Optimal posterior, $Q^{\ast}_{\text{lin}, \chi^2}$,  for uniform prior \label{appdx_secn:QlinChi2_unifP}}

\begin{theorem}[Optimal posterior on an ordered subset support] \label{appdx_thm:OptQlinChi2.subset.unifP}
When prior is uniform distribution on $\mathcal{H}$, among all the posteriors with support as subset of $\mathcal{H}$ of size exactly $H'$, the optimal/best posterior denoted by $Q^{\ast}_{\text{lin}, \chi^2}(H')$ has the support on the ordered subset $\mathcal{H'}_{\text{ord}} = \lbrace \lbrace \hat{l}_i \rbrace_{i=1}^{H'} \vert \hat{l}_1 \leq \hat{l}_2 \leq \ldots \leq \hat{l}_{H'}\rbrace$ consisting of smallest $H'$ values in $\mathcal{H}$. The optimal posterior weights are determined as follows:
\begin{equation}
    q^{\ast}_{i, \text{lin}, \chi^2} (H')=
    \begin{cases} \left[ 1 + \frac{\left( \frac{\sum_{i = 1}^{H'}\hat{l}_i}{H'} - \hat{l}_{i} \right)}{\sqrt{\frac{H}{H'4m\delta} - \hat{var}_{H'}(\hat{l})}} \right] \frac{1}{H} \quad &i=1, \ldots, H' \\
0 \quad &i= H'+1, \ldots, H,
    \end{cases}
\label{appdx_eqn:OptQlinChi2.subset.unifP}
\end{equation}
where $\hat{var}_{H'}(\hat{l})=  \frac{1}{H'}\sum\limits_{i =1}^{H'} \left( \frac{\sum_{i =1}^{H'} \hat{l}_i}{H'}- \hat{l}_i\right)^2 = \frac{\sum_{i =1}^{H'} \hat{l}^2_i}{H'} - \left(\frac{\sum_{i =1}^{H'} \hat{l}_i}{H'}\right)^2$ is the variance of the values in $\mathcal{H'}_{\text{ord}}$. We assume that the subset size $H'$ is such that $\frac{H}{H'4m\delta} - \hat{var}_{H'}(\hat{l}) > 0$ so that $Q^{\ast}_{\text{lin}, \chi^2}(H')$ is defined and for feasibility, we require $q^{\ast}_{i, \text{lin}, \chi^2}(H') > 0$ for $i=1, \ldots, H'$.
\end{theorem}
\begin{proof}
Under the uniform prior set, that is, when $p_i = \frac{1}{H}$ for all $i = 1, \ldots, H$, we obtain $Q^{\ast}_{\text{lin}, \chi^2}(H')$ using the partial KKT system (ignoring the positivity constraints) and the proof technique in Theorem \ref{appdx_thm:optQlinChi2.subset} above.
\end{proof}

\begin{figure}
    \centering
    \includegraphics[width = 0.7\textwidth]{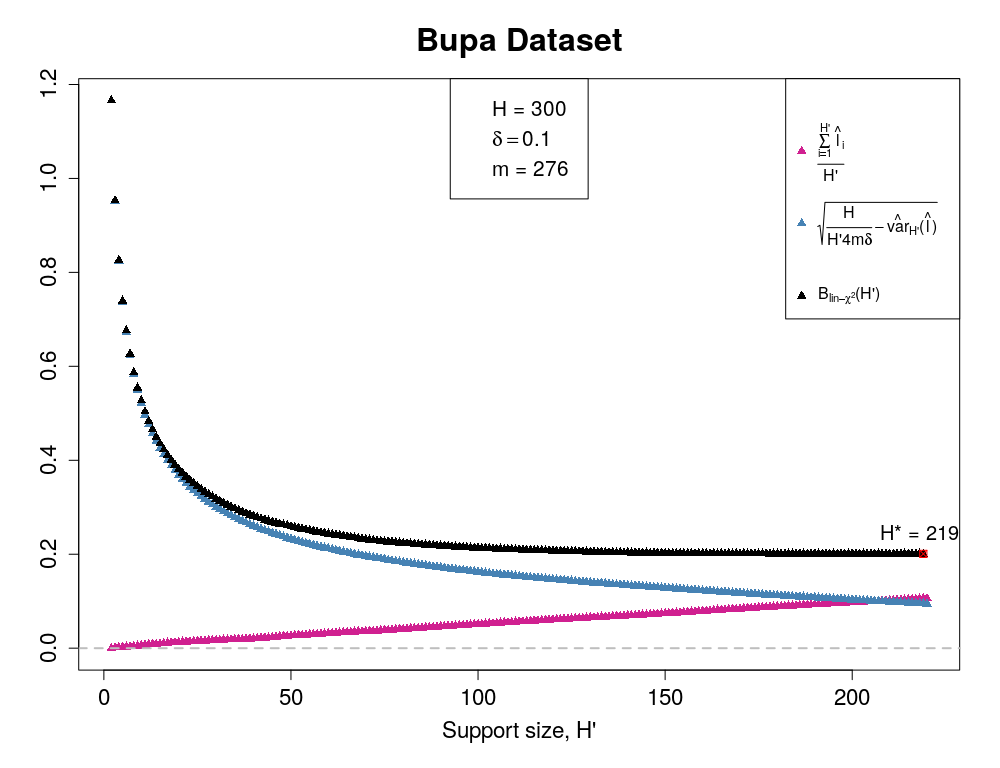}
    \caption{Nature of $B^{\ast}_{\text{lin}, \chi^2}(H')$ and its components as the subset size $H'$ varies}
    \label{appdx_fig:BlinChi2_Hdash}
\end{figure}

\begin{theorem} \label{appdx_thm:decreasing.optBlinChi2.Hdash}
The bound value for the linear distance function
\begin{equation}
B^{\ast}_{\text{lin}, \chi^2}(H') := \frac{\sum_{i =1}^{H'} \hat{l}_i}{H'} + \sqrt{\frac{H}{H'4m\delta} - \hat{var}_{H'}(\hat{l})} \label{appdx_eqn:optBlinChi2.Hdash}
\end{equation} 
of an optimal posterior $Q^{\ast}_{\text{lin}, \chi^2}(H')$ on an ordered subset of size $H'$, is decreasing function of $H'$ for $H' \leq H^{\ast}$. Here $H^{\ast}$ is the size of the ordered subset which forms the support of the globally optimal posterior $Q^{\ast}_{\text{lin}, \chi^2}$.
\end{theorem}
\begin{proof}
Consider an ordered subset of size $H' \in \lbrace 2, \ldots, H \rbrace$ such that $\frac{H}{H'4m\delta} - \hat{var}_{H'}(\hat{l}) > 0$ so that the posterior $Q^{\ast}_{\text{lin}, \chi^2}(H')$ given in \eqref{appdx_eqn:OptQlinChi2.subset.unifP} is defined. Suppose further that all the elements of this posterior are positive, so that it is feasible and hence optimal posterior on the considered ordered subset of size $H'$. The bound value at this optimal posterior can be computed as:
\begin{align}
B^{\ast}_{\text{lin}, \chi^2}(H') :&= B_{\text{lin}, \chi^2} \left( Q^{\ast}_{\text{lin}, \chi^2}(H') \right)= \sum_{i =1}^{H'}\hat{l}_iq^{\ast}_{i,\text{lin}, \chi^2}(H') + \sqrt{\frac{H}{4m\delta} \left(\sum_{i = 1}^{H'}\left(q^{\ast}_{i,\text{lin}, \chi^2}(H') \right)^2 \right)},
\end{align}
where
\begin{align}
\sum_{i =1}^{H'}\hat{l}_iq^{\ast}_{i,\text{lin}, \chi^2}(H') &= \sum_{i =1}^{H'}\hat{l}_i\left[ 1 + \frac{\left( \frac{\sum_{i = 1}^{H'}\hat{l}_i}{H'} - \hat{l}_{i} \right)}{\sqrt{\frac{H}{H'4m\delta} - \hat{var}_{H'}(\hat{l})}} \right] \frac{1}{H'} \nonumber \\
&=  \frac{\sum_{i =1}^{H'}\hat{l}_i}{H'} + \frac{\sum_{i =1}^{H'} \frac{\hat{l}_i}{H'}\left( \frac{\sum_{i = 1}^{H'}\hat{l}_i}{H'} - \hat{l}_{i} \right)}{\sqrt{\frac{H}{H'4m\delta} - \hat{var}_{H'}(\hat{l})}} \nonumber \\
&= \frac{\sum_{i =1}^{H'} \hat{l}_i}{H'} + \frac{ \left( \frac{\sum_{i = 1}^{H'}\hat{l}_i}{H'} \right)^2 - \frac{\sum_{i =1}^{H'}\hat{l}^2_i}{H'} }{\sqrt{\frac{H}{H'4m\delta} - \hat{var}_{H'}(\hat{l})}} \nonumber \\
&=  \frac{\sum_{i =1}^{H'}\hat{l}_i}{H'} - \frac{ \hat{var}_{H'}(\hat{l}) }{\sqrt{\frac{H}{H'4m\delta} - \hat{var}_{H'}(\hat{l})}}.
\end{align}
For evaluating the second term, we simplify the $\chi^2$-divergence factor:
\begin{align}
\sum_{i = 1}^{H'}\left(q^{\ast}_{i,\text{lin}, \chi^2}(H') \right)^2  &= \sum_{i=1}^{H'}\left[ 1 + \frac{\left( \frac{\sum_{i = 1}^{H'}\hat{l}_i}{H'} - \hat{l}_{i} \right)}{\sqrt{\frac{H}{H'4m\delta} - \hat{var}_{H'}(\hat{l})}} \right]^2 \frac{1}{(H')^2} \nonumber \\
&= \frac{1}{(H')^2} \sum_{i=1}^{H'}\left[ 1 + \frac{\left( \frac{\sum_{i = 1}^{H'}\hat{l}_i}{H'} - \hat{l}_{i} \right)^2}{\frac{H}{H'4m\delta} - \hat{var}_{H'}(\hat{l})} + \frac{2\left( \frac{\sum_{i = 1}^{H'}\hat{l}_i}{H'} - \hat{l}_{i} \right)}{\sqrt{\frac{H}{H'4m\delta} - \hat{var}_{H'}(\hat{l})}} \right]  \nonumber\\
&= \frac{1}{(H')^2} \left[ H' + \frac{\sum\limits_{i = 1}^{H'}\left( \frac{\sum_{i = 1}^{H'}\hat{l}_i}{H'} - \hat{l}_{i} \right)^2}{\frac{H}{H'4m\delta} - \hat{var}_{H'}(\hat{l})} + \frac{2 \sum\limits_{i = 1}^{H'} \left( \frac{\sum_{i = 1}^{H'}\hat{l}_i}{H'} - \hat{l}_{i} \right)}{\sqrt{\frac{H}{H'4m\delta} - \hat{var}_{H'}(\hat{l})}} \right] \nonumber \\
&= \frac{1}{H'} \left[ 1 + \frac{\frac{1}{H'}\sum\limits_{i = 1}^{H'}\left( \frac{\sum_{i = 1}^{H'}\hat{l}_i}{H'} - \hat{l}_{i} \right)^2}{\frac{H}{H'4m\delta} - \hat{var}_{H'}(\hat{l})} + \frac{\frac{2}{H'} \left( \cancel{\sum_{i = 1}^{H'}\hat{l}_i} - \cancel{\sum_{i = 1}^{H'}\hat{l}_{i}} \right)}{\sqrt{\frac{H}{H'4m\delta} - \hat{var}_{H'}(\hat{l})}} \right] \nonumber \\
&= \frac{1}{H'} \left[ 1 + \frac{\hat{var}_{H'}(\hat{l})}{\frac{H}{H'4m\delta} - \hat{var}_{H'}(\hat{l})} \right] \nonumber \\
&= \frac{1}{H'}\left(\frac{\frac{H}{H'4m\delta}}{\frac{H}{H'4m\delta} - \hat{var}_{H'}(\hat{l})}\right)
\end{align}
Substituting values of $q^{\ast}_{i,\text{lin}, \chi^2}(H')$ from Theorem \ref{appdx_thm:OptQlinChi2.subset.unifP}, the bound value becomes:
\begin{align}
B^{\ast}_{\text{lin}, \chi^2}(H') 
&=  \sum_{i =1}^{H'}\hat{l}_iq^{\ast}_{i,\text{lin}, \chi^2}(H') + \sqrt{\frac{H}{4m\delta} \left(\sum_{i = 1}^{H'}\left(q^{\ast}_{i,\text{lin}, \chi^2}(H') \right)^2 \right)} \nonumber \\
&=  \frac{\sum_{i =1}^{H'}\hat{l}_i}{H'} - \frac{ \hat{var}_{H'}(\hat{l}) }{\sqrt{\frac{H}{H'4m\delta} - \hat{var}_{H'}(\hat{l})}} + \sqrt{ \frac{H}{H'4m\delta} \left(\frac{\frac{H}{H'4m\delta}}{\frac{H}{H'4m\delta} - \hat{var}_{H'}(\hat{l})} \right)}  \nonumber \\
&=  \frac{\sum_{i =1}^{H'}\hat{l}_i}{H'} + \frac{\frac{H}{H'4m\delta} - \hat{var}_{H'}(\hat{l}) }{\sqrt{\frac{H}{H'4m\delta} - \hat{var}_{H'}(\hat{l})}}  \nonumber \\
&=  \frac{\sum_{i =1}^{H'}\hat{l}_i}{H'} + \sqrt{\frac{H}{H'4m\delta} - \hat{var}_{H'}(\hat{l})}.
\end{align}
This bound value $B^{\ast}_{\text{lin}, \chi^2}(H')$ is the sum of an increasing term, $ \frac{\sum_{i = 1}^{H'}\hat{l}_i}{H'}$, and a decreasing term, $\sqrt{\frac{H}{H'4m\delta} - \hat{var}_{H'}\left( \hat{l} \right)}$. $B^{\ast}_{\text{lin}, \chi^2}(H')$ can be shown to be a decreasing function of $H'$, as illustrated in Figure \ref{appdx_fig:BlinChi2_Hdash}.
\end{proof}

\begin{algorithm}[H]
\DontPrintSemicolon
\KwIn{$ m, \delta, H, \lbrace \hat{l}_i \rbrace_{i =1}^H$}
\KwOut{$Q^{\ast}_{\text{lin}, \chi^2}$}
flag $\gets 0$ \;
\For{$H' = 2, \ldots, H$}{
$\hat{var}_{H'}(\hat{l}) \gets \frac{1}{H'} \sum\limits_{i =1}^{H'} \left( \frac{\sum_{i =1}^{H'} \hat{l}_i}{H'}- \hat{l}_i\right)^2 $ \;
$DeltaBnd \gets \sqrt{\frac{H}{H'4m\delta} - \hat{var}_{H'}\left( \hat{l} \right)}$ \;
\If{$DeltaBnd < 0$ \label{step:undefinedQlinChi2}}{
flag $\gets 1$ \;
\KwBreak
}

\For{$i=1, \ldots, H'$}{
$    q^{\ast}_{i, \text{lin}, \chi^2}  \gets \left[ 1 + \frac{\left( \frac{\sum_{i = 1}^{H'}\hat{l}_i}{H'} - \hat{l}_{i} \right)}{DeltaBnd} \right] \frac{1}{H} $ \;
\If{$q^{\ast}_{i, \text{lin}, \chi^2} < 0$ \label{step:infeasibleQlinChi2}}{
flag $\gets 1$\;
\KwBreak
}
}
\If{flag $= 1$}{\KwBreak}
}
\eIf{flag $= 1$}{$H^{\ast} \gets H' - 1 $ \;
\For{$i=1, \ldots, H^{\ast}$}{
$    q^{\ast}_{i, \text{lin}, \chi^2}  \gets \left[ 1 + \frac{\left( \frac{\sum_{i = 1}^{H'}\hat{l}_i}{H'} - \hat{l}_{i} \right)}{DeltaBnd} \right] \frac{1}{H} $ \;
}
\For{$i=H^{\ast}+1, \ldots, H$}{
$  q^{\ast}_{i, \text{lin}, \chi^2} \gets 0$ \;
}
}{$H^{\ast} \gets H'$}

\Return $Q^{\ast}_{\text{lin}, \chi^2} \gets \left(q^{\ast}_{1, \text{lin}, \chi^2}, \ldots, q^{\ast}_{H, \text{lin}, \chi^2}\right)$ 
\caption{\textsc{OptQ lin-$\chi^2$ For Uniform Prior}: Algorithm for finding optimal posterior for the PAC-Bayesian bound with linear distance function and $\chi^2$-divergence when prior is uniform distribution \label{appdx_algo:optQsubset.unifP}}
\end{algorithm}

\subsubsection{Correctness of Algorithm \textsc{OptQ lin-$\chi^2$ For Uniform Prior}}
We want to determine the globally optimal posterior $Q^{\ast}_{\text{lin}, \chi^2}$ that has the minimum bound value $B_{\text{lin}, \chi^2}(Q)$ over the $H$-dimensional probability simplex, $\Delta^H$. Using the result of Theorem 1 in the main paper, we can confine the search to a much smaller space of posteriors with support  on a family of increasing ordered subsets of $\mathcal{H}$. These ordered subsets are defined by their size. For example, an ordered subset of size $H' \in \lbrace 1, \ldots, H \rbrace$ comprises of the \textit{lowest} $H'$ values in the set $\lbrace\hat{l}_i \rbrace_{i =1}^{H}$. Thus, the restricted space of posteriors, say $\Delta^{\text{ord}} \subset \Delta^H$, is a union of convex sets of posteriors with supports on the ordered subsets defined above. Due to increasing subset relation between consecutive supports, this union itself is a convex set. Therefore, as a consequence of Theorem \ref{appdx_thm:convexity.BlinChi2}, the bound function $B_{\text{lin}, \chi^2}(Q)$ is convex on the set of posteriors, $\Delta^{\text{ord}}$ as well, which contains the global minimum. The search space $\Delta^{\text{ord}}$ is a restriction of the simplex $\Delta^H$, yet consists of uncountably many posteriors on the ordered subsets. We refine the search further by localizing to optimal posteriors on each of the increasing ordered subsets and comparing their bound values to find the minimum.  As identified by \eqref{appdx_eqn:optBlinChi2.Hdash}, these bound values, $B^{\ast}_{\text{lin}, \chi^2}$ are functions of the subset size $H'$. $B^{\ast}_{\text{lin}, \chi^2}(H')$ is defined only for those values of $H'$ where $\sqrt{\frac{H}{H'4m\delta} - \hat{var}_{H'}\left( \hat{l} \right)} > 0$. We ignore those $H'$ values where this condition is not met (Line \ref{step:undefinedQlinChi2} in Algorithm \ref{appdx_algo:optQsubset.unifP}). Further, we also need to verify that for given $H'$, the optimal posterior $Q^{\ast}_{\text{lin}, \chi^2}(H')$ in \eqref{appdx_eqn:OptQlinChi2.subset.unifP} satisfies positivity constraints, as done in Line \ref{step:infeasibleQlinChi2} of the algorithm.  Hence an exponential search on restricted posterior space is simplified to a finite \textit{linear} search on the support size. We denote the support size of $Q^{\ast}_{\text{lin}, \chi^2}$ by $H^{\ast} \in [H]$. Therefore, for finding the optimal posterior $Q^{\ast}_{\text{lin}, \chi^2}$ in the restricted posterior space $\Delta^{\text{ord}}$, it is sufficient to search for $H^{\ast}$ in the set $\lbrace 1, \ldots, H \rbrace $ corresponding to support sizes.

\paragraph{Warm start for searching optimal support size, $H^{\ast}$}
We can reduce the sequential search for the optimal support size, $H^{\ast}$ on $\lbrace 1, \ldots, H \rbrace $ by using a warm start value for $H'$. Then we can reach $H^{\ast}$ by identifying a direction which leads to decrease in the bound value, $B^{\ast}_{\text{lin}, \chi^2}(H')$ as long as the corresponding posterior is defined and non-negative.  This bound value is the sum of an increasing term, $ \frac{\sum_{i = 1}^{H'}\hat{l}_i}{H'}$, and a decreasing term, $\sqrt{\frac{H}{H'4m\delta} - \hat{var}_{H'}\left( \hat{l} \right)}$. We expect $H^{\ast}$ to lie in the neighbourhood of the point of intersection of the these two components of $B^{\ast}_{\text{lin}, \chi^2}(H')$. The point of intersection can be obtained by equating the two terms and solving for $H'$:
\begin{align}
&\frac{\sum_{i = 1}^{H'}\hat{l}_i}{H'} = \sqrt{\frac{H}{H'4m\delta} - \hat{var}_{H'}\left( \hat{l} \right)} \nonumber \\
\Rightarrow \; &\sum_{i =1}^{H'} \hat{l}_i^2  = \frac{H}{4m\delta} \label{appdx_eqn:warmstart.HdashlinChi2}
\end{align}
The value of $H'$ which satisfies \eqref{appdx_eqn:warmstart.HdashlinChi2} can be used as a warm start point of the Algorithm \ref{appdx_algo:optQsubset.unifP}. If $H'$ satisfies feasibility conditions for $Q^{\ast}_{\text{lin}, \chi^2}(H')$, we compare its bound value, $B^{\ast}_{\text{lin}, \chi^2}(H')$ with its neighbours $B^{\ast}_{\text{lin}, \chi^2}(H' - 1)$ and $B^{\ast}_{\text{lin}, \chi^2}(H' + 1)$ to identify a descent direction until feasibility is violated. Otherwise, if $B^{\ast}_{\text{lin}, \chi^2}(H')$ is infeasible or undefined, we keep decrementing $H'$ till we reach a feasible $Q^{\ast}_{\text{lin}, \chi^2}(H')$.

\section{Optimal PAC-Bayesian Posterior using Squared Distance Function \label{appdx_secn:sqChi2}}
The distance function of our interest is the squared distance function:
$\phi_{\text{sq}} \left(\hat{l}, l \right) =  \left(\hat{l} - l \right)^2$ for $\hat{l}, l \in [0,1]$. The PAC-Bayesian bound for squared distance function with chi-squared divergence can be stated as:
\begin{equation}
\mathbb{P}_S \left \lbrace \left(\mathbb{E}_Q [\hat{l}] - \mathbb{E}_Q[l] \right)^2 \leq \sqrt{\left[ \chi^2(Q||P) + 1 \right] \left(\frac{\mathcal{I}^{R}_{{\text{sq}}}(m,2)}{\delta} \right)} \right \rbrace \geq 1 - \delta. \label{appdx_eqn:PACB_sqChi}
\end{equation}
The above statement gives the following probabilistic upper bound on the true risk of an stochastic classifier governed by a distribution $Q$ on $\mathcal{H}$:
\begin{equation}
B_{sq, \chi^2}(Q) := \mathbb{E}_Q [\hat{l}] + \sqrt[4]{\left[ \chi^2(Q||P) + 1 \right] \left(\frac{\mathcal{I}^{R}_{{\text{sq}}}(m,2)}{\delta} \right)}
\end{equation}
We first to need to identify the constant $\mathcal{I}^{R}_{\text{sq}}(m, 2)$ for a given sample size $m$.

\begin{lem} \label{appdx_lem:I_R_sq}
For a given sample size, $m$, $l^\ast = 0.5$ is the maximizer of $\mathcal{I}^{R}_{\text{sq}}(m, 2, l) := \sum_{k=0}^{m}{m \choose k} l^k (1-l)^{m-k}e^{2m\left(\frac{k}{m}-l\right)^4} $ for $l \in [0, 1]$.
\end{lem}
\begin{proof} We have
\begin{align*}
\mathcal{I}^{R}_{\text{sq}}(m, 2) 
&= \sup_{l \in [0, 1]} \left[ \sum_{k = 0}^m \binom{m}{k} l^k (1- l)^{m - k} \left(\frac{k}{m} - l \right)^{4}\right], \\
&= \sup_{l \in [0, 1]} \frac{1}{m^4}\left[ \sum_{k = 0}^m \binom{m}{k} l^k (1- l)^{m - k} \left(k - ml \right)^{4}\right].
\end{align*}
The quantity to be maximized in the above expression is the fourth central moment of a Binomial distribution, and can be computed in terms of its first and second central moments:
\begin{align*}
\mathcal{I}^{R}_{\text{sq}}(m, 2) 
&= \sup_{l \in [0, 1]} \frac{1}{m^4}\left \lbrace m [l (1 - l)^4 + l^4(1-l) ] + 3m(m-1)l(1-l) ]\right \rbrace \\
&= \frac{1}{m^3} \sup_{l \in [0, 1]} l(1-l)\left \lbrace [(1-l)^3 + l^3] + 3(m-1)]\right \rbrace \\
&= \frac{1}{m^3} \sup_{l \in [0, 1]} l(1-l)\left \lbrace (1-l+l) [(1-l)^2 + l^2 - l(1-l)] + 3(m-1)]\right \rbrace \\
&= \frac{1}{m^3} \sup_{l \in [0, 1]} l(1-l)\left[3l^2 -3l + 1 + 3(m-1)\right] \\
&= \frac{1}{m^3} \sup_{l \in [0, 1]} \underbrace{\left[-3l^4 + 6l^3 - 4l^2 + l + 3(m-1)(l-l^2) \right]}_{\mathcal{I}^{R}_{\text{sq}}(m, 2, l)} \\
\text{Thus, }\mathcal{I}^{R}_{\text{sq}}(m, 2) &= \frac{1}{m^3} \sup_{l \in [0,1]} \mathcal{I}^{R}_{\text{sq}}(m, 2, l).
\end{align*} 
$\mathcal{I}^{R}_{\text{sq}}(m, 2, l)$ is a smooth, continuous function of $l \in [0, 1]$ and the maximum can be obtained via derivative test. $\mathcal{I}^{R}_{\text{sq}}(m, 2, l)$ is a concave function of $l \in [0,1]$. We observe that the unique maximum is attained at $l^{\ast} = \frac{1}{2}$ for all $m \geq 2$ (See Figure \ref{appdx_fig:I_R_sq}).
\begin{figure}
\includegraphics[scale=0.45]{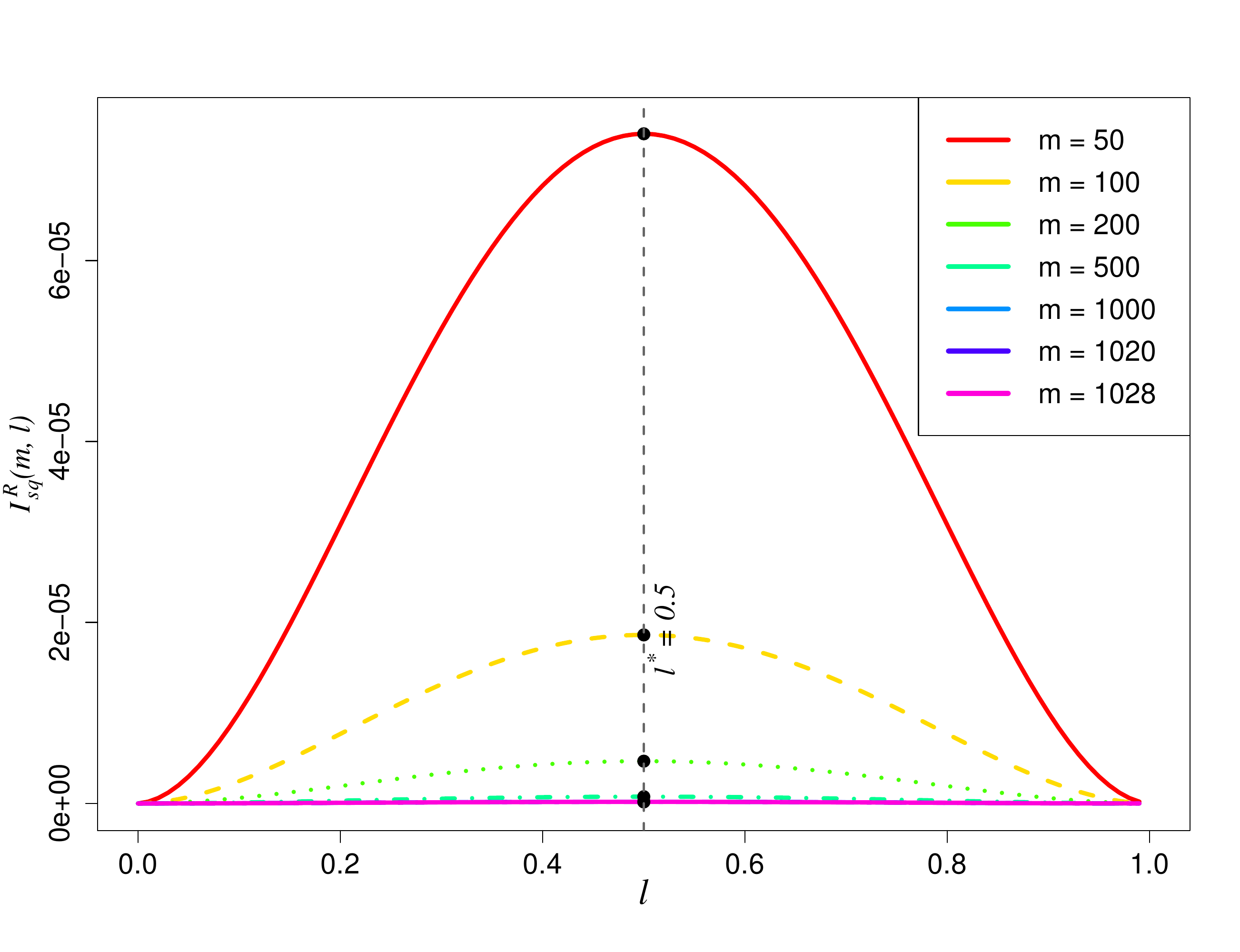}
\caption{Plot of the function $\mathcal{I}^{R}_{\text{sq}}(m,2,l) := \sup_{l \in [0, 1]} \sum_{k = 0}^m \binom{m}{k} l^k (1- l)^{m - k} \left(\frac{k}{m} - l \right)^{4}$ for different sample sizes, $m$. $\mathcal{I}^{R}_{\text{sq}}(m,2,l)$ is symmetric about $l = 0.5$, its maximizer. \label{appdx_fig:I_R_sq}}
\end{figure}

Thus,
\begin{align*}
\mathcal{I}^{R}_{\text{sq}}(m, 2) &= \frac{1}{m^3} \mathcal{I}^{R}_{\text{sq}}(m, 2, \frac{1}{2})= \frac{12m -11}{16m^3}
\end{align*}
Hence, proved.
\end{proof}

Thus, the PAC-Bayesian bound for squared distance function with chi-squared divergence can be stated as: 
\begin{equation}
\mathbb{P}_S \left \lbrace \left(\mathbb{E}_Q [\hat{l}] - \mathbb{E}_Q[l] \right)^2 \leq \sqrt{\left[ \chi^2(Q||P) + 1 \right] \left(\frac{12m -11}{16m^3\delta} \right)} \right \rbrace \geq 1 - \delta. \label{appdx_eqn:PACB_sqChi}
\end{equation}

\begin{theorem}
For a finite set of classifiers, $\mathcal{H}$, PAC-Bayesian upper bound on the averaged true risk based on squared distance function when chi-squared divergence is used as a measure of divergence between the prior and the posterior is given by:
\begin{equation}
  B_{sq, \chi^2}(Q) = \sum_{i = 1}^H \hat{l}_i q_i + \sqrt[4]{\left( \sum_{i = 1}^H \frac{q_i^2}{p_i} \right) \left(\frac{12m -11}{16m^3\delta} \right)}
\end{equation}
\end{theorem}
\begin{proof}
Using the PAC-Bayesian statement in \eqref{appdx_eqn:PACB_sqChi} for the case of a finite classifier set, we can obtain the above form of $B_{sq, \chi^2}(Q)$, where the expression $\frac{12m -11}{16m^3\delta} $ has been identified via Lemma \ref{appdx_lem:I_R_sq}.
\end{proof}

\subsection{The bound minimization problem}
We want to determine the optimal posterior $Q^{\ast}_{sq, \chi^2}$ which minimizes the upper bound, $B_{sq, \chi^2}(Q)$. When classifier space $\mathcal{H}$ is a finite set, say $\mathcal{H} = \lbrace h_i\rbrace_{i =1}^H$, this optimization problem can be described as:
\begin{equation} \label{appdx_eqn:Bsqchi2OP}
\begin{split}
\min_{q_1, \ldots, q_H} &\sum_{i = 1}^H \hat{l}_i q_i + \sqrt[4]{\left( \sum_{i = 1}^H \frac{q_i^2}{p_i} \right) \left(\frac{12m -11}{16m^3\delta} \right)} \\
\text{s. t.} \; & \sum_{i = 1}^H q_i = 1 \\
& q_i \geq 0 \quad \forall i = 1, \ldots, H.
\end{split}
\end{equation}

\subsection{Non-convexity of the bound function \label{appdx_secn:nonconvexity.BsqChi2}}
If the bound function, $B_{sq, \chi^2}(Q)$ turns out to be convex, then it has a unique minimizer which can be easily obtained using the KKT conditions.  We investigate whether this bound function is convex in $Q$ using the first order conditions for convexity.
\begin{theorem}
The bound function, $B_{sq, \chi^2}(Q) = \sum_{i = 1}^H \hat{l}_i q_i + \sqrt[4]{\left( \sum_{i = 1}^H \frac{q_i^2}{p_i} \right) \left(\frac{12m -11}{16m^3\delta} \right)} $ is non-convex.
\end{theorem}
\begin{proof}
We use the first order condition to verify convexity of our bound function. We need to check if the following condition holds for any pair of distributions $Q$ and $Q'$ on classsifier space $\mathcal{H}$:
\begin{align}
&B_{\text{sq}, \chi^2}(Q') \geq B_{\text{sq}, \chi^2}(Q) + \langle\nabla B_{\text{sq}, \chi^2}(Q), Q'-Q\rangle \nonumber\\
\Rightarrow \; &\sum_{i = 1}^H \hat{l}_i q'_i + \sqrt[4]{\left( \sum_{i = 1}^H \frac{q'^{2}_i}{p_i} \right) \left(\frac{12m -11}{16m^3\delta} \right)} 
\geq 
\sum_{i=1}^{H} \hat{l}_iq'_i - \sum_{i=1}^{H} \hat{l}_iq_i  \nonumber\\
&\hspace{1cm} + \sqrt[4]{\frac{12m - 11}{16m^3 \delta}} \cdot \frac{\left( \sum_{i = 1}^H\frac{q_iq'_i}{p_i} - \sum_{i = 1}^H\frac{q_i^2}{p_i}\right)}{2\left(\sum_{i = 1}^H \frac{q_i^2}{p_i} \right)^{3/4}} 
+ \sum_{i=1}^{H} \hat{l}_iq_i + \sqrt[4]{\left( \sum_{i = 1}^H \frac{q_i^2}{p_i} \right) \left(\frac{12m -11}{16m^3\delta} \right)} \nonumber \\
\Rightarrow \; & \left(\sum_{i = 1}^H \frac{q_i^2}{p_i} \right)^{3/4} \left( \sum_{i = 1}^H \frac{q'^{2}_i}{p_i} \right)^{1/4}
\geq 
\left(\sum_{i = 1}^H\frac{q_iq'_i}{p_i} + \sum_{i = 1}^H \frac{q_i^2}{p_i} \right) \bigg/2 \label{appdx_eqn:1ConvexCondn_BsqChi2}
\end{align}
Notice that this inequality does not depend on $\hat{l}_i$ values. We have counter examples which violate this convexity condition. Consider $H = 10$ and $P$ to be uniform distribution on $\mathcal{H}$. If $Q$ is a degenerate distribution with $q_8 = 1$ and
\begin{multline*}
Q' = ( 0.1538802, 0.1199569, 0.04226614, 0.06115894, 0.06160916, 0.07520679, \\ 0.1450413, 0.2345929, 0.01762696, 0.08866069),
\end{multline*}
then we have LHS =  6.087086 and RHS =  6.172964 for \eqref{appdx_eqn:1ConvexCondn_BsqChi2} violating the convexity property. Thus, we can claim that $B_{\text{sq}, \chi^2}$ is a non-convex function.
\end{proof}

\begin{rem} Computationally this bound minimization problem is observed to have single local minimum. The quasi-convexity of this bound function is holds under a condition identified in Propostion \ref{appdx_propn:quasiconvex.condn.BsqChi2}.
\end{rem}

We are interested in checking whether $B_{\text{sq},\chi^2}(Q)$ is strictly quasi-convex. If so, we can  claim that a local optimal solution will be a global optimal solution \cite{bazaraa2013nonlinear}.
\begin{defn} \cite{bazaraa2013nonlinear}
Let $f: E \longrightarrow \mathbb{R}$ where $E$ is a non-empty convex set in $\mathbb{R}^n$. A function $f$ is strictly quasi convex if, for each $\mathbf{x}_1,\mathbf{x}_2 \in E$ with $f(\mathbf{x}_1) \neq f(\mathbf{x}_2)$, we have 
\begin{equation}
    f[\alpha \mathbf{x}_1 + (1-\alpha)\mathbf{x}_2] < \max(f(\mathbf{x}_1),f(\mathbf{x}_2)) \quad \forall \alpha \in (0,1). \label{appdx_defn:quasiconvex}
\end{equation} 
\end{defn}
\begin{theorem} \cite{bazaraa2013nonlinear}
Let $f: E \longrightarrow \mathbb{R}$ be strictly quasi-convex. Consider the problem of minimizing $f(\mathbf{x})$ subject to $\mathbf{x} \in E$, where $E$ is a non-empty convex set in $\mathbb{R}^n$. If $\bar{\mathbf{x}}$ is a local optimal solution, then $\bar{\mathbf{x}}$ is also a global optimal solution. \label{appdx_thm:quasiconvex_localglobal}
\end{theorem}

\begin{propn} \label{appdx_propn:quasiconvex.condn.BsqChi2}
The bound function $B_{\text{sq},\chi^2}(Q)$ is strictly quasi-convex if the following condition holds for any $Q, Q'$ for each $\alpha \in (0,1)$:
\begin{multline*}
 \bigg(\sqrt[4]{\frac{12m-11}{16m^3 \delta}}\bigg)\left[ \sqrt[4]{\sum_{i=1}^{H} \frac{(\alpha q_i + (1-\alpha)q'_i)^2}{p_i}} - \sqrt[4]{\sum_{i=1}^{H} \frac{q_i^2}{p_i}} \right] < (1-\alpha)(E_{Q} [\hat{l}]-E_{Q'} [\hat{l}])
\end{multline*}
and hence a local minimum to the bound minimization problem \eqref{appdx_eqn:Bsqchi2OP} is also a global minimum.
\end{propn}
\begin{proof}
$B_{\text{sq},\chi^2}(Q)$ is defined on the simplex $\Delta^H$ which is a non-empty convex set in $\mathbb{R}^H$. For quasiconvexity, we need to show that for each $Q \neq Q' \in \Delta^H$ with $B_{\text{sq},\chi^2}(Q) \neq B_{\text{sq},\chi^2}(Q') $, the following holds: 
\begin{equation*}
    B_{\text{sq},\chi^2}[\alpha Q + (1-\alpha)Q'] < \max(B_{\text{sq},\chi^2}(Q),B_{\text{sq},\chi^2}(Q')) \quad \forall \alpha \in (0,1).
\end{equation*}
That is equivalent to showing:
\begin{multline}
E_{\alpha Q + (1-\alpha)Q'} [\hat{l}] + \bigg(\sqrt[4]{\frac{12m-11}{16m^3 \delta}}\bigg)\sqrt[4]{\sum\limits_{i=1}^{H} \frac{(\alpha q_i + (1-\alpha)q'_i)^2}{p_i}}  \nonumber \\
<  \max \left \lbrace E_{Q} [\hat{l}] + \bigg(\sqrt[4]{\frac{12m-11}{16m^3 \delta}}\bigg)\sqrt[4]{\sum_{i=1}^{H} \frac{q_i^2}{p_i}}, E_{Q'} [\hat{l}] + \bigg(\sqrt[4]{\frac{12m-11}{16m^3 \delta}}\bigg)\sqrt[4]{\sum_{i=1}^{H} \frac{q{'^2}_i}{p_i}}\right \rbrace
\end{multline}
We assume that  $B_{\text{sq},\chi^2}(Q) > B_{\text{sq},\chi^2}(Q')$. This implies that we need to show that $ B_{\text{sq},\chi^2}(\alpha Q + (1-\alpha)Q') < B_{\text{sq},\chi^2}(Q)$. We consider 4 cases as follows:
\begin{enumerate}
\item[Case I]: $E_{Q}[\hat{l}] = E_{Q'}[\hat{l}]$ and $\sum\limits_{i=1}^{H} \frac{q^2_i}{p_i} = \sum\limits_{i=1}^{H}\frac{q{'^2}_i}{p_i}$
 then we have
$$ E_{\alpha Q + (1-\alpha)Q'}[\hat{l}] = \alpha E_{Q}[\hat{l}] +(1-\alpha)E_{Q'}[\hat{l}] = E_{Q}[\hat{l}] = E_{Q'}[\hat{l}].$$
Thus to show that $B_{\text{sq},\chi^2}(\alpha Q + (1-\alpha)Q') < B_{\text{sq},\chi^2}(Q)$, we have to show the following for any $Q, Q'$ for each $\alpha \in (0, 1)$:
\begin{multline*}
    E_{\alpha Q + (1-\alpha)Q'}[\hat{l}] + \sqrt[4]{\frac{12m-11}{16m^3 \delta}} \sqrt[4]{\sum_{i=1}^{H} \frac{(\alpha q_i + (1-\alpha)q'_i)^2}{p_i}} \\
    < E_{Q} [\hat{l}] + \bigg(\sqrt[4]{\frac{12m-11}{16m^3 \delta}}\bigg)\sqrt[4]{\sum_{i=1}^{H} \frac{q_i^2}{p_i}}
\end{multline*}
This is equivalent to showing that for any $Q, Q'$ for each $\alpha \in (0, 1)$,
\begin{equation}
    \sqrt[4]{\sum_{i=1}^{H} \frac{(\alpha q_i + (1-\alpha)q'_i)^2}{p_i}} < \sqrt[4]{\sum_{i=1}^{H} \frac{q_i^2}{p_i}}
\end{equation}

Consider the LHS as given below:
\begin{eqnarray}\nonumber
\sum_{i=1}^{H} \frac{(\alpha q_i + (1-\alpha)q'_i)^2}{p_i} & = & \sum_{i=1}^{H} \frac{(\alpha^2q_i^2 + (1-\alpha)^2 q{'^2}_i + 2\alpha(1-\alpha)q_iq'_i)}{p_i} \\ \nonumber 
&<& \sum_{i=1}^{H} \frac{(\alpha^2q_i^2 + (1-\alpha)^2 q{'^2}_i + 2\alpha(1-\alpha)q_i^2)}{p_i}, \\
&& \left(\text{since} \sum_{i=1}^{H} q_iq'_i < \sum_{i=1}^{H} q_i^2 \forall Q \neq Q' \right) \nonumber\\ \nonumber
&=& \sum_{i=1}^{H} \frac{q_i^2 - (1-\alpha)^2 q_i^2 + (1-\alpha)^2 q{'^2}_i)}{p_i}\\ \nonumber
&=& \sum_{i=1}^{H} \frac{q_i^2}{p_i} + (1-\alpha)^2 \left[\sum_{i=1}^{H}\frac{q{'^2}_i}{p_i} - \sum_{i=1}^{H}\frac{q_i^2}{p_i}\right] \nonumber \\
&<& \sum_{i=1}^{H}\frac{q_i^2}{p_i}
\end{eqnarray}
Since $\sqrt[4]{x}$ is an increasing function of $x$, we have the proof of strict quasi-convexity for this case.

\item[Case II]: $\sum\limits_{i=1}^{H}\frac{q_i^2}{p_i} = \sum_{i=1}^{H}\frac{q{'^2}_i}{p_i}$ and $E_{Q}[\hat{l}] > E_{Q'}[\hat{l}]$
 then we have, $$E_{\alpha Q + (1-\alpha)Q'}[\hat{l}] = \alpha E_{Q}[\hat{l}] +(1-\alpha)E_{Q'}[\hat{l}]  <  E_{Q}[\hat{l}]$$  
By previous argument, 
\begin{eqnarray}\nonumber
\sum_{i=1}^{H} \frac{(\alpha q_i + (1-\alpha)q'_i)^2}{p_i} & < & \sum_{i=1}^{H} \frac{q_i^2}{p_i} + (1-\alpha)^2 \left[\sum\limits_{i=1}^{H}\frac{q{'^2}_i}{p_i} - \sum\limits_{i=1}^{H}\frac{q_i^2}{p_i}\right] 
\end{eqnarray}
As the second term on the RHS of above inequality is zero by assumption, we get
\begin{multline*}
    E_{\alpha Q + (1-\alpha)Q'}[\hat{l}] + \sqrt[4]{\frac{12m-11}{16m^3 \delta}} \sqrt[4]{\sum_{i=1}^{H} \frac{(\alpha q_i + (1-\alpha)q'_i)^2}{p_i}} \\
    < E_{Q} [\hat{l}] + \bigg(\sqrt[4]{\frac{12m-11}{16m^3 \delta}}\bigg)\sqrt[4]{\sum_{i=1}^{H} \frac{q_i^2}{p_i}}
\end{multline*}
This implies that
\begin{equation*}
    \Leftrightarrow B_{\text{sq},\chi^2}[\alpha Q + (1-\alpha)Q'] < B_{\text{sq},\chi^2}(Q) =  \max(B_{\text{sq},\chi^2}(Q),B_{\text{sq},\chi^2}(Q')) 
\end{equation*}
Therefore, $B_{\text{sq},\chi^2}(Q)$ is strictly quasi-convex in this case too.

\item[Case III]: $E_{Q}[\hat{l}] > E_{Q'}[\hat{l}]$ and $\sum_{i=1}^{H}\frac{q_i^2}{p_i} > \sum_{i=1}^{H}\frac{q{'^2}_i}{p_i}$
, then, as earlier we have,
\begin{equation*}
    E_{\alpha Q + (1-\alpha)Q'}[\hat{l}] < E_{Q}[\hat{l}]
\end{equation*}
And,
\begin{eqnarray}\nonumber
\sum_{i=1}^{H} \frac{(\alpha q_i + (1-\alpha)q'_i)^2}{p_i}  <  \sum_{i=1}^{H} \frac{q_i^2}{p_i} + (1-\alpha)^2 \left[\sum_{i=1}^{H}\frac{q{'^2}_i}{p_i} - \sum_{i=1}^{H}\frac{q_i^2}{p_i}\right]  < \sum_{i=1}^{H} \frac{q_i^2}{p_i}
\end{eqnarray}
The last inequality follows since we have assumed \mbox{$\sum_{i=1}^{H}\frac{q{'^2}_i}{p_i} - \sum_{i=1}^{H}\frac{q_i^2}{p_i} < 0$.} Hence, it follows that $B_{\text{sq},\chi^2}[\alpha Q + (1-\alpha)Q'] < B_{\text{sq},\chi^2}(Q) =  \max\{B_{\text{sq},\chi^2}(Q),B_{\text{sq},\chi^2}(Q')\}$.
Therefore, $B_{\text{sq},\chi^2}(Q)$ is strictly quasi-convex in this case.

The condition for quasiconvexity  holds easily in all the above three cases. The next case requires an added assumption.

\item[Case IV]: $E_{Q}[\hat{l}] > E_{Q'}[\hat{l}]$ and $\sum_{i=1}^{H}\frac{q_i^2}{p_i} < \sum_{i=1}^{H}\frac{q{'^2}_i}{p_i}$ with $B_{\text{sq},\chi^2}(Q)>B_{\text{sq},\chi^2}(Q')$
then we get the following inequality,
\begin{eqnarray}\nonumber
E_{Q} [\hat{l}] + \bigg(\sqrt[4]{\frac{12m-11}{16m^3 \delta}}\bigg)\sqrt[4]{\sum_{i=1}^{H} \frac{q_i^2}{p_i}} &> & E_{Q'} [\hat{l}] + \bigg(\sqrt[4]{\frac{12m-11}{16m^3 \delta}}\bigg)\sqrt[4]{\sum_{i=1}^{H} \frac{q{'^2}_i}{p_i}} \\ \nonumber
\Longleftrightarrow 
(E_Q[\hat{l}]-E_{Q'}[\hat{l}]) &>& \sqrt[4]{\frac{12m-11}{16m^3 \delta}} \bigg[\sqrt[4]{\sum_{i=1}^{H} \frac{q{'^2}_i}{p_i}}-\sqrt[4]{\sum_{i=1}^{H} \frac{q_i^2}{p_i}} \bigg]
\end{eqnarray}
Hence, we have to show that, $ B_{\text{sq},\chi^2}(\alpha Q + (1-\alpha)Q') < B_{\text{sq},\chi^2}(Q)$
\begin{multline*}
    E_{\alpha Q + (1-\alpha)Q'}[\hat{l}] + \sqrt[4]{\frac{12m-11}{16m^3 \delta}} \sqrt[4]{\sum_{i=1}^{H} \frac{(\alpha q_i + (1-\alpha)q'_i)^2}{p_i}} \\
    < E_{Q} [\hat{l}] + \bigg(\sqrt[4]{\frac{12m-11}{16m^3 \delta}}\bigg)\sqrt[4]{\sum_{i=1}^{H} \frac{q_i^2}{p_i}}
\end{multline*}
This is equivalent to proving that
\begin{multline*}
 \bigg(\sqrt[4]{\frac{12m-11}{16m^3 \delta}}\bigg)\left[ \sqrt[4]{\sum_{i=1}^{H} \frac{(\alpha q_i + (1-\alpha)q'_i)^2}{p_i}} - \sqrt[4]{\sum_{i=1}^{H} \frac{q_i^2}{p_i}} \right] \\
 < E_{Q} [\hat{l}]- (\alpha E_{Q} [\hat{l}] + (1-\alpha)E_{Q'} [\hat{l}]) 
\end{multline*}
That is, we need to show that for any $Q, Q'$ for each $\alpha \in (0, 1)$:
\begin{multline*}
 \bigg(\sqrt[4]{\frac{12m-11}{16m^3 \delta}}\bigg)\left[ \sqrt[4]{\sum_{i=1}^{H} \frac{(\alpha q_i + (1-\alpha)q'_i)^2}{p_i}} - \sqrt[4]{\sum_{i=1}^{H} \frac{q_i^2}{p_i}} \right] \\
< (1-\alpha)(E_{Q} [\hat{l}]-E_{Q'} [\hat{l}])
\end{multline*}
The above holds due to the assumption in the theorem statement.
\end{enumerate}
Thus, under the given condition, $B_{\text{sq}, \chi^2}$ is strictly quasi-convex and admits a global minimum which can be identified based on KKT conditions.
\end{proof} 

\begin{rem}
The condition that for any $Q, Q'$ for each $\alpha \in (0, 1)$:
\begin{multline*}
 \bigg(\sqrt[4]{\frac{12m-11}{16m^3 \delta}}\bigg)\left[ \sqrt[4]{\sum_{i=1}^{H} \frac{(\alpha q_i + (1-\alpha)q'_i)^2}{p_i}} - \sqrt[4]{\sum_{i=1}^{H} \frac{q_i^2}{p_i}} \right] < (1-\alpha)(E_{Q} [\hat{l}]-E_{Q'} [\hat{l}])
\end{multline*}
is required to complete the proof of quasi-convexity of $B_{\text{sq}, \chi^2}(Q)$ for the case when $E_{Q}[\hat{l}] > E_{Q'}[\hat{l}]$ and $\sum_{i=1}^{H}\frac{q_i^2}{p_i} < \sum_{i=1}^{H}\frac{q{'^2}_i}{p_i}$. 
We haven't been able to verify that this condition will always hold for any pair $(Q, Q')$.
Other cases are easy to prove.
\end{rem}

\subsection{The posterior based on fixed point scheme, $Q^{\ast}_{\text{sq}, \chi^2}$ \label{appdx_secn:optQsqChi2}}
The Lagrangian function corresponding to the optimization problem \eqref{appdx_eqn:Bsqchi2OP} is:
\begin{equation}
\mathcal{L}_{sq, \chi^2}(Q, \mu) := \hat{l}_i q_i + \sqrt[4]{\left( \sum_{i = 1}^H \frac{q_i^2}{p_i} \right) \left(\frac{12m -11}{16m^3\delta} \right)} - \mu \left( \sum_{i = 1}^H q_i - 1 \right)
\end{equation}
At optimality, posterior $Q$ should set the derivatives of this Lagrangian function $\mathcal{L}_{sq, \chi^2}(Q, \mu)$ to zero. 
Setting the derivative of $\mathcal{L}_{sq, \chi^2}$ with respect to $q_i$'s as zero,
We get:
\begin{align}
&\frac{\partial \mathcal{L}_{sq, \chi^2}(Q, \mu)}{\partial \, q_i} = 0  \quad \forall i = 1, \ldots, H \nonumber\\
\Rightarrow &\hat{l}_i + \left(\sqrt[4]{\frac{12m -11}{16m^3\delta}} \times \frac{1}{4 \left( \sum_{i = 1}^H  \frac{q_i^2}{p_i}\right)^{\frac{3}{4}}} \times \frac{2q_i}{p_i} \right) - \mu = 0 \nonumber\\
\Rightarrow &
\frac{q_i}{\left( \sum_{i = 1}^H  \frac{q_i^2}{p_i}\right)^{\frac{3}{4}}} = \frac{2p_i(\mu - \hat{l}_i)}{\sqrt[4]{\frac{12m -11}{16m^3\delta}}}  \quad \forall i = 1, \ldots, H \label{appdx_eqn:qLagrange_sqChi2}
\end{align}
And now, setting the derivative of $\mathcal{L}_{sq, \chi^2}$ with respect to $\mu$ as zero, we get:
\begin{align}
&\frac{\partial \mathcal{L}_{sq, \chi^2}(Q, \mu)}{\partial \, \mu} = 0  \nonumber\\
\Rightarrow &\sum_{i = 1}^H q_i - 1 = 0 \nonumber\\
\Rightarrow &
\sum_{i = 1}^H \left( \left[\frac{2\left( \sum_{i = 1}^H  \frac{q_i^2}{p_i}\right)^{\frac{3}{4}}}{\sqrt[4]{\frac{12m -11}{16m^3\delta}} } \right](\mu - \hat{l}_i)p_i \right)  = 1 \label{eqn:Lagrangemultiplier_sqChi2}
\\
\Rightarrow &\mu = \sum_{i = 1}^H \hat{l}_i p_i + \frac{\sqrt[4]{\frac{12m -11}{16m^3\delta}}}{2\left( \sum_{i = 1}^H  \frac{q_i^2}{p_i}\right)^{\frac{3}{4}}}
\end{align}
We get the following fixed point equation in $q_i$'s:
\begin{align}
q_i =  \left[ \frac{2 \left( \sum_{i = 1}^H  \frac{q_i^2}{p_i}\right)^{\frac{3}{4}}}{\sqrt[4]{\frac{12m -11}{16m^3\delta}}} \left(\sum_{i = 1}^H \hat{l}_i p_i - \hat{l}_i \right) + 1 \right] p_i \; \forall i = 1, \ldots, H
\end{align}

\begin{theorem}[Optimal posterior on an ordered subset support] \label{appdx_thm:OptQsqChi2.subset}
When prior is uniform distribution on $\mathcal{H}$, among all the posteriors with support as subset of size exactly $H'$, the best posterior denoted by $Q^{\ast}_{\text{sq}, \chi^2}(H')$ has the support on the ordered subset $\mathcal{H'}_{\text{ord}} = \lbrace \lbrace \hat{l}_i \rbrace_{i=1}^{H'} \vert \hat{l}_1 \leq \hat{l}_2 \leq \ldots \leq \hat{l}_{H'}\rbrace$ consisting of smallest $H'$ values in $\mathcal{H}$. The optimal posterior weights are determined as the solution to the following fixed point equation:
\begin{equation}
    q^{FP}_{i, \text{sq}, \chi^2} (H')=
    \begin{cases} \frac{1}{H'} + \frac{2 \left( \sum_{i = 1}^{H'}  (q^{FP}_{i, \text{sq}, \chi^2} (H'))^2 \right)^{\frac{3}{4}}}{\sqrt[4]{\frac{(12m -11)H}{16m^3\delta}}} \left(\frac{\sum_{i = 1}^{H'} \hat{l}_i}{H'} - \hat{l}_i \right)  \quad &i=1, \ldots, H' \\
0 \quad &i= H'+1, \ldots, H.
    \end{cases}
\label{appdx_eqn:OptQsqChi2.subset}
\end{equation}
under the assumption that for a given $H'$, \eqref{appdx_eqn:OptQsqChi2.subset} converges to a fixed point solution and for feasibility, we require $q^{FP}_{i, \text{sq}, \chi^2}(H') > 0$ for $i=1, \ldots, H'$.
\end{theorem}

\section{Optimal PAC-Bayesian Posterior using KL-distance } \label{appdx_secn:klChi2}
The PAC-Bayesian bound using the distance function $kl(\hat{l} , l) = \hat{l} \ln \left( \frac{\hat{l}}{l}\right) + (1 - \hat{l}) \ln \left( \frac{1 - \hat{l}}{1 - l}\right)$ (for any $\hat{l}, l \in (0, 1)$) is obtained as:
\begin{equation}
\mathbb{P}_S \left \lbrace \forall Q \text{ on } \mathcal{H}: \;kl\left(\mathbb{E}_Q [\hat{l}], \mathbb{E}_Q[l] \right) \leq \sqrt{\frac{(\chi^2(Q||P) + 1) \mathcal{I}^{R}_{\text{kl}}(m, 2)}{\delta}} \right \rbrace \geq 1 - \delta.
\end{equation}

The upper bound on the averaged true risk $\mathbb{E}_Q[l]$ corresponding to the above PAC-Bayesian theorem is obtained as:
\begin{equation}
B_{\text{kl}, \chi^2}(Q) = \sup_{r \in (0, 1)} \left\lbrace r : kl\left(\mathbb{E}_Q [\hat{l}], r \right) \leq \sqrt{\frac{(\chi^2(Q||P) + 1) \mathcal{I}^{R}_{\text{kl}}(m, 2)}{\delta}} \right\rbrace
\end{equation}

An inverse $kl( \cdot,  \cdot)$ function does not exist since it is not a monotone function, and so the bound $B_{\text{kl}, \chi^2}(Q) $ does not have an explicit form. 
However, we can employ a numerical root finding algorithm such as that described in \cite{PACBIntervals} (Algo. (\textsc{KLroots})) to obtain $B_{\text{kl}, \chi^2}(Q) $ for a given instance of system parameters.

We first need to compute the constant $\mathcal{I}^{R}_{\text{kl}}(m, 2) := \sum\limits_{k=0}^{m}{m \choose k} l^k (1-l)^{m-k} \left(kl\left(\frac{k}{m},l\right)\right)^2$ in order to determine the bound value.

\begin{figure}
\centering
\includegraphics[width = 0.7\textwidth]{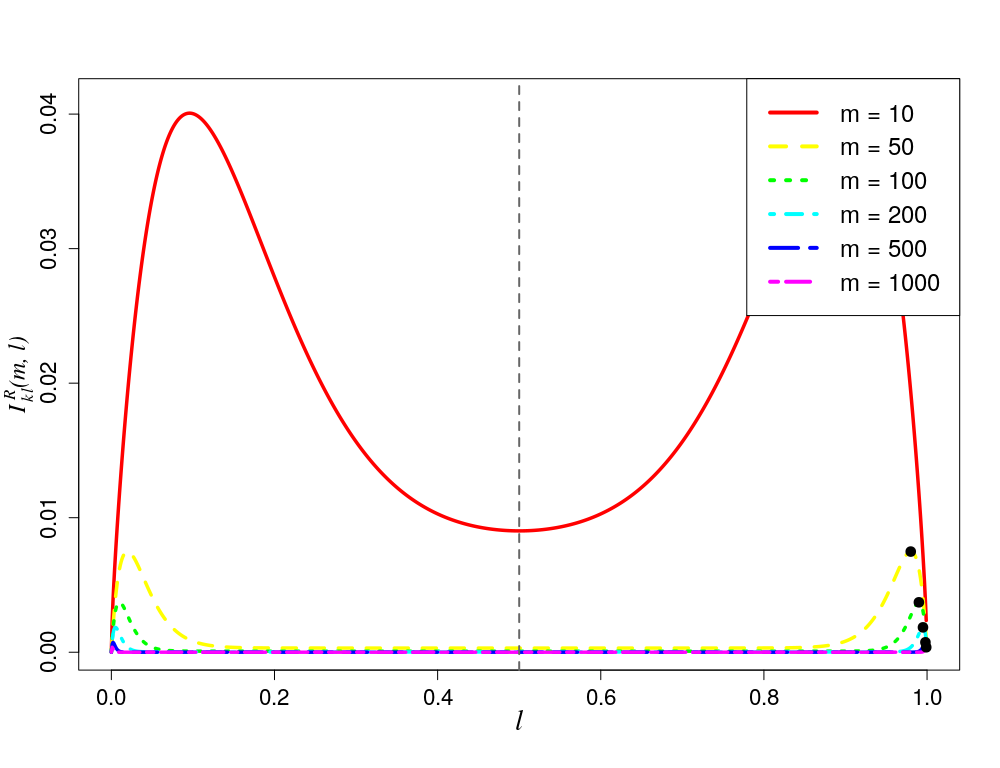}
\caption[Plot of the function $\mathcal{I}^R_{\text{kl}} (m, l)$ for PAC-Bayesian bound with squared distance function and KL divergence measure]{Plot of the function $\mathcal{I}^R_{\text{kl}} (m, l) = \sum\limits_{k=0}^{m}{m \choose k} l^k (1-l)^{m-k} \left(kl\left(\frac{k}{m},l\right)\right)^2$ as a function of the true risk $l \in [0, 1]$ for different values of the sample size, $m$ represented by different curves in the above graph. We observe that the function $\mathcal{I}^R_{\text{kl}} (m, l)$ is bimodal and symmetric about $l = 0.5$. We are interested in the quantity $\mathcal{I}^R_{\text{kl}} (m) = \sup_{l \in [0,1]} \mathcal{I}^R_{\text{kl}} (m, l)$ as a function of $m$ which we identify graphically (and mark it by a $\bullet$ on each curve). \label{appdx_fig:I_R_kl}}
\end{figure}

\begin{table}[ht]
\centering
\caption{Table of $\mathcal{I}^R_{\text{kl}} (m) = \sup_{l \in [0,1]}\sum\limits_{k=0}^{m}{m \choose k} l^k (1-l)^{m-k} \left(kl\left(\frac{k}{m},l\right)\right)^2$ values for different sample sizes, $m$. We notice that $\mathcal{I}^R_{\text{kl}} (m)$ decreases as $m$ increases. For a given $m$,$l^{\ast}(m)$ denotes the value of $l \in [0,1]$ at which the supremum is attained. We observe that $l^{\ast}(m) \rightarrow 1$ as $m$ grows beyond 1000. \label{appdx_tab:I_R_kl}}
\begin{tabular}{r r r }
 \toprule
\makecell{Sample \\ size, $m$} & $l^{\ast} (m) $ & $I^R_\text{kl}(m)$ \\ 
\midrule
50 &  0.98 & 0.0074799 \\ 
  100 &  0.99 & 0.0037092 \\ 
  200 & 0.995 & 0.0018470 \\ 
  500 & 0.998 & 0.0007369 \\ 
  1000 & 0.999 & 0.0003682 \\ 
  1020 & 0.999 & 0.0003609 \\ 
  1028 & 0.999 & 0.0003580 \\ 
\bottomrule
\end{tabular}
\end{table}
For $m > 1028$, computation is difficult due to storage limitations in the range of floating point numbers -- gives  $I^R_\text{kl}(m)$ as NaN. We notice that $\mathcal{I}^R_{\text{kl}} (m)$ decreases with $m$ and hence, we can use $I^R_\text{kl}(1028)$ as an upper approximation for $I^R_\text{kl}(m)$ for $m > 1028$.

\subsection{The KL-distance bound minimization problem}
For a finite classifier space $\mathcal{H} = \lbrace h_i\rbrace_{i =1}^H$, this optimization problem can be described as:
\begin{subequations}  \label{appdx_eqn:BklChi2OP}
\begin{align}
&\min_{q_1, \ldots, q_H, r} r\\
\hspace{-4cm}\text{s.t.} \quad & \left(\sum_{i = 1}^H \hat{l}_iq_i \right) \ln\left( \frac{\sum\limits_{i = 1}^H \hat{l}_iq_i}{r}\right) + \left(1 - \sum_{i = 1}^H \hat{l}_iq_i \right) \ln\left( \frac{ 1  - \sum\limits_{i = 1}^H \hat{l}_iq_i}{1 - r}\right) 
=\sqrt{\frac{\left(\sum_{i=1}^H \frac{q_i^2}{p_i}\right) \mathcal{I}^{R}_{\text{kl}}(m, 2)}{\delta}} \label{appdx_eqn:klChi2DCcons} \\
& r \geq \sum_{i = 1}^H \hat{l}_iq_i  \label{appdx_eqn:r>EQl}\\
& \sum_{i = 1}^H q_i = 1  \label{appdx_eqn:sumQ}\\
& q_i \geq 0 \; \forall i = 1, \ldots, H  \label{appdx_eqn:qi>0}
\end{align} 
\end{subequations}
Here, $r$ is the right root of $kl\left(\mathbb{E}_Q [\hat{l}], r \right) = \sqrt{\frac{(\chi^2(Q||P) + 1) \mathcal{I}^{R}_{\text{kl}}(m, 2)}{\delta}}$ for a given $\mathbb{E}_Q [\hat{l}]$. The above is known to be a non-convex problem with a difference of convex (DC) equality constraint \eqref{appdx_eqn:klChi2DCcons}. The constraint \eqref{appdx_eqn:r>EQl} is a strict inequality which is relaxed for modelling purpose.

\subsection{The posterior based on fixed point scheme, $Q^{\text{FP}}_{\text{kl}, \chi^2}$ \label{appdx_secn:optQklChi2}}
We derive FP equation for KL-distance based bound optimization problem below:
\begin{theorem}[Optimal posterior on an ordered subset support] \label{appdx_thm:qFP_klChi2}
Among all the posteriors with support as subset of size exactly $H'$, a stationary point $Q^{\text{FP}}_{\text{kl}, \chi^2}(H')$ can be obtained as the solution to the following fixed point equation:
\begin{equation}
q_i = \frac{1}{Z_{\text{kl}, \chi^2}} p_i \left(\sum_{i=1}^{H'} \frac{q_i^2}{p_i} \right) \left \lbrace 1 + \frac{\left(\hat{l}_i - \sum_{i=1}^{H'}\hat{l}_iq_i \right)}{\sqrt{\frac{\left(\sum_{i=1}^{H'} \frac{q_i^2}{p_i}\right) \mathcal{I}^{R}_{\text{kl}}(m, 2)}{\delta}} } \left[ \ln\left( \frac{(1 - r) \sum_{i = 1}^{H'} \hat{l}_i q_i}{r(1-\sum_{i = 1}^{H'} \hat{l}_i q_i)}\right) \right] \right \rbrace 
\label{appdx_eqn:qFP_klChi2}
\end{equation}
where $Z_{\text{kl}, \chi^2}$ is a suitable normalization constant and $r$ is the solution to \eqref{appdx_eqn:klChi2DCcons} and \eqref{appdx_eqn:r>EQl} for a given $Q = (q_1, \ldots, q_H)$.
\end{theorem}
\begin{proof}
The Lagrangian function for \eqref{appdx_eqn:BklChi2OP} can be written as follows:
\begin{multline}
\mathcal{L}_{\text{kl}, \chi^2} = r - \beta_0 \left[ \left(\sum_{i = 1}^{H'} \hat{l}_iq_i \right) \ln\left( \frac{\sum_{i = 1}^{H'} \hat{l}_iq_i}{r}\right) + \left(1 - \sum_{i = 1}^{H'} \hat{l}_iq_i \right) \ln\left( \frac{ 1  - \sum_{i = 1}^{H'} \hat{l}_iq_i}{1 - r}\right) \right.\\
\left. - \sqrt{\frac{\left(\sum_{i=1}^{H'} \frac{q_i^2}{p_i}\right) \mathcal{I}^{R}_{\text{kl}}(m, 2)}{\delta}} \right] -\beta_1 \left( r - \sum_{i = 1}^{H'} \hat{l}_iq_i \right) - \mu_0 \left(\sum_{i=1}^{H'} q_i -1 \right) - \sum_{i = 1}^{H'} \mu_i q_i
\end{multline}
Due to the strict inequality constraint \eqref{appdx_eqn:r>EQl}, complementary slackness conditions for a stationary point imply that the Lagrange multiplier $\beta_1$ should vanish at optimality.

Differentiating $\mathcal{L}_{\text{kl}, \chi^2}$ with respect to primal variables $r$ and $q_i$s, and also with respect to dual variable $\mu_0$, we get:
\begin{align}
\frac{\partial \mathcal{L}_{\text{kl}, \chi^2}}{\partial r} &= 1- \beta_0\left[- \left(\frac{\sum_{i=1}^{H}\hat{l}_iq_i}{r}\right) + \left( \frac{1-\sum_{i = 1}^{H'} \hat{l}_i q_i}{1-r}\right)\right] \label{appdx_eqn:deriver_klChi2}\\
\frac{\partial \mathcal{L}_{\text{kl}, \chi^2}}{\partial q_i} &= -\beta_0 \left[ \hat{l}_i \ln\left( \frac{\sum_{i = 1}^{H'} \hat{l}_i q_i}{r}\right) + \cancel{\hat{l}_i} - \hat{l}_i \ln \left( \frac{1-\sum_{i = 1}^{H'} \hat{l}_i q_i}{1-r}\right) - \cancel{\hat{l}_i} - \sqrt{\frac{\mathcal{I}^{R}_{\text{kl}}(m, 2)}{\delta \left(\sum_{i=1}^{H'} \frac{q_i^2}{p_i}\right)}} \cdot \frac{q_i}{p_i}\right] \nonumber \\
&\hspace{8cm}  - \mu_0 -\mu_i \quad \forall i = 1, \ldots, H \label{appdx_eqn:deriveq_klChi2}\\
\frac{\partial \mathcal{L}_{\text{kl}, \chi^2}}{\partial \mu_0} &= \sum_{i=1}^{H}q_i -1 \label{appdx_eqn:derivmu0_klChi2}
\end{align}
At an optimal solution, these derivatives should be set to zero. Let us first consider the derivative \eqref{appdx_eqn:deriver_klChi2} and set it to zero. That is,
\begin{align}
&\dfrac{\partial \mathcal{L}_{\text{kl}, \chi^2}}{\partial r} = 0 \nonumber\\
\Rightarrow &1- \beta_0\left[- \left(\frac{\sum_{i=1}^{H}\hat{l}_iq_i}{r}\right) + \left( \frac{1-\sum_{i = 1}^{H'} \hat{l}_i q_i}{1-r}\right)\right] = 0 \nonumber\\[3mm]
\Rightarrow & \beta_0\left[ \frac{-\sum\limits_{i=1}^{H}\hat{l}_iq_i + \cancel{r \left(\sum\limits_{i=1}^{H}q_i\hat{l}_i \right)}  + r - \cancel{r \left(\sum\limits_{i=1}^{H}q_i\hat{l}_i \right)} }{r(1-r)} \right] =1 \nonumber\\
\Rightarrow &\beta_0 = \frac{r(1-r)}{r - \sum\limits_{i=1}^{H}\hat{l}_iq_i} >0 \label{appdx_eqn:beta0LagrangeklChi2}
\end{align}
The denominator in above is strictly positive since $r>\sum\limits_{i=1}^{H}q_i\hat{l}_i$ . The inequality constraint in \eqref{appdx_eqn:BklChi2OP} also implies that $r \in (0, 1)$, which means that the numerator term is also strictly positive. Hence, we have $\beta_0 > 0$ which is a feasible value for the Lagrange parameter.

Next consider the derivative \eqref{appdx_eqn:deriveq_klChi2} of the Lagrange $\mathcal{L}_{\text{kl}, \chi^2}$. We multiply it with $q_i$ and set it zero to get:
\begin{align}
&q_i\frac{\partial \mathcal{L}_{\text{kl}, \chi^2}}{\partial q_i} = 0 \quad \forall i = 1, \ldots, H \nonumber\\
\Rightarrow & -\beta_0 \left\{ \hat{l}_iq_i \left[ \ln\left( \frac{\sum_{i = 1}^{H'} \hat{l}_i q_i}{r}\right) - \ln \left( \frac{1-\sum_{i = 1}^{H'} \hat{l}_i q_i}{1-r}\right)\right] - \sqrt{\frac{\mathcal{I}^{R}_{\text{kl}}(m, 2)}{\delta \left(\sum_{i=1}^{H'} \frac{q_i^2}{p_i}\right)}} \cdot \frac{q^2_i}{p_i} \right\}  - \mu_0q_i -\mu_iq_i =0 \label{appdx_eqn:qmultiplyLagrangederive}
\end{align}
where $\mu_iq_i = 0$ due to complementary slackness conditions, since $\mu_i$ is the Lagrange multiplier for the constraint $q_i \geq 0$ in \eqref{appdx_eqn:BklChi2OP}. Since we are interested in finding the best posterior on the ordered subset of size $H' \leq H$, only first $H'$ values of the  distribution $Q = (q_1, \ldots, q_H)$ will take strictly positive values. Therefore summing \eqref{appdx_eqn:qmultiplyLagrangederive} over $i = 1, \ldots, H'$, we get:
\begin{multline*}
\sum_{i = 1}^{H'}q_i \frac{\partial \mathcal{L}_{\text{kl}, \chi^2}}{\partial q_i} = -\beta_0 \left\{\sum_{i = 1}^{H'} \hat{l}_iq_i \left[ \ln\left( \frac{\sum_{i = 1}^{H'} \hat{l}_i q_i}{r}\right) - \ln \left( \frac{1-\sum_{i = 1}^{H'} \hat{l}_i q_i}{1-r}\right)\right]  \right. \\
\left. - \sqrt{\frac{\mathcal{I}^{R}_{\text{kl}}(m, 2)}{\delta \left(\sum_{i=1}^{H'} \frac{q_i^2}{p_i}\right)}} \cdot \sum_{i = 1}^{H'}\frac{q^2_i}{p_i} \right\}  - \mu_0 \left( \sum_{i =1}^{H'}q_i \right) -\sum_{i =1}^{H'}\mu_iq_i =0 
\end{multline*}
Since $\sum_{i =1}^{H'}q_i = 1$, we get:
\begin{multline}
\mu_0 = -\beta_0 \left\{\sum_{i = 1}^{H'} \hat{l}_iq_i \left[ \ln\left( \frac{\sum_{i = 1}^{H'} \hat{l}_i q_i}{r}\right) - \ln \left( \frac{1-\sum_{i = 1}^{H'} \hat{l}_i q_i}{1-r}\right)\right] - \sqrt{\frac{\mathcal{I}^{R}_{\text{kl}}(m, 2)}{\delta \left(\sum_{i=1}^{H'} \frac{q_i^2}{p_i}\right)}} \cdot \sum_{i = 1}^{H'}\frac{q^2_i}{p_i} \right\} \label{appdx_eqn:mu0LagrangeklChi2}
\end{multline}
Then using \eqref{appdx_eqn:deriveq_klChi2} and above \eqref{appdx_eqn:mu0LagrangeklChi2}, we get:
\begin{multline}
\cancel{-\beta_0} \left[ \hat{l}_i \ln\left( \frac{\sum_{i = 1}^{H'} \hat{l}_i q_i}{r}\right) - \hat{l}_i \ln \left( \frac{1-\sum_{i = 1}^{H'} \hat{l}_i q_i}{1-r}\right) - \sqrt{\frac{\mathcal{I}^{R}_{\text{kl}}(m, 2)}{\delta \left(\sum_{i=1}^{H'} \frac{q_i^2}{p_i}\right)}} \cdot \frac{q_i}{p_i}\right] \nonumber \\
= \cancel{-\beta_0} \left\{\sum_{i = 1}^{H'} \hat{l}_iq_i \left[ \ln\left( \frac{\sum_{i = 1}^{H'} \hat{l}_i q_i}{r}\right) - \ln \left( \frac{1-\sum_{i = 1}^{H'} \hat{l}_i q_i}{1-r}\right)\right] - \sqrt{\frac{\mathcal{I}^{R}_{\text{kl}}(m, 2)}{\delta \left(\sum_{i=1}^{H'} \frac{q_i^2}{p_i}\right)}} \cdot \sum_{i = 1}^{H'}\frac{q^2_i}{p_i} \right\} 
\end{multline}
\begin{align}
\Rightarrow \; &\sqrt{\frac{\mathcal{I}^{R}_{\text{kl}}(m, 2)}{\delta \left(\sum_{i=1}^{H'} \frac{q_i^2}{p_i}\right)}} \cdot \left( \frac{q_i}{p_i} - \sum_{i = 1}^{H'}\frac{q^2_i}{p_i} \right) =
 \left( \hat{l}_i - \sum_{i = 1}^{H'} \hat{l}_iq_i \right) \left[ \ln\left( \frac{\sum_{i = 1}^{H'} \hat{l}_i q_i}{r}\right) - \ln \left( \frac{1-\sum_{i = 1}^{H'} \hat{l}_i q_i}{1-r}\right)\right] \nonumber \\
\Rightarrow \;  &q_i = p_i \times \left[ \sum_{i = 1}^{H'}\frac{q^2_i}{p_i} + \frac{\left( \hat{l}_i - \sum_{i = 1}^{H'} \hat{l}_iq_i \right)}{\sqrt{\frac{\mathcal{I}^{R}_{\text{kl}}(m, 2)}{\delta \left(\sum_{i=1}^{H'} \frac{q_i^2}{p_i}\right)}}}
\left[ \ln\left( \frac{\sum_{i = 1}^{H'} \hat{l}_i q_i}{r}\right) - \ln \left( \frac{1-\sum_{i = 1}^{H'} \hat{l}_i q_i}{1-r}\right)\right] \right] 
 \forall i = 1, \ldots, H' \nonumber \\
 \Rightarrow \;  &q_i = p_i  \left( \sum_{i = 1}^{H'}\frac{q^2_i}{p_i}\right) \left\lbrace 1 + \frac{\left( \hat{l}_i - \sum_{i = 1}^{H'} \hat{l}_iq_i \right)}{\sqrt{\frac{\left(\sum_{i=1}^{H'} \frac{q_i^2}{p_i}\right) \mathcal{I}^{R}_{\text{kl}}(m, 2)}{\delta }}}
 \left[ \ln\left( \frac{\sum_{i = 1}^{H'} \hat{l}_i q_i}{r}\right) - \ln \left( \frac{1-\sum_{i = 1}^{H'} \hat{l}_i q_i}{1-r}\right)\right] \right \rbrace
 \forall i = 1, \ldots, H'
\end{align}
For feasibility, we need $\sum_{i =1}^{H} q_i = 1$. Therefore by using a suitable normalization constant $Z_{\text{kl}, \chi^2}$, we get the normalized equation in $q_i$s:
\begin{equation}
q_i = \frac{1}{Z_{\text{kl}, \chi^2}} p_i \left(\sum_{i=1}^{H'} \frac{q_i^2}{p_i} \right) \left \lbrace 1 + \frac{\left(  \hat{l}_i - \sum_{i=1}^{H}\hat{l}_iq_i \right)}{\sqrt{\frac{\left(\sum_{i=1}^{H'} \frac{q_i^2}{p_i}\right) \mathcal{I}^{R}_{\text{kl}}(m, 2)}{\delta}} } \left[ \ln\left( \frac{(1 - r) \sum_{i = 1}^{H'} \hat{l}_i q_i}{r(1-\sum_{i = 1}^{H'} \hat{l}_i q_i)}\right) \right] \right \rbrace \forall i = 1, \ldots, H'
\end{equation}
The above is the fixed point equation (FPE) which identifies a stationary point of the bound minimization problem \eqref{appdx_eqn:BklChi2OP}.
\end{proof}

\begin{cor}[Optimal posterior on an ordered subset support]
When prior is uniform distribution on $\mathcal{H}$, among all the posteriors with support as subset of size exactly $H'$, a stationary point $Q^{\text{FP}}_{\text{kl}, \chi^2}(H')$ for \eqref{appdx_eqn:BklChi2OP} can be obtained as the solution to the following fixed point equation:
\begin{equation}
q_i = \frac{1}{Z_{\text{kl}, \chi^2}} \left(\sum_{i=1}^{H'} q_i^2 \right) \left \lbrace 1 + \frac{\left(\hat{l}_i  -\sum_{i=1}^{H'}\hat{l}_iq_i \right)}{\sqrt{\frac{H \left(\sum_{i=1}^{H'} q_i^2\right) \mathcal{I}^{R}_{\text{kl}}(m, 2)}{\delta}} } \left[ \ln\left( \frac{(1 - r) \sum_{i = 1}^{H'} \hat{l}_i q_i}{r(1-\sum_{i = 1}^{H'} \hat{l}_i q_i)}\right) \right] \right \rbrace 
\label{appdx_eqn:OptQklChi2.subset.unifP}
\end{equation}
for $i = 1, \ldots, H$ where $Z_{\text{kl}, \chi^2}$ is a suitable normalization constant and $r$ is the solution to \eqref{appdx_eqn:klChi2DCcons} and \eqref{appdx_eqn:r>EQl} for a given $Q = (q_1, \ldots, q_H)$.
\end{cor}

\begin{lem}
When all the classifiers have same empirical risk (all $\hat{l}_i$s are same), the optimal posterior for the bound minimization problem \eqref{appdx_eqn:BklChi2OP} is $Q \equiv P$.
\end{lem}

KL-distance based bound minimization is non-convex with multiple stationary points which makes it difficult to identify the global minimum even by FP scheme. The iterative root finding algorithm adds to the computational complexity of the bound minimization algorithm. 

\subsection{Convex-concave procedure for finding a local solution for minimization of $B_{\text{kl}, \chi^2}(Q)$} \label{appdx_secn:klChi2CCP}
We have seen that our optimization problem \eqref{appdx_eqn:BklChi2OP} for finding the bound $B_{\text{kl}, \chi^2}(Q)$ consists of a linear objective function and linear constraints, except for the constraint \eqref{appdx_eqn:klChi2DCcons}, which takes the form:
\begin{eqnarray}
&kl\left(\mathbb{E}_Q [\hat{l}], r \right) = \sqrt{ \frac{(\chi^2(Q||P) + 1) \mathcal{I}^R_{\text{kl}}(m, 2)}{\delta} }\\
&\hspace{-1cm} \Leftrightarrow \left(\sum_{i = 1}^H \hat{l}_iq_i \right) \ln\left( \frac{\sum\limits_{i = 1}^H \hat{l}_iq_i}{r}\right) + \left(1 - \sum\limits_{i = 1}^H \hat{l}_iq_i \right) \ln\left( \frac{ 1  - \sum\limits_{i = 1}^H \hat{l}_iq_i}{1 - r}\right) =  \sqrt{ \frac{\left(\sum_{i =1}^H \frac{q^2_i}{p_i} \right) \mathcal{I}^R_{\text{kl}}(m, 2)}{\delta} }
\end{eqnarray}
We know that $\chi^2[Q||P]$ is jointly convex in both its arguments \cite{VanErvenHarremoes2014}. Based on the proof of Theorem \ref{appdx_thm:convexity.BlinChi2}, we have that $\sqrt{\left(\sum_{i =1}^H \frac{q^2_i}{p_i} \right) }$ is a convex function of $Q$. And hence, for given system parameters $m$ and $\delta$, right hand side of the above constraint is a convex function of $Q$. The left hand side is a composition of two functions: $\mathbb{E}_Q[\hat{l}]$ (a linear function) and $kl(p, q)$ (a jointly convex function). The superposition of a convex function and an affine mapping is convex, provided that it is finite at least at one point \cite{anatoli2015notesconvexity,boyd2004convex}. Hence, it is established that $kl\left(\mathbb{E}_Q [\hat{l}], r \right)$ is convex in its arguments $(Q, r)$. This implies that the constraint \eqref{appdx_eqn:klChi2DCcons} is a difference of convex (DC) function and the associated optimization problem \eqref{appdx_eqn:BklChi2OP} is a DC program.

Reformulating the original problem \eqref{appdx_eqn:BklChi2OP} in terms of all inequality constraints of the form $f(x) - g(x) \leq 0$, we have:
\begin{subequations} \label{appdx_eqn:BklChi2OP.preCCP}
\begin{align} 
&\min_{q_1, \ldots, q_H, r} r\\
&\hspace{-8mm} \left(\sum_{i = 1}^H \hat{l}_iq_i \right) \ln\left( \tfrac{\sum\limits_{i = 1}^H \hat{l}_iq_i}{r}\right) + \left(1 - \sum\limits_{i = 1}^H \hat{l}_iq_i \right) \ln\left( \tfrac{ 1  - \sum\limits_{i = 1}^H \hat{l}_iq_i}{1 - r}\right) - \sqrt{ \tfrac{\left(\sum_{i =1}^H \frac{q^2_i}{p_i} \right) \mathcal{I}^R_{\text{kl}}(m, 2)}{\delta} } \leq 0 \label{appdx_eqn:klChi2DCcons.preCCP}\\ 
&\hspace{-8mm} \sqrt{ \tfrac{\left(\sum_{i =1}^H \frac{q^2_i}{p_i} \right) \mathcal{I}^R_{\text{kl}}(m, 2)}{\delta} } - \left(\sum_{i = 1}^H \hat{l}_iq_i \right) \ln\left( \tfrac{\sum\limits_{i = 1}^H \hat{l}_iq_i}{r}\right) + \left(1 - \sum\limits_{i = 1}^H \hat{l}_iq_i \right) \ln\left( \tfrac{ 1  - \sum\limits_{i = 1}^H \hat{l}_iq_i}{1 - r}\right) \leq 0 \label{appdx_eqn:klChi2DCcons.reverse.preCCP}\\
& \sum_{i = 1}^H \hat{l}_iq_i - r \leq 0 \label{appdx_eqn:r>EQlhat.klChi2.preCCP}\\
& \sum_{i = 1}^H q_i = 1 \\
& -q_i \leq 0 \; \forall i = 1, \ldots, H 
\end{align} 
\end{subequations}

To apply the convex-concave procedure (CCP), we determine the approximations to the DC functions \eqref{appdx_eqn:klChi2DCcons.preCCP} and  \eqref{appdx_eqn:klChi2DCcons.reverse.preCCP}, at a point $(Q^0, r^0)$ which is feasible to \eqref{appdx_eqn:BklChi2OP.preCCP}, and equivalently to \eqref{appdx_eqn:BklChi2OP}. Let $\widehat{kC}_1((Q,r); (Q^0, r^0))$ denote the linear under-approximation to the function $kC_1(Q, r) := \sqrt{ \tfrac{\left(\sum_{i =1}^H \frac{q^2_i}{p_i} \right) \mathcal{I}^R_{\text{kl}}(m, 2)}{\delta} }$ in \eqref{appdx_eqn:klChi2DCcons.preCCP} at $(Q^0, r^0)$.

\begin{align*}
\widehat{kC}_1((Q,r); (Q^0, r^0)) &:= kC_1(Q^0, r^0) + \langle \nabla kC_1(Q^0, r^0), \left((Q - Q^0), (r  - r^0)\right) \rangle\\
&= \sqrt{ \tfrac{\left(\sum_{i =1}^H \frac{(q^0_i)^2}{p_i} \right) \mathcal{I}^R_{\text{kl}}(m, 2)}{\delta} } +  \left( \sum_{i = 1}^{H} \frac{\partial  kC_1 }{\partial q_i} \bigg\vert_{q_i = q^0_i} \cdot (q_i - q^0_i)  \right) + 0 \cdot (r  - r^0) \\
&= \sqrt{ \tfrac{\left(\sum_{i =1}^H \frac{(q^0_i)^2}{p_i} \right) \mathcal{I}^R_{\text{kl}}(m, 2)}{\delta} } + \sqrt{ \tfrac{\mathcal{I}^R_{\text{kl}}(m, 2)}{\delta \left(\sum_{i =1}^H \frac{(q^0_i)^2}{p_i} \right)} } \left(\sum_{i = 1}^{H} \frac{q^0_i}{p_i} (q_i - q^0_i) \right) \displaybreak[3] \\
&= \cancel{\sqrt{ \tfrac{\left(\sum_{i =1}^H \frac{(q^0_i)^2}{p_i} \right) \mathcal{I}^R_{\text{kl}}(m, 2)}{\delta} }} + \sqrt{ \tfrac{\mathcal{I}^R_{\text{kl}}(m, 2)}{\delta \left(\sum_{i =1}^H \frac{(q^0_i)^2}{p_i} \right)} } \left(\sum_{i = 1}^{H} \frac{q^0_i q_i}{p_i} \right) - \cancel{\sqrt{ \tfrac{\left(\sum_{i =1}^H \frac{(q^0_i)^2}{p_i} \right) \mathcal{I}^R_{\text{kl}}(m, 2)}{\delta} }} \\
&= \sqrt{ \tfrac{\mathcal{I}^R_{\text{kl}}(m, 2)}{\delta \left(\sum_{i =1}^H \frac{(q^0_i)^2}{p_i} \right)} } \left(\sum_{i = 1}^{H} \frac{q^0_i q_i}{p_i} \right)
\end{align*}

Recall the linear under-approximation $\widehat{kK}_2((Q,r); (Q^0, r^0))$ to the function $kK_2(Q, r) :=kl\left(\sum_{i = 1}^H \hat{l}_iq_i, r \right)$ at $(Q^0, r^0)$.
\begin{align*}
\widehat{kK}_2((Q,r); (Q^0, r^0)) &:= kK_2(Q^0, r^0) + \langle \nabla kK_2(Q^0, r^0), \left(Q - Q^0, r  - r^0\right) \rangle\\
&= \ln{\left(\frac{1-\sum_{i=1}^{H} \hat{l}_i q_i^0}{1-r^0}\right)} + \sum_{i=1}^{H} \hat{l}_iq_i \left[\ln{\left(\frac{\sum_{i=1}^{H} \hat{l}_i q^0_i}{r^0} \right)} - \ln{\left(\frac{1-\sum_{i=1}^{H} \hat{l}_i q^0_i}{1-r^0} \right)} \right] \\
	&\hspace{3cm}+ \left[\frac{r^0 - \sum_{i=1}^{H} \hat{l}_i q^0_i }{r^0 (1 - r^0)}\right] (r - r^0) 
\end{align*}

Using the above linear approximations $\widehat{kC}_1((Q,r); (Q^0, r^0))$ and $\widehat{kK}_2((Q,r); (Q^0, r^0))$ in \eqref{appdx_eqn:klChi2DCcons.preCCP} and \eqref{appdx_eqn:klChi2DCcons.reverse.preCCP}, we can invoke the CCP procedure described in \cite{CCP2016LippBoyd} to get a local minimizer to the KL-distance based bound minimization problem \eqref{appdx_eqn:BklChi2OP}.

\section{Computational Illustrations for SVMs} \label{appdx_secn:SVM_PACB}
The datasets that we have considered for our computations, the scheme used to generate classifiers and compute risk values are same as the ones considered in \cite{sahu2019PACBKLarxiv}. On this set of base classifiers, we compare the optimal PAC-Bayesians posteriors for the case of $\chi^2$-divergence, obtained using the FP scheme and the solver for the different $\phi$ functions considered.



\begin{table}[]
{ \setlength{\tabcolsep}{0.5mm}
\footnotesize
\begin{tabular}{|c|c|c|c|c|c|c|c|c|c|c|c|}
\hline
\diagbox{\textbf{Dataset}}{\textbf{H}} & \multicolumn{2}{c|}{\textbf{50}} & \multicolumn{2}{c|}{\textbf{200}} & \multicolumn{2}{c|}{\textbf{500}} & \multicolumn{2}{c|}{\textbf{1000}} & \multicolumn{2}{c|}{\textbf{1990}} \\  \cline{2-11}
\makecell{(Validation set \\size, $v$)}& $B^{FP}_{\text{kl}, \chi^2}$ &  $B^{solver}_{\text{kl}, \chi^2}$  &  $B^{FP}_{\text{kl}, \chi^2}$ &  $B^{solver}_{\text{kl}, \chi^2}$ & $B^{FP}_{\text{kl}, \chi^2}$ &  $B^{solver}_{\text{kl}, \chi^2}$ & $B^{FP}_{\text{kl}, \chi^2}$ &  $B^{solver}_{\text{kl}, \chi^2}$  & $B^{FP}_{\text{kl}, \chi^2}$ &  $B^{solver}_{\text{kl}, \chi^2}$ \\ \hline
\makecell{Spambase \\ $(v = 1840)$ }&0.43076 &0.84883(I) &0.48821 & 0.96762(I)&0.52451 & 0.58991(E)& 0.54945&0.61292(E) &0.57082 &0.63182(E)
 \\ \hline
 \makecell{Bupa \\ $( v = 138)$ }&0.59598 & 0.98459(I)&0.67963 &0.99969(I) & 0.74690& 0.82615(E)&0.79420 &0.86137(E)&0.83864&0.90011(E)
 \\ \hline
 \makecell{Mammographic \\ $( v = 332)$ } &0.55152 & 0.57227(I)&0.56886 &0.99729(I) & 0.58527& 0.99988(I)& 0.60076&0.66770(E) &0.62780 &0.70106(E)
 \\ \hline
\makecell{Wdbc \\ $( v = 227)$ }&0.44717 &0.49717 & 0.45708& 0.45708& 0.46697& 0.9999(I) &0.47589 &0.9999(R) &0.49508 &0.59303(E)
 \\ \hline
\makecell{Banknote \\ $( v = 549)$ }&0.22776 & 0.22776&0.23557 &0.23557 &0.24665 &0.30710(R) & 0.25209& 0.41775 (v) &0.26038 &0.45979(M)
 \\ \hline
\makecell{Mushroom \\ $( v = 2257)$ } &0.17498 &0.17498 &0.17498 & 0.17498&0.17657 & 0.17657 &0.18097 &0.18639(R) & 0.18660&0.25445(M)
 \\ \hline
\makecell{Ionosphere \\ $( v = 140)$ }&0.47554 &0.47544 &0.49484 & 0.49484& 0.51114&0.97958(M)  &0.57477 &0.70459(E) &0.67610 &0.80135(E)
 \\ \hline
\makecell{Waveform \\ $( v = 1323)$ }& 0.27118& 0.27118& 0.05185&0.05185 &0.28334 &0.98959(I) &0.28612 &0.99833(I) &0.28978 &0.28978
 \\ \hline
\makecell{Haberman \\ $( v = 122)$ }& 0.68001& 0.93379(I)& 0.69962& 0.99270(I)&0.72220 &0.81692(E) & 0.73266& 0.80557(E)&0.73790 &0.80836(E)
 \\ \hline
\end{tabular}}
\caption[Bound values for kl-$\chi^2$ case]{\textbf{Bound values for kl-$\chi^2$ case}: Comparing the bound values $B^{FP}_{\text{kl},\chi^2}$ and $B^{solver}_{\text{kl},\chi^2}$ for the posterior obtained via fixed point equation \eqref{appdx_eqn:OptQklChi2.subset.unifP} with the linear search algorithm \ref{appdx_algo:optQsubset.unifP} and the optimal posterior for minimizing the PAC-Bayesian bound \eqref{appdx_eqn:BklChi2OP} for the KL-distance with chi-squared divergence between the prior and posterior distributions. The fixed point equation always converges and identifies the local minimum output by the \texttt{Ipopt} solver, even when the solver fails to identify a solution for reasons like local infeasibility (I), Restoration Phase Failed (R), Maximum number of iterations exceeded (M), Unknown error (E), etc.}
\label{appdx_tab:Bnd.klChi2}
\end{table}

\begin{table}[]
{ \setlength{\tabcolsep}{0.5mm}
\footnotesize
\begin{tabular}{|c|c|c|c|c|c|c|c|c|c|c|c|}
\hline
\diagbox{\textbf{Dataset}}{\textbf{H}} & \multicolumn{2}{c|}{\textbf{50}} & \multicolumn{2}{c|}{\textbf{200}} & \multicolumn{2}{c|}{\textbf{500}} & \multicolumn{2}{c|}{\textbf{1000}} & \multicolumn{2}{c|}{\textbf{1990}} \\  \cline{2-11}
\makecell{(Test set \\size, $t$)}& $T^{FP}_{\text{kl}, \chi^2}$ &  $T^{solver}_{\text{kl}, \chi^2}$  &  $T^{FP}_{\text{kl}, \chi^2}$ &  $T^{solver}_{\text{kl}, \chi^2}$ & $T^{FP}_{\text{kl}, \chi^2}$ &  $T^{solver}_{\text{kl}, \chi^2}$ & $T^{FP}_{\text{kl}, \chi^2}$ &  $T^{solver}_{\text{kl}, \chi^2}$  & $T^{FP}_{\text{kl}, \chi^2}$ &  $T^{solver}_{\text{kl}, \chi^2}$ \\ \hline
\makecell{Spambase \\ $( t =921)$ } &0.17387 &0.16069(I) &0.20623 & 0.16069(I) &0.23489 & 0.24251 (E) &0.25548 &0.26379(E) &0.27345 &0.27927(E)
 \\ \hline
 \makecell{Bupa \\ $( t= 69)$ } & 0.16630& 0.08695(I) & 0.22385& 0.08695(I) &0.28570 &0.30806(E) & 0.34015&0.37587(E) &0.42343 &0.47202(E)
 \\ \hline
 \makecell{Mammographic \\ $( t =166)$ } & 0.21192&0.21030(I) &0.21986 &0.20481(I) &0.22636 &0.20481(I) &0.22813 &0.22768(E) &0.23266 &0.23402(E)
 \\ \hline
\makecell{Wdbc \\ $( t =115)$ }&0.06490 &0.06956(I) &0.06025 &0.06025&0.06579 &0.06956(I) & 0.07061& 0.06942(R)& 0.07808&0.08818(E)
 \\ \hline
\makecell{Banknote \\ $( t =275)$ } & 0.00036& 0.00036&0.00203 &0.00203 &0.00424 & 0.00235(R)& 0.00569 & 0.002401(M)  & 0.00629 & 0.00359(M)
 \\ \hline
\makecell{Mushroom \\ $( t = 1129)$ } &0 &0 &$4.42  e{-06}$ & $4.42  e{-06}$ &0.00021 &0.00021 &0.00055 &0.00044(R) &0.00092 &0.00047(M)
 \\ \hline
\makecell{Ionosphere \\ $( t =71)$ } & 0.04561& 0.04561& 0.04380&0.04380 & 0.04864&0.03991(M)  & 0.07364&0.14996(E) &0.07973 &0.25377(E)
 \\ \hline
\makecell{Waveform \\ $( t =662)$ }&0.05210 &0.05210 & 0.27939& 0.27939 &0.05198 &0.05898(I) & 0.05154& 0.05891(I) &0.05120 &0.05120
 \\ \hline
\makecell{Haberman \\ $( t =62)$ }&0.29362 &0.28629(I)& 0.28702& 0.28586(I) &0.28849 & 0.28681(E)& 0.28936&0.28943(E) &0.28983 &0.28997(E)
 \\ \hline
\end{tabular}}
\caption[Test error rates for kl-$\chi^2$ case]{\textbf{Test error rates for $kl-\chi^2$ case}: Comparing the test error rates $T^{FP}_{\text{kl},\chi^2}$ and $T^{solver}_{\text{kl},\chi^2}$  when using the posterior obtained via fixed point equation \eqref{appdx_eqn:OptQklChi2.subset.unifP} with the linear search algorithm \ref{appdx_algo:optQsubset.unifP} and the optimal posterior for minimizing the PAC-Bayesian bound \eqref{appdx_eqn:BklChi2OP} for the KL-distance with chi-squared divergence between the prior and posterior distributions. The fixed point equation always converges and identifies the local minimum output by the \texttt{Ipopt} solver, even when the solver fails to identify a solution for reasons like local infeasibility (I), Restoration Phase Failed (R), Maximum number of iterations exceeded (M), Unknown error (E), etc.}
\label{appdx_tab:TestErr.klChi2}
\end{table}

\begin{table}[ht]
\centering 
{ \footnotesize
\setlength{\tabcolsep}{0.1em}
\begin{tabular}{ |c | c c c | c c c|}
\hline
\textbf{Dataset} & \multicolumn{3}{c|}{\textbf{PAC-Bayesian Bound}} & \multicolumn{3}{c|}{\textbf{Average Test Error}} \\
\cline{2-7}
  & $B^{FP}_{\text{kl}, \chi^2}$ & Range($B^{CCP}_{\text{kl}, \chi^2}$) & Mean($B^{CCP}_{\text{kl}, \chi^2}$) & $T^{FP}_{\text{kl}, \chi^2}$ & Range($T^{CCP}_{\text{kl}, \chi^2}$) & Mean($T^{CCP}_{\text{kl}, \chi^2}$)
\\ 
\hline
Spambase  & 0.43076 & [0.46072, 0.53931] & 0.48658 $\pm$ 0.01139 & 0.17387 & [0.16480, 0.17548] & 0.17107 $\pm$ 0.00175\\
Bupa &  0.59598 & [0.63927, 0.77591] & 0.67911 $\pm$ 0.01681 & 0.16630 & [0.14297, 0.18139] & 0.16408 $\pm$ 0.00565\\
\makecell{Mammographic} & 0.55152 & [0.58868, 0.70839] & 0.62348 $\pm$ 0.01459 & 0.21192 & [0.20581, 0.21725] & 0.21178 $\pm$ 0.00186 \\
Wdbc &0.44717 & [0.48998, 0.60177] & 0.53482 $\pm$ 0.01825 & 0.06490 & [0.05921, 0.07138] & 0.06491 $\pm$ 0.00174\\
Banknote & 0.22776 & NA & NA & 0.00036 & NA & NA \\
Mushroom & 0.17498& NA & NA & 0 & NA & NA\\
Ionosphere & 0.47534 & [0.51826, 0.67369] & 0.57483 $\pm$ 0.02047 & 0.04561 & [0.03955, 0.05269] & 0.04557 $\pm$ 0.00213\\
Waveform & 0.27118 & [0.30395, 0.38119] & 0.33086 $\pm$ 0.01219 & 0.05210 & [0.04983, 0.05436] & 0.05215 $\pm$ 0.00069\\
Haberman & 0.68001 & [0.73588, 0.82356] & 0.77063 $\pm$ 0.01396 & 0.29362 & [0.28378, 0.30268] & 0.29376 $\pm$ 0.00290\\ 
\hline
\end{tabular}
}
\caption[Comparing bound values and test error rates of FP based posterior with CCP based posterior for kl-$\chi^2$ case]{We compare the bound values and test error rates of the optimal posterior obtained via Fixed Point (FP) scheme and the posterior based on Convex-Concave Procedure (CCP) for minimizing the PAC-Bayesian bound $B_{\text{kl}, \chi^2}$ based on KL-distance function with $\chi^2$-divergence measure. The CCP based posteriors are identified by the bound minimization model described in Section \ref{appdx_secn:klChi2CCP}. The bound values and test error rates for FP scheme based solution are denoted by $B^{FP}_{\text{kl}, \chi^2}$ and $T^{FP}_{\text{kl}, \chi^2}$. Similarly, the bound values and test error rates of the CCP based posterior are denoted by $B^{CCP}_{\text{kl}, \chi^2}$ and $T^{CCP}_{\text{kl}, \chi^2}$. For computations, we consider SVM classifiers generated on nine datasets from UCI repository \cite{UCI:2017} using the scheme in Section 7 of the main paper (also considered in \cite{sahu2019PACBKLarxiv}) for $H=50$ values in $\Lambda = \lbrace 0.1, 0.11, \ldots, \rbrace$. We run the CCP procedure for 1000 different initializations of posterior $Q^0$ (as done in \cite{CCP2016LippBoyd}). The range, mean and standard deviation of the bound values and average test error rates of the CCP based posteriors obtained by these 1000 initializations are tabulated above. We notice that $B^{FP}_{\text{kl}, \chi^2}$ is always better than $B^{CCP}_{\text{kl}, \chi^2}$ and $T^{FP}_{\text{kl}, \chi^2}$ is comparable with mean value of $T^{CCP}_{\text{kl}, \chi^2}$ for different datasets considered. This might be because FP scheme identifies the global minimum for kl-$\chi^2$ based bound minimization problem, whereas CCP converges to a local solution or a stationary point. `NA' denotes the cases where the CCP cannot provide linear approximation to $kl(\mathbb{E}_Q[\hat{l}], r)$ because a subgradient cannot be determined when $\mathbb{E}_Q[\hat{l}]$ takes the boundary value zero. Such cases usually occur for almost separable datasets -- Banknote and Mushroom, where the quantity $\mathbb{E}_Q[\hat{l}] = 0$ for any distribution $Q$ since all $\hat{l}_i$s take value zero for $i = 1, \ldots, 50$.}
\label{appdx_tab:Bnd.TestErr.FPnCCP.klChi2}
\end{table}

\subsection{Illustration of various optimal posteriors, $Q^{\ast}_{\phi, \chi^2}$ }
\begin{figure}
\centering
\includegraphics[width=0.48\textwidth]{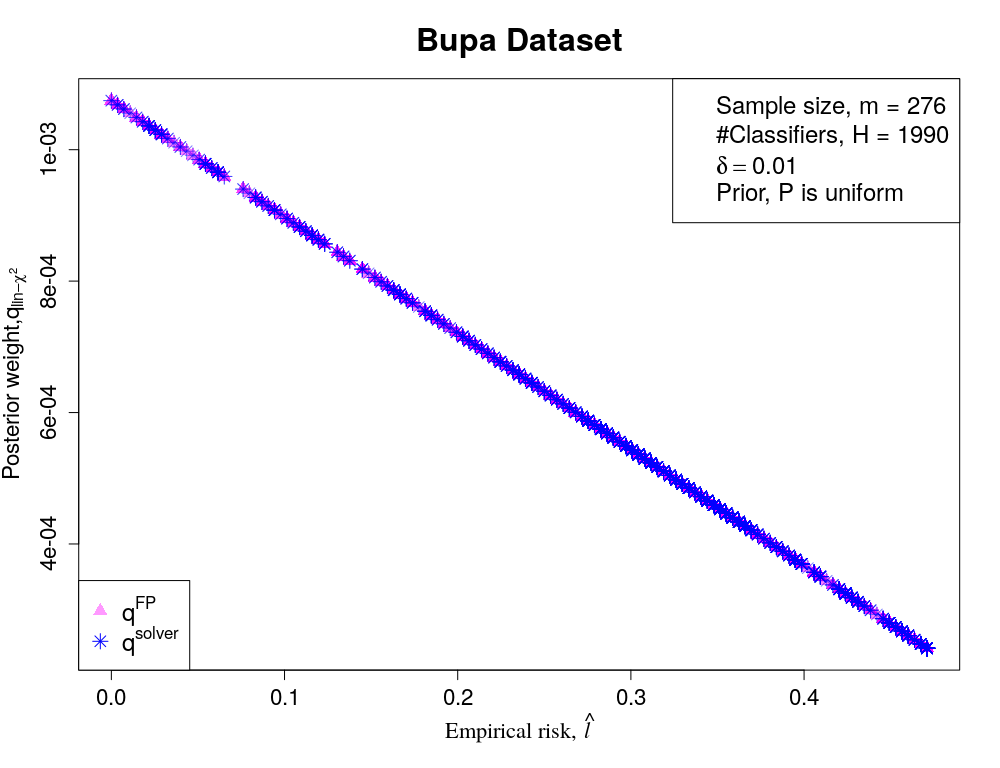}
\includegraphics[width=0.48\textwidth]{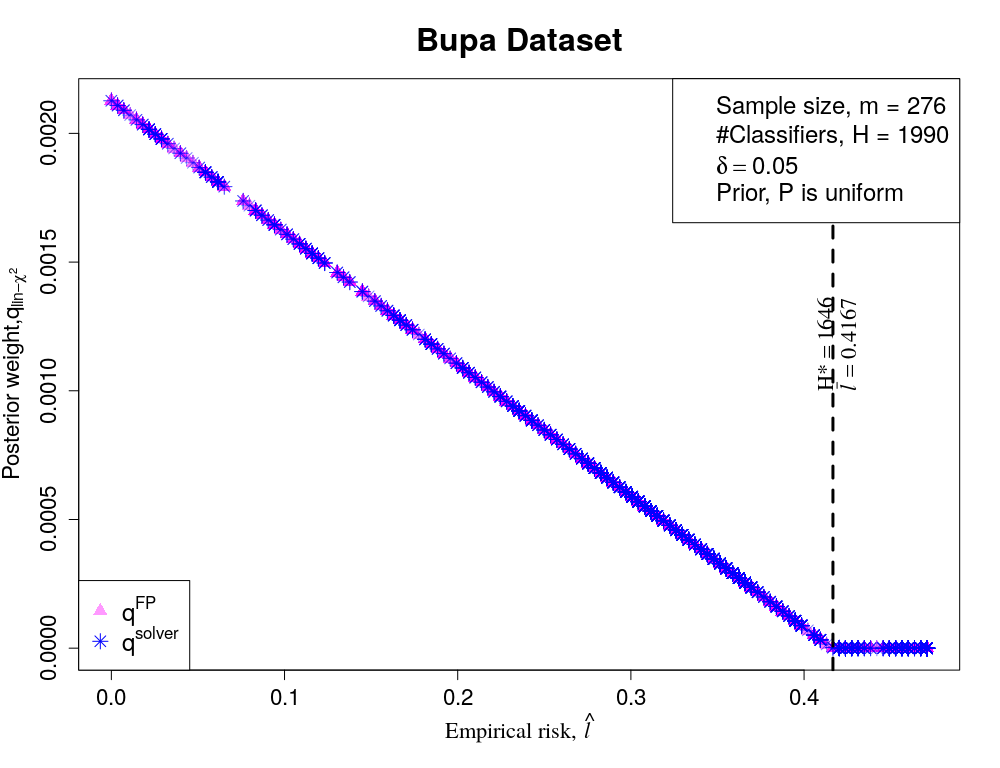}
\includegraphics[width=0.48\textwidth]{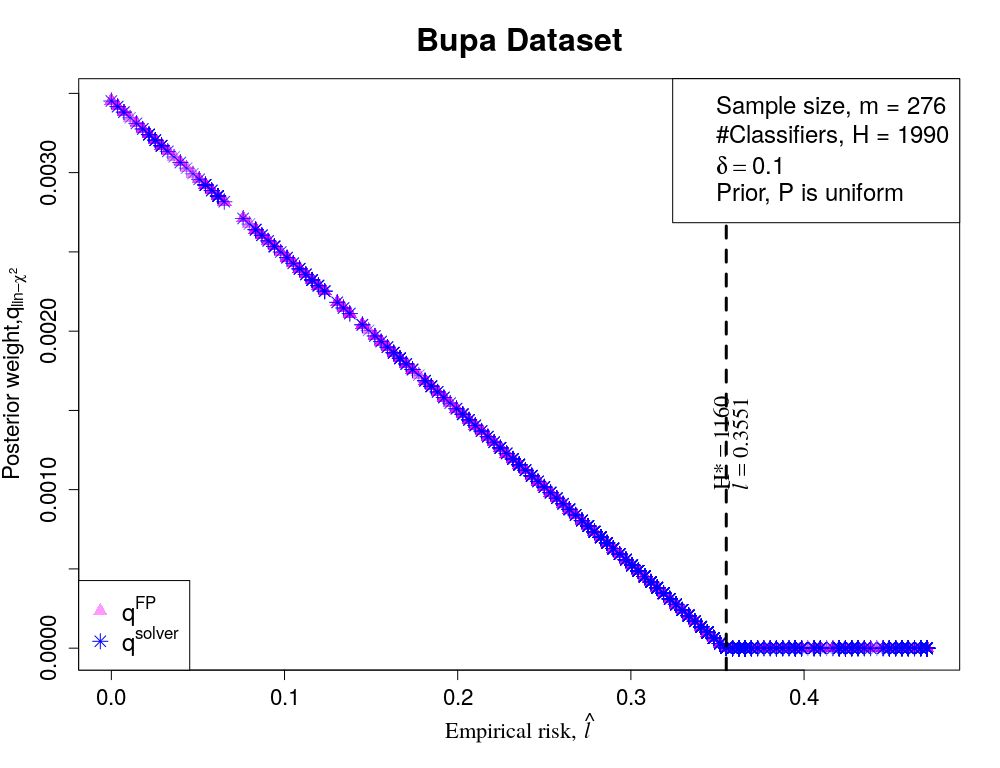}
\includegraphics[width=0.48\textwidth]{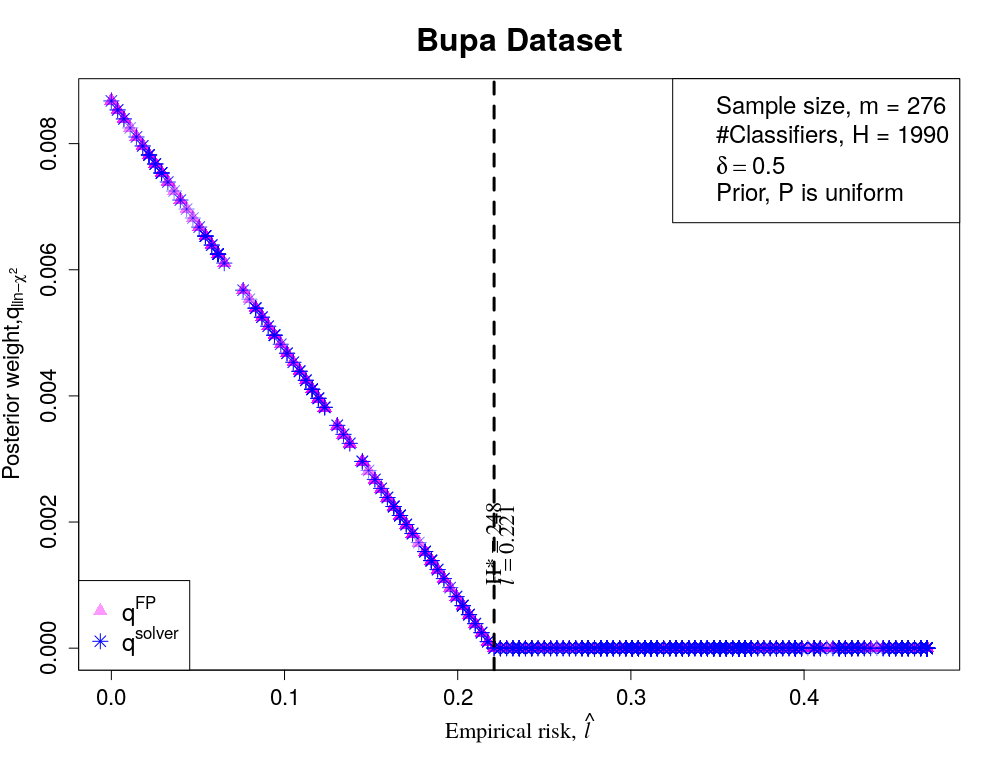}
\caption{\textbf{Illustration of variation of subset support for $Q^{\ast}_{\text{lin}, \chi^2}$ as the PAC-Bayesian confidence level $\delta$ changes}. We consider Bupa dataset (345 samples, 6 features) with a training sample of size $m = 276$ for determining our SVM classifiers using $H=1990$ regularization parameter values from the set $\Lambda = \lbrace 0.01, \ldots, 20 \rbrace$. For a uniform prior distribution, the optimal posterior for linear distance function, $Q^{\ast}_{\text{lin}, \chi^2}$ is computed via \texttt{Ipopt} solver on the full simplex ($Q^{solver}$) as well fixed point equation (FPE) \eqref{appdx_eqn:OptQlinChi2.subset.unifP} on increasing ordered subsets of $\mathcal{H}$, denoted by $Q^{FP}$. We observe that the FPE correctly identifies the global minimum. In case of uniform prior, the posterior weights $q^{\ast}_{i, \text{lin}, \chi^2}$ are negatively proportional to the empirical risk $\hat{l}_i$ values of the classifiers in the ordered support set. $H^{\ast}$ denotes the support size of the optimal posterior and $\bar{l} := \hat{l}_{H^{\ast}}$ denotes the value of the empirical risk beyond which the posterior weights are zero. For given fixed parameters, namely $H$ and $m$, we consider three values of confidence level $\delta = 0.01, 0.05, 0.1$. The optimal subset size $H^{\ast}$ decreases as $\delta$ increases, allowing sparser posteriors. \label{appdx_fig:OptQlinChi2.subset.unifP}}
\end{figure}

\begin{figure}
\centering
\includegraphics[width=0.48\textwidth]{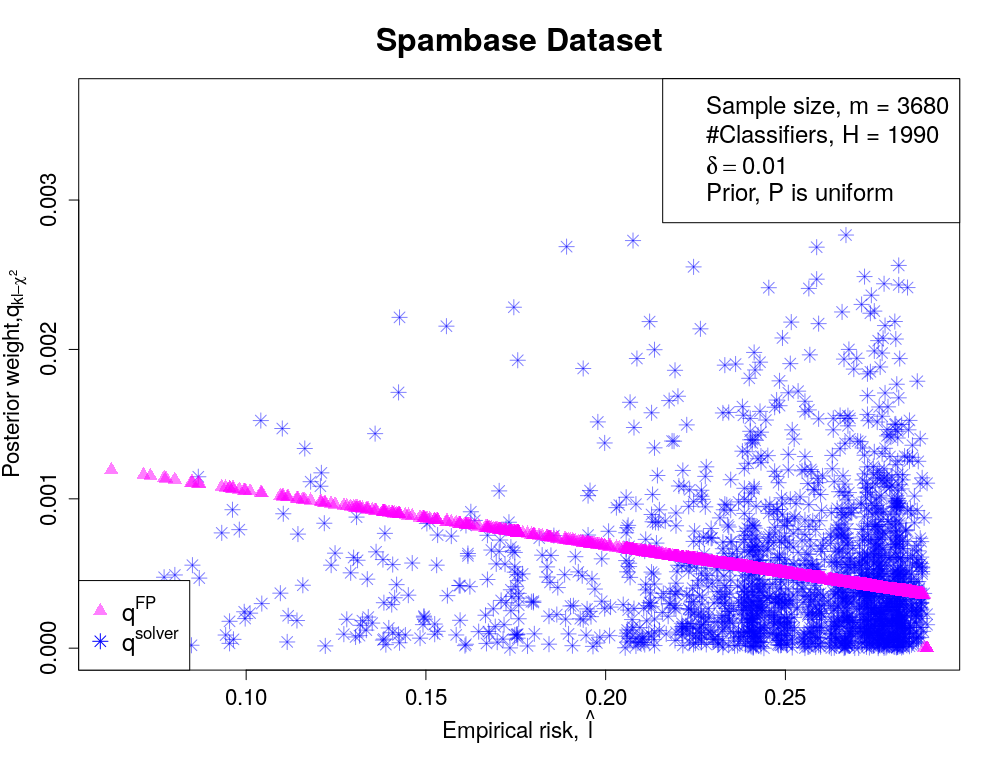}
\includegraphics[width=0.48\textwidth]{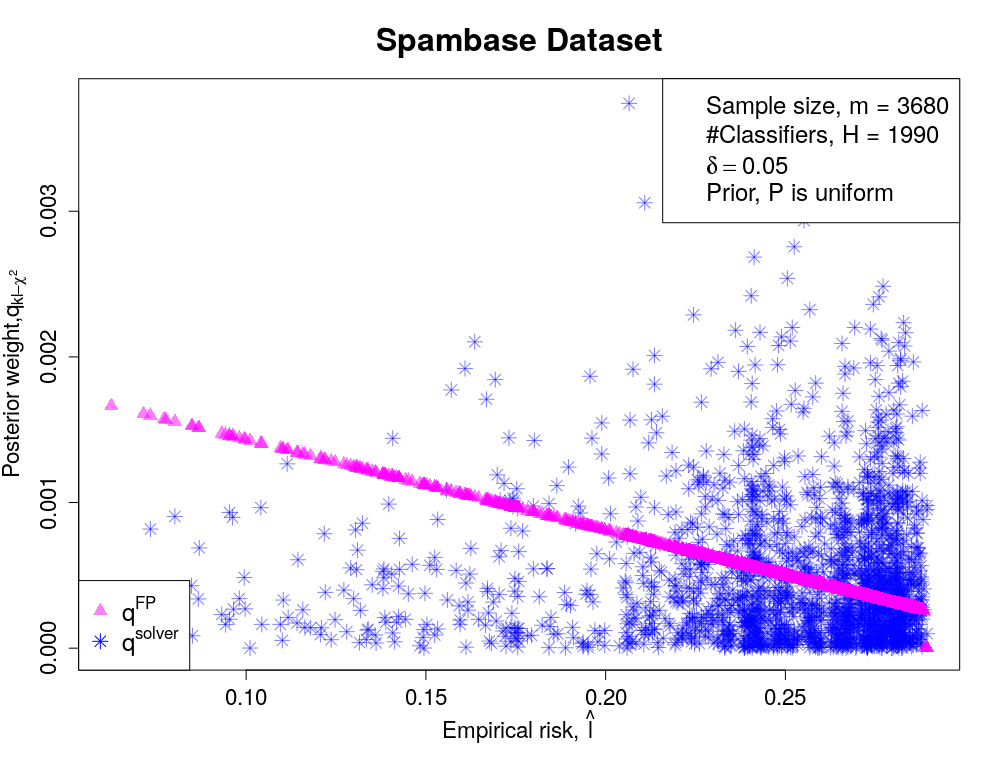}
\includegraphics[width=0.48\textwidth]{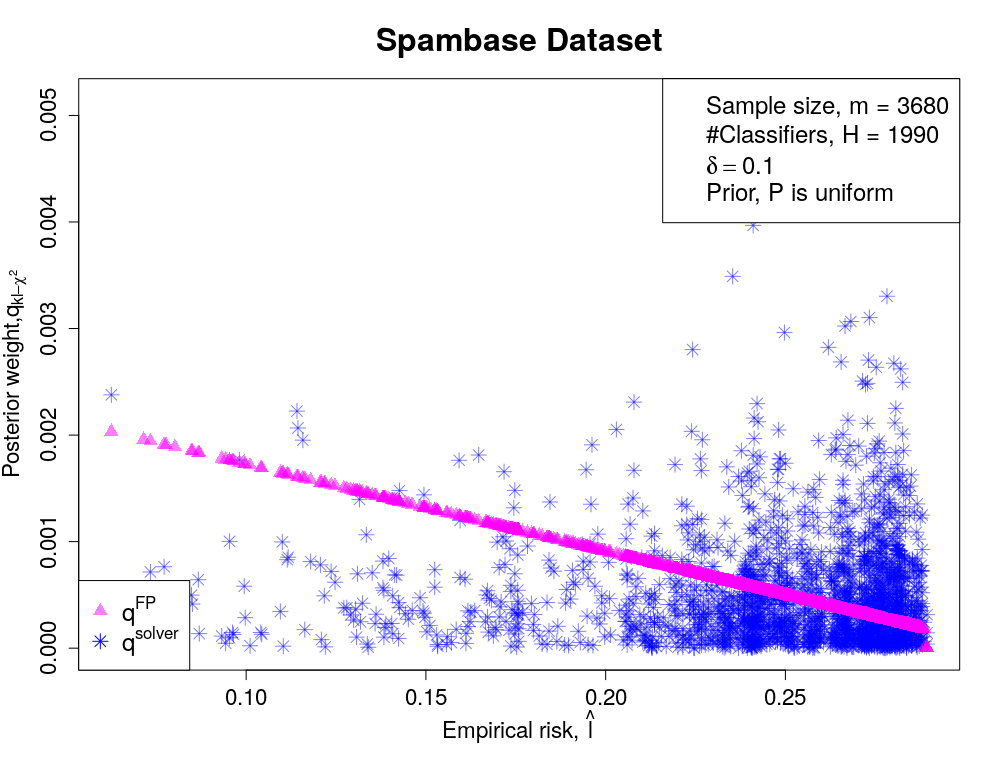}
\caption{\textbf{Illustration of (almost) full support for $Q^{FP}_{\text{kl}, \chi^2}$}. We consider Spambase dataset (4601 samples, 57 features) with a training sample of size $m = 3680$ for determining our SVM classifiers using regularization parameter values from the set $\Lambda = \lbrace 0.01, \ldots, 20 \rbrace$. For a uniform prior distribution, the optimal posterior for linear distance function, $Q^{\ast}_{\text{kl}, \chi^2}$ is computed via \texttt{Ipopt} solver on the full simplex ($Q^{solver}$) as well fixed point equation (FPE) \eqref{eqn:qFP_klChi2.unifP} on increasing ordered subsets of $\mathcal{H}$, denoted by $Q^{FP}$. We observe that the FPE always converges to a stationary point even when the solver throws up error and outputs infeasible solutions. In case of uniform prior, the posterior weights $q^{FP}_{i, \text{kl}, \chi^2}$ are negatively proportional to the empirical risk $\hat{l}_i$ values of the classifiers in the ordered support set. We notice that, the support size for $Q^{FP}$, denoted by $H^{FP} = 1986$ for all the three values of confidence level $\delta = 0.01, 0.05, 0.1$; implying that $Q^{FP}$ for kl-$\chi^2$ case has almost full support. \label{appdx_fig:QFPklChi2.fullsupport.unifP}}
\end{figure}

We present some graphs to illustrate the nature of the optimal posteriors that we have computed under the framework mentioned above. Figure \ref{appdx_fig:OptQlinChi2.subset.unifP} depicts the role of the confidence level $\delta$ in determining the optimal support size $H^{\ast}$ for $Q^{\ast}_{\text{lin}, \chi^2}$ in case of uniform prior. Figure \ref{appdx_fig:QFPklChi2.fullsupport.unifP} shows that the stationary point $Q^{FP}_{\text{kl}, \chi^2}$ obtained by the fixed point scheme has almost full support, and that the fixed point equation \eqref{eqn:qFP_klChi2.unifP} always converges to a solution even when the solver throws up an error due to issues such as \textbf{M}' (Maximum Number of iterations exceeded) or `\textbf{I}' (Locally infeasible solution) or `\textbf{E}' (Unknown Error)  or `\textbf{R}' (Restoration Phase Failed). Please see Table \ref{appdx_tab:Bnd.klChi2} for such examples.

\subsection{Sparsity and Concentration of Optimal Posteriors, $Q^{\ast}_{\phi, \chi^2}$ \label{appdx_secn:sparsityconc}}
We determine the optimal posteriors $Q^{\ast}_{\phi, \chi^2}$ for different distance functions, $\phi$s and compare their bound values and test error rates. To understand the differences between the nature of these posteriors for different choices of $\phi$, we need to compare these vectors of posterior weights. For a large $H$, as we have considered, it is difficult to compare these high-dimensional probability weight vectors elementwise. We use different measures for capturing the information from these posteriors. 

To measure the sparsity of the posteriors, we first compute their cumulative distribution functions (CDFs) denoted by $F_{Q^{\ast}_{\phi, \chi^2}} (\cdot)$. We consider three different significance levels $\alpha \in \{0.8, 0.9, 0.95\}$ and identify the number of classifiers $N_{\phi, \chi^2} (\alpha)$ out of $H = 1990$, required by $F_{Q^{\ast}_{\phi, \chi^2}}(\cdot)$ to achieve the given significance level $\alpha$. That is
\begin{equation}
N_{\phi, \chi^2} (\alpha) := \min \lbrace i \vert F_{Q^{\ast}_{\phi, \chi^2}}(i) \geq \alpha  \rbrace, \quad \alpha \in \{0.8, 0.9, 0.95\}.
\end{equation}

For a given level $\alpha$, a low $N_{\phi, \chi^2} (\alpha)$ indicates that the distribution $Q^{\ast}_{\phi, \chi^2}$ is sparse. In our computations, we observe that for the three significance levels $\alpha \in \{0.8, 0.9, 0.95\}$, the optimal posterior $Q^{\ast}_{\text{kl}, \chi^2}$ has large $N_{\text{kl}, \chi^2} (\alpha)$ values, implying almost full support. Whereas $Q^{\ast}_{\text{sq}, \chi^2}$ is sparse as reflected by low $N_{\text{sq}, \chi^2} (\alpha)$ values. (Please see Table \ref{appdx_tab:allQChi2.HHI} for the computed values.)

We quantify the level of concentration of a posterior distribution on its support via Herfindahl-Hirschman Index (HHI) \cite{Hirschman1945HHI,HHIwiki}. HHI is a prominent index in the economics literature, widely used for measuring the contribution of a sector in the economy or market share of a firm in the industry. It is an indicator of the amount of competition among the firms. It is defined as the square root of the sum of the squares of the contribution/market shares of the firms. It turns out that for any probability vector like our posteriors, HHI is equivalent to the $\ell_2$ norm of the given probability vector. A high HHI score indicates high concentration of probabilities and vice versa. HHI scores for $Q^{\ast}_{\phi, \chi^2}$ are given in Table \ref{appdx_tab:allQChi2.HHI}. We observe that $Q^{\ast}_{\text{sq}, \chi^2}$ has relatively high HHI score, indicating higher concentration compared to $Q^{\ast}_{\text{lin}, \chi^2}$ and $Q^{\ast}_{\text{kl}, \chi^2}$. The differences in concentration levels are remarkable in case of datasets with highly varying empirical risk values.

\begin{table}[ht]
\caption{\small Herfindahl - Hirschman Index (HHI) for measuring concentration of the PAC-Bayesian optimal posteriors $Q^{\ast}_{\phi, \chi^2}$ and number of classifiers, $N_{\phi, \chi^2}(\alpha)$, required by the cumulative distribution functions (CDFs), $F_{Q^{\ast}_{\phi, \chi^2}}(\cdot)$  to achieve a given significance level: $\alpha \in \lbrace 0.8, 0.9, 0.95 \rbrace$. The posteriors were computed on a set of SVM classifiers generated using $H = 1000$ regularization parameter values at confidence level $\delta = 0.01$. See Figure \ref{appdx_fig:allQChi2.cdf} for a visual illustration of the stochastic dominance and concentration levels of the posteriors $Q^{\ast}_{\phi, \chi^2}$ via CDFs, $F_{Q^{\ast}_{\phi, \chi^2}}(\cdot)$. Among the three posteriors, $Q^{\ast}_{\text{sq}, \chi^2}$ has the lowest values of $N_{\phi, \chi^2}(\alpha)$ for different values of $\alpha$ and the fastest convergence to have CDF as 1 across different datasets. We observe that $Q^{\ast}_{\text{sq}, \chi^2}$ has relatively high HHI score, indicating higher concentration compared to $Q^{\ast}_{\text{lin}, \chi^2}$ and $Q^{\ast}_{\text{kl}, \chi^2}$ in case of datasets with high variation in empirical risk values, such as Spambase and Bupa. In Spambase dataset HHI score for $Q^{\ast}_{\text{sq}, \chi^2}$ is almost 5 times the scores for $Q^{\ast}_{\text{lin}, \chi^2}$ and $Q^{\ast}_{\text{kl}, \chi^2}$. Whereas in Bupa dataset HHI score for $Q^{\ast}_{\text{sq}, \chi^2}$ is almost 3 times the scores for $Q^{\ast}_{\text{lin}, \chi^2}$ and $Q^{\ast}_{\text{kl}, \chi^2}$. In case of almost separable datasets such as Wdbc or Banknote, HHI is same for all three posteriors and even $N_{\phi, \chi^2}(\alpha)$ for different $\alpha$ levels are almost the same, indicating that the posteriors are comparable. Figure \ref{appdx_fig:allQChi2.cdf} provides a visual confirmation to this claim where the posteriors are observed to be overlapping for the two datasets -- Wdbc and Banknote. \label{appdx_tab:allQChi2.HHI}}
\centering
\begin{tabular}{c c c c c}
\toprule
Dataset &\makecell{Significance \\ level, $\alpha$} & \multicolumn{3}{c}{\makecell{Number of classifiers, $N_{\phi, \chi^2}(\alpha)$, required by \\ the cumulative distribution function (CDF), $F_{Q^{\ast}_{\phi, \chi^2}}(\cdot)$\\ to achieve a given significance level, $\alpha$}} \\
& &  $N_{\text{lin}, \chi^2}(\alpha)$ &  $N_{\text{sq}, \chi^2}(\alpha)$ & $N_{\text{kl}, \chi^2}(\alpha)$ \\
\midrule
\multirow{4}{*}{\makecell{Mammo-\\-graphic \\ $(v = 332)$ }}& 0.8 & 782 & 729 & 773 \\ 
& 0.9 & 890 & 856 & 884\\ 
& 0.95 & 944 & 925 & 941\\ 
\cmidrule{2-5}
&  \textbf{HHI} & 0.0317 & 0.0325 & 0.0318\\
\midrule
\multirow{4}{*}{\makecell{Spambase \\ $(m = 1840)$ }}& 0.8 & 735 & 66 & 749 \\ 
& 0.9 & 864 & 87 & 872 \\ 
& 0.95 & 931 & 104 & 935 \\ 
\cmidrule{2-5}
& \textbf{HHI} & 0.033 & 0.1071 & 0.0325 \\
\midrule
\midrule
\multirow{4}{*}{\makecell{Wdbc \\ $(m = 227)$ }}& 0.8 & 795 & 782 & 789\\ 
& 0.9 & 897 & 890 & 894\\ 
& 0.95 & 948 & 945 & 947 \\ 
\cmidrule{2-5}
& \textbf{HHI} & 0.0316 & 0.0317 & 0.0317\\
\midrule
\midrule
\multirow{4}{*}{\makecell{Bupa \\ $(m = 138)$ }}& 0.8 & 737 & 347 & 668 \\ 
& 0.9 & 858 & 472 & 798 \\ 
& 0.95 & 927 & 561 & 883\\ 
\cmidrule{2-5}
& \textbf{HHI} & 0.0324 & 0.0492 & 0.0342\\ 
\midrule
\midrule
\multirow{4}{*}{\makecell{Banknote \\ $(m = 549)$ }}& 0.8 & 798 & 794 & 793 \\ 
& 0.9 & 899 & 897 & 896 \\ 
& 0.95 & 949 & 948 & 947 \\ 
\cmidrule{2-5}
& \textbf{HHI} & 0.0316 & 0.0317 & 0.0317\\
\bottomrule
\end{tabular}
\end{table}

\begin{figure} \label{appdx_fig:allQChi2.cdf}
\centering
\includegraphics[width = 0.48\textwidth]{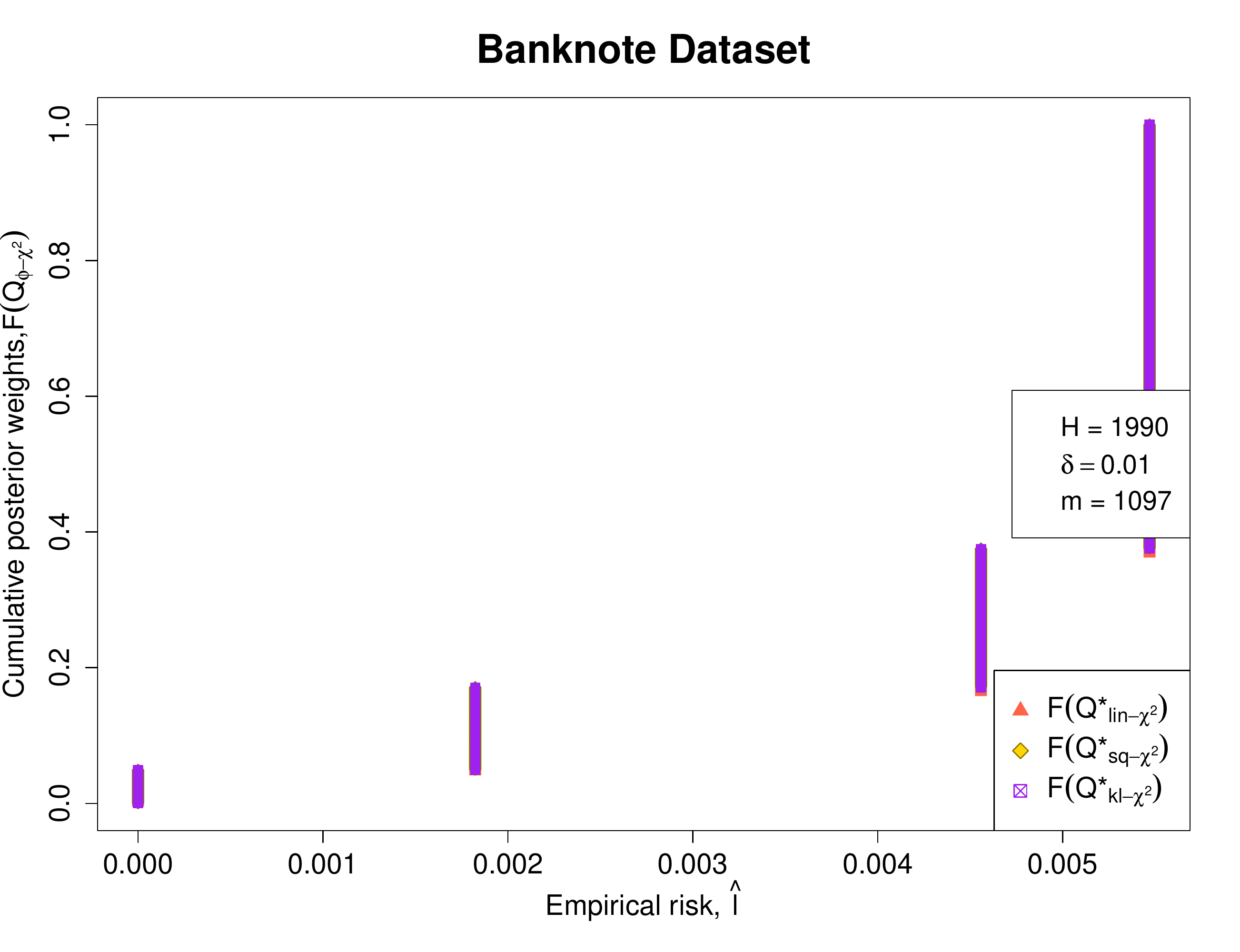}
\includegraphics[width = 0.48\textwidth]{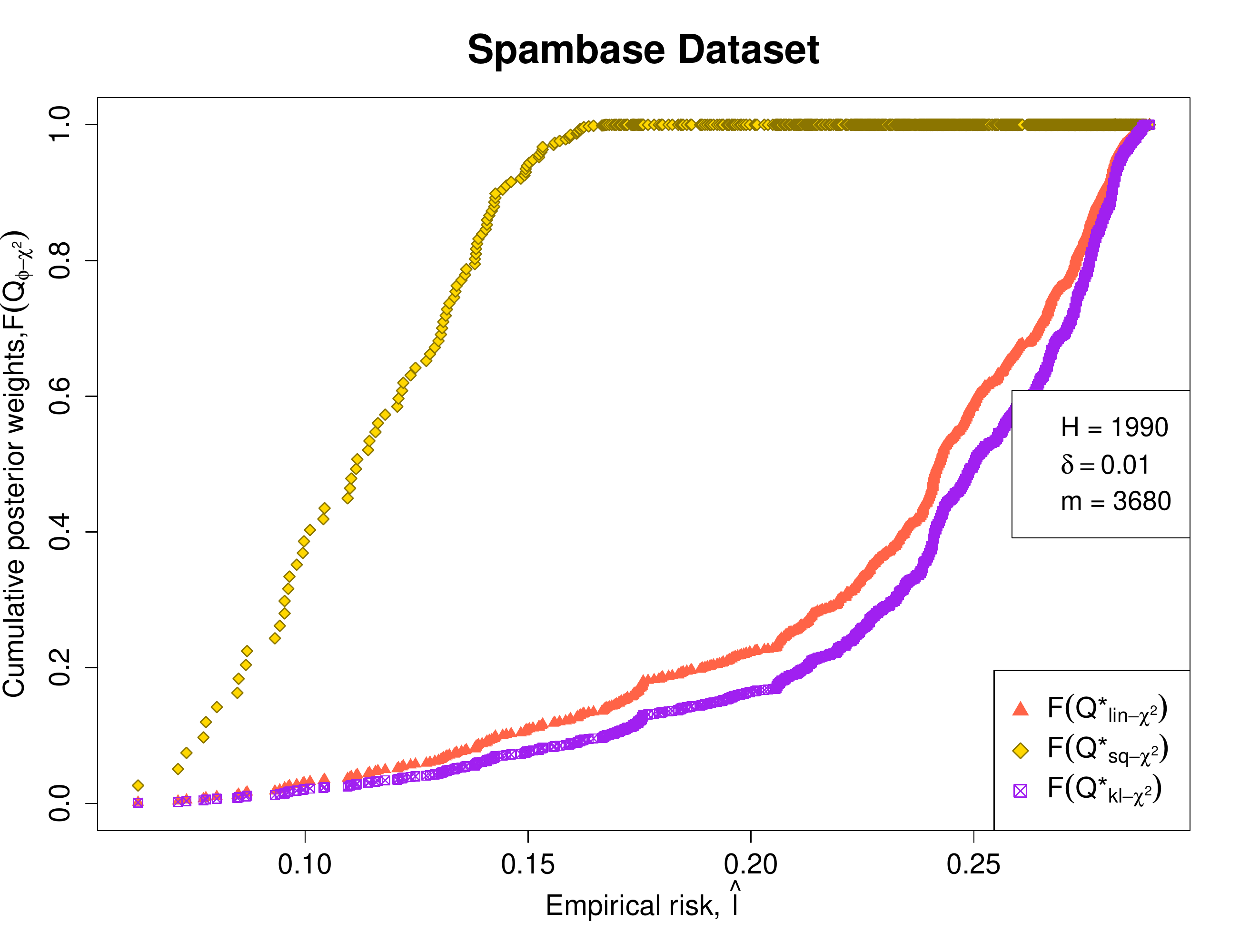}
\includegraphics[width = 0.48\textwidth]{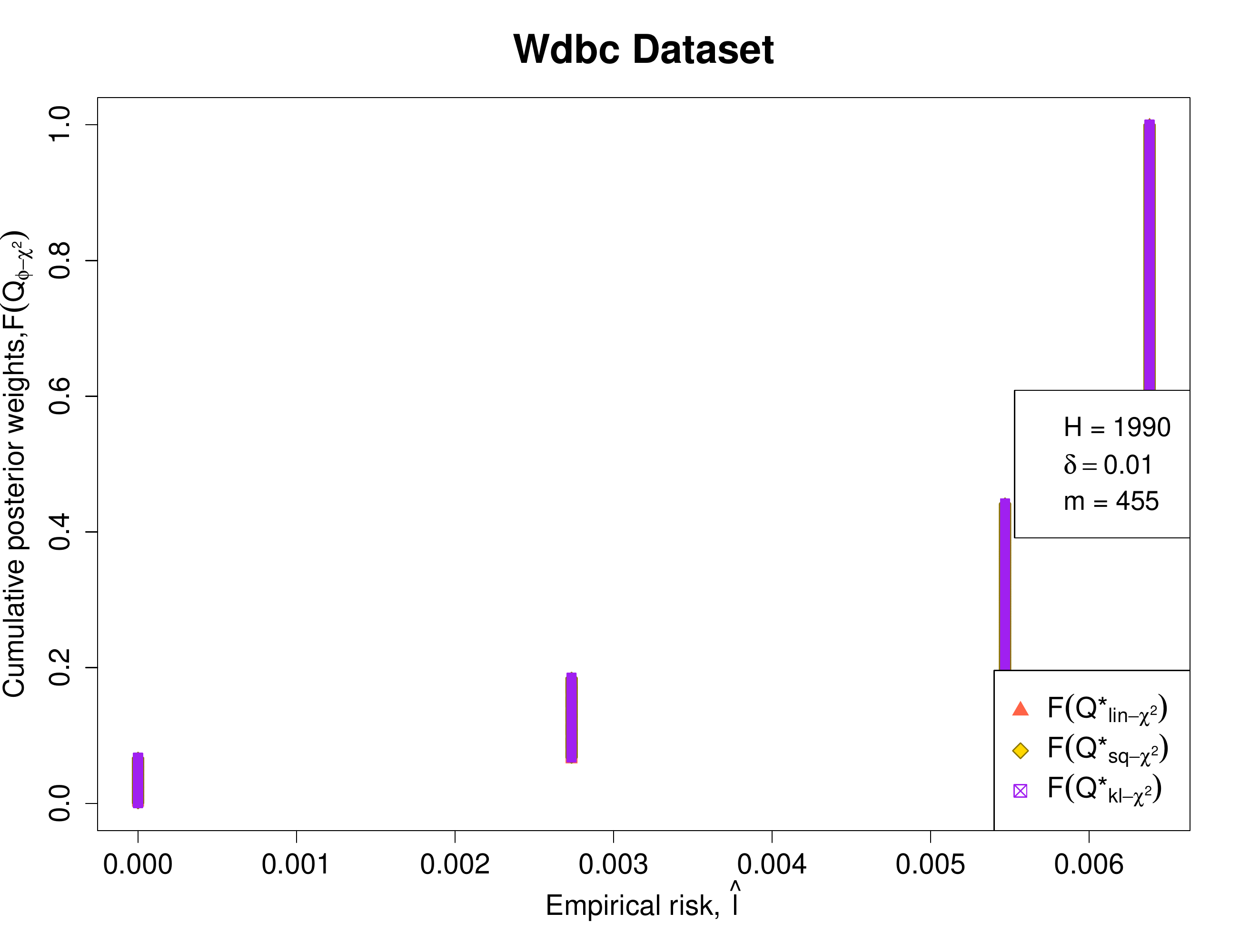}
\includegraphics[width = 0.48\textwidth]{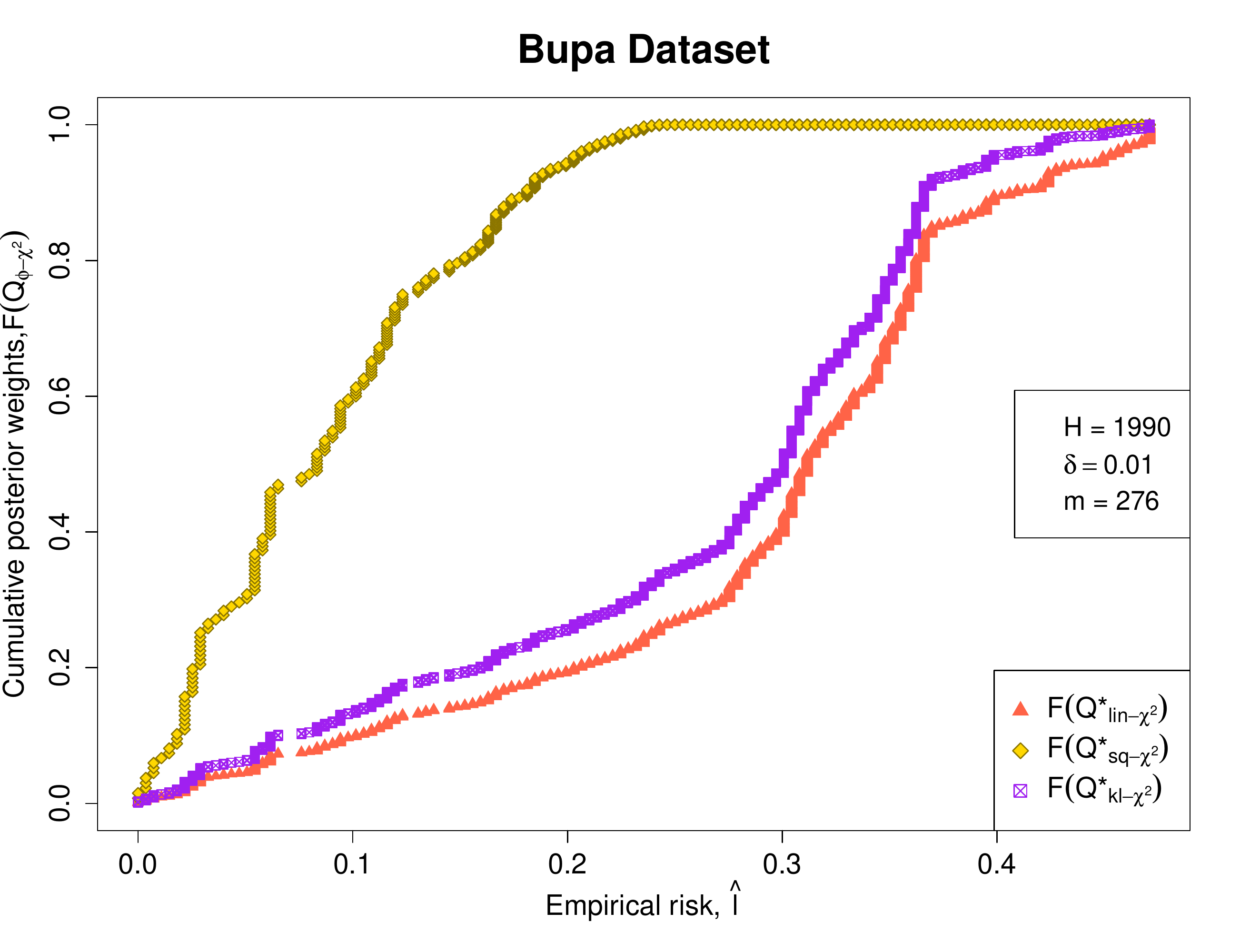}
\includegraphics[width = 0.48\textwidth]{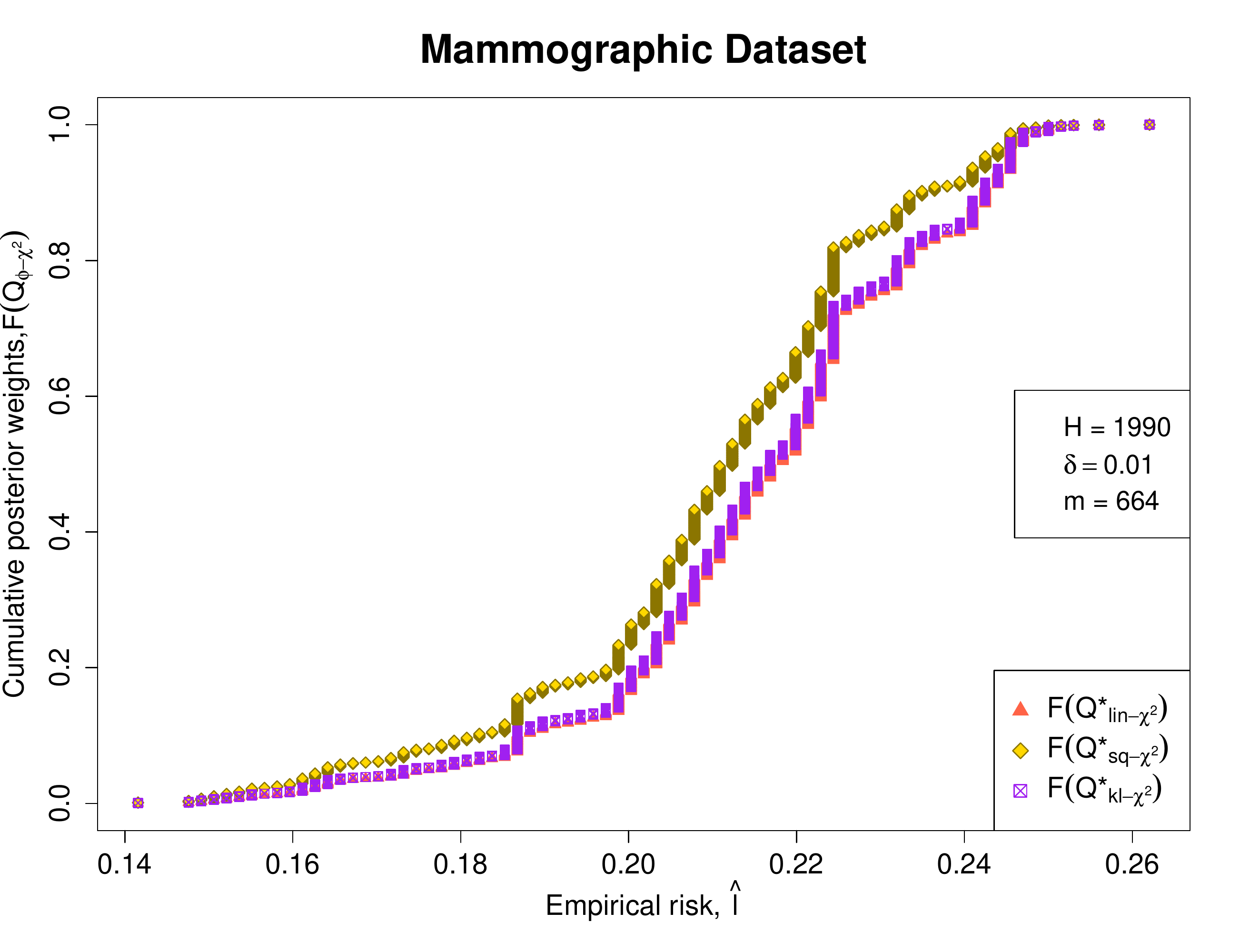}
\caption{Cumulative Distribution Functions (CDFs), $F_{Q^{\ast}_{\phi, \chi^2}}(\cdot)$ for the different PAC-Bayesian optimal posteriors. We observe that, squared distance based posteriors ($Q^{\ast}_{\text{sq}, \chi^2}$) stochastically dominates the other two ($Q^{\ast}_{\text{lin}, \chi^2}$ and $Q^{\ast}_{\text{kl}, \chi^2}$) for Spambase, Bupa and Mammographic datasets. For almost separable datasets (Wdbc and Banknote), the three posteriors have overlapping CDFs. Refer to Table \ref{appdx_tab:allQChi2.HHI} for discussion regarding concentration levels and sparsity. \label{appdx_fig:allQChi2.cdf}}
\end{figure}

\end{document}